\documentclass[11pt]{article}
\usepackage{amssymb,amsmath,amsfonts,amsthm,array,bm,bbm,color}%
\usepackage{mathrsfs}
\usepackage{graphicx}
\providecommand{\U}[1]{\protect\rule{.1in}{.1in}}

\usepackage[
a4paper,
left=25mm,
right=25mm,
top=22mm,
bottom=25mm
]{geometry}

\newcommand{\I}{\hspace{0.2mm}\text{I}\hspace{0.2mm}}
\newcommand{\II}{\text{I \hspace{-2.8mm} I} }
\newcommand{\III}{\text{I \hspace{-2.9mm} I \hspace{-2.9mm} I}}

\numberwithin{equation}{section}

\newtheorem{theorem}{Theorem}[section]
\newtheorem{lemma}[theorem]{Lemma}
\newtheorem{corollary}[theorem]{Corollary}
\newtheorem{proposition}[theorem]{Proposition}
\newtheorem{remark}[theorem]{Remark}
\newtheorem{example}[theorem]{Example}
\newtheorem{definition}[theorem]{Definition}
\newtheorem{assumption}[theorem]{Assumption}

\newcommand\Id{\mathrm{Id}}
\newcommand{\ind}{\mathbf{1}}

\DeclareMathOperator*{\esssup}{ess\,sup}

\def\<{\langle}
\def\>{\rangle}

\def\d{{\rm d}}
\def\d{\, {\rm d}}
\def\i{{\rm i}}

\def\E{\mathbb{E}}
\def\N{\mathbb{N}}
\def\P{\mathbb{P}}
\def\R{\mathbb{R}}

\def\T{\mathbb{T}}

\def\Z{\mathbb{Z}}

\def\cF{\mathcal{F}}

\def\cH{\mathcal{H}}
\def\cX{\mathcal{X}}

\begin{document}

\title{Anomalous regularity and scaling limit by rough transport noise for stochastic 2D fluid models on the torus}
\author{Dejun Luo\footnote{Email: luodj@amss.ac.cn. SKLMS, Academy of Mathematics and Systems Science, Chinese Academy of Sciences, Beijing 100190, China, and School of Mathematical Sciences, University of Chinese Academy of Sciences, Beijing 100049, China}
  \quad Bin Tang\footnote{Email: tangbin@njust.edu.cn. School of Mathematics and Statistics, Nanjing University of Science and Technology, Nanjing, Jiangsu, 210094, China}
}

\maketitle

\vspace{-15pt}

\begin{abstract}
We consider the stochastic 2D Euler and mSQG equations on the torus, with $L^2$ initial data and driven by rough Kraichnan transport noise, for which we construct weak solutions enjoying anomalous regularity and dissipation properties. Moreover, we show that weak solutions in this regularity class are pathwise unique, thus establishing probabilistically strong well-posedness. We further obtain the scaling limit from these equations to the deterministic Navier--Stokes/mSQG equations, with respect to positive Sobolev norms. The proofs rest on the analysis of energy spectrum of the vorticity and on a uniform weighted discrete Poincar\'e inequality on the lattice, the latter being crucial for making the positive norms accessible.
\end{abstract}



\section{Introduction}\label{sec: introduction}

We consider the stochastic two-dimensional (2D) modified Surface Quasi-Geostrophic (mSQG) equations on the torus $\T^2$, driven by rough Kraichnan-type transport noise:
\begin{equation}\label{stoch-2D-Euler}
  \d\omega+ u\cdot \nabla\omega\, \d t+ \circ \d W_t \cdot \nabla\omega=0,
  \quad u= \nabla^\perp (-\Delta)^{-\beta/2}\omega, \quad \beta\in(1,2].
\end{equation}
Here $\omega$ and $u$ are the vorticity and velocity fields; the endpoint $\beta=2$ corresponds to the classical 2D Euler equation. The noise $W=W(t,x)$ is a time-dependent vector field on $\T^2$, Brownian in time, colored and divergence-free in space, see Section \ref{subsec: noise} for more details. The notation $\circ \d$ means the Stratonovich stochastic differential which, according to the Wong-Zakai principle, is a physical choice of random perturbation in fluid equations. Inspired by recent works \cite{BGalMau25, CogMau26, JiaoLuo26} on regularization by noise for 2D fluid models, we show that for $\omega_0\in L^2(\T^2)$, the above equation admits a probabilistically strong and analytically weak solution which is pathwise unique, thanks to the anomalous regularity due to the rough noise. We work on the torus $\T^2$, while the previous papers treat the setting of $\R^2$ where the Kraichnan noise enjoys nice rotation symmetry.

Since the solutions possess anomalous regularity properties which are uniform over a class of rough Kraichnan-type noises, we also obtain the scaling limit, with a quantitative convergence rate, to the deterministic Navier--Stokes/mSQG equation in much stronger Sobolev spaces: 
\[
  L^2\big(\Omega; L^2([0,T]; H^{1-}(\T^2)) \cap L^{\infty}([0,T]; L^2(\T^2)) \big),
\]
thus improving some results in \cite{Flandoli_2024_Quantitative}. 
The proof of uniform anomalous regularity relies on a weighted discrete Poincaré inequality on the lattice, proved in Section \ref{sec: poincare inequality}, which extends the result of \cite[Proposition 3.1]{Rowan_2025_Obukhovcorrsin}.

\subsection{Motivations} \label{subsec-motivations}

The well-known Yudovich theory \cite{Yudovich_1995_Uniqueness} asserts that for any $\omega_0\in (L^1\cap L^\infty)(\R^2)$, the vorticity form of 2D Euler equation 
  \begin{equation}\label{2D-Euler}
  \partial_t \omega + u\cdot \nabla\omega = f, \quad \omega|_{t=0} = \omega_0
  \end{equation}
admits a unique solution $\omega \in L^\infty([0,T]; (L^1\cap L^\infty)(\R^2))$ for a suitable source term $f$. When $\omega_0\in (L^1\cap L^p)(\R^2)$ with $p\in (1,\infty)$, the existence of weak solutions is known since the work of DiPerna and Majda \cite{DiPMaj}; however, the uniqueness of these weak solutions remains open for a long time. In recent years, several nonuniqueness results of weak solutions appeared; for instance, Vishik \cite{Vishik_2018a, Vishik_2018b} constructed a forcing term $f\in L^1([0,T]; (L^1\cap L^p)(\R^2))$ such that \eqref{2D-Euler} with $\omega_0 =0$ has infinitely many nonunique weak solutions, see also \cite{ABCDG24} for a more streamlined presentation of the construction. Vishik's approach was adapted in \cite{ABC22} to provide nonunique Leray-Hopf solutions to the 3D Navier--Stokes equations with a suitable forcing term and zero initial condition. We refer to the introductions of \cite{CogMau26, GalLuo25} for more references on weak solutions to \eqref{2D-Euler}. 

In the last two decades, there have been intensive studies on the theory of regularization by noise, see the lecture notes by Flandoli \cite{Flandoli11} for a comprehensive summary of early works. Nowadays, regularization by noise for finite dimensional equations is quite well understood, see e.g. \cite{Zvonkin74, Ver80, KryRok05, FedFla13, Krylov21, RocZhao25}; some of these results have also been extended to the infinite dimensional setting, cf. \cite{DaPFla10, DaPFPR13, DaPFRV16}, but the nonlinearity treated in these papers are not in the direction of fluid equations. In the celebrated paper \cite{FGP10}, Flandoli et al. proved that if the drift vector field $b$ is spatially H\"older continuous, then a usual $d$-dimensional Brownian motion $W= W(t)$ regularizes the linear transport equation, such that the stochastic equation
  $$\d \rho + b\cdot \nabla \rho\d t + \circ\d W \cdot \nabla \rho =0 $$
admits a unique weak solution. However, as mentioned in \cite[Section 6.2]{FGP10}, such a space-independent noise does not work for regularizing nonlinear fluid equations like 2D Euler, and one needs noises with nontrivial space correlations. 

As an attempt to establishing regularization by noise for fluid models,  Flandoli \cite{Flandoli_2015_Open} proposed to study the regularized stochastic 2D Euler equation
  \begin{equation}\label{eq-regul-2D-Euler}
  \d\omega+ (-\Delta)^{-\epsilon} ({\rm curl}^{-1} \omega) \cdot \nabla\omega\, \d t+ \circ \d W_t \cdot \nabla\omega=0, \quad \omega_0\in L^2(\T^2)
  \end{equation}
with some small $\epsilon>0$; the problem is to find suitable space-time noise $W(t,x)$ such that the above equation enjoys uniqueness in law. A possible idea, as discussed in \cite{Flandoli_2015_Open}, is to transform \eqref{eq-regul-2D-Euler} into a stochastic linear transport equation via the Girsanov theorem, i.e. absorbing the nonlinearity into the noise part. The remaining tasks are to prove well-posedness of the stochastic linear transport equation and then deduce uniqueness in law of solutions to \eqref{eq-regul-2D-Euler}. This has been achieved in \cite{GalLuo25}, see also \cite{BFM10, Bianchi13, BBF14} for early applications of the idea to dyadic model of turbulence and 3D Leray-$\alpha$ model of Euler equations. Although the basic strategy is clear, we would like to mention that several difficulties arise when implementing the above procedure, since the noise $W(t,x)$ has to be relatively rough in space variable in order to absorb the nonlinearity; however, spatial roughness of noise results in difficulties in proving the well-posedness of stochastic linear transport equation. Therefore, one has to find a delicate balance between the regularity properties of noise and the nonlinearity. 

In the recent work \cite{CogMau26}, Coghi and Maurelli made a breakthrough by considering the stochastic 2D Euler equation on $\R^2$ (\eqref{stoch-2D-Euler} with $\beta=2$), driven by Kraichnan transport noise with parameter $\alpha\in (0,1)$ which determines spatial regularity of the noise. For any $\omega_0 \in (L^1\cap L^p)(\R^2)$ with $p> 3/2$, they showed that if $\frac{2}{p\wedge 2} -1 <\alpha< 1- \frac{1}{p\wedge 2} $, then \eqref{stoch-2D-Euler} has analytical weak solutions which are pathwise unique, thus establishing probabilistically strong well-posedness. A key ingredient in their proof is to construct weak solutions enjoying $H^{-\alpha}$ anomalous regularity in $H^{-1}$-estimate of vorticity. However, since rougher noise is expected to have stronger regularizing effects, the lower bound on $\alpha$ is unnatural; this restriction was relaxed by Jiao and Luo \cite{JiaoLuo25b} to cover all $p>1$ and $\alpha< 1- \frac{1}{p\wedge 2}$. Coghi and Maurelli's approach was soon extended in \cite{BGalMau25, JiaoLuo25a} to the stochastic mSQG equations with Kraichnan noise, for suitable choices of the parameters $\alpha$ defining the noise, $\beta$ relating fluid velocity and vorticity, and $p$ concerning integrability of initial data. Later, Jiao and Luo \cite{JiaoLuo26} applied the method in \cite{CogMau26, JiaoLuo25b} to the 2D Boussinesq system with zero viscosity and positive thermal diffusivity, where the coupling between the vorticity and fluid temperature brings difficulties in establishing uniform $L^p$-estimates on vorticity. Recently, Galeati et al. \cite{GGM26} obtained anomalous regularity for Kraichnan's passive scalar model in general $H^{-s}$-estimates including $s\in (0,1)$. Bagnara and Galeati \cite{BagGal26} remove the limitation on $\alpha\in (0,1/2)$ and their results cover the regularizing properties of Kraichnan noise with $\alpha=2/3$, a choice which is more relevant in the context of turbulence modelling. The above works always deal with incompressible Kraichnan noise. Drivas et al. \cite{Drivas_2025_Anomalous} made a thorough investigation of the anomalous regularity properties of passive scalars by Kraichnan noise also in the compressible case; in particular, their results revealed the sharpness of the gain of $1-\alpha$ Sobolev regularity, see also \cite{GGP26} for more recent extensions to the active scalar models. Finally, we point out the paper \cite{BGroMau} which studies the Kazantsev-Kraichnan model for passive magnetic field with transport-stretching noises, thus extending some results of \cite{GGM26} to the 3D case.

The papers mentioned in the last paragraph deal all with the Kraichnan noise on the full space, leaving the interesting problem of extending similar results to the torus setting. In fact, as discussed in \cite[Remark 1.7]{Drivas_2025_Anomalous}, Kraichnan noise on the full space enjoys nice rotational symmetry which is missing on the torus, hence the proofs of several results in \cite{Drivas_2025_Anomalous} do not work for the torus case, where new ideas are needed for establishing anomalous regularity and dissipation. The anomalous diffusion of Kraichnan model on the torus was obtained by Rowan \cite{Rowan_2024_Anomalous} with a PDE-based approach, by studying the averaged two-point correlation function of passive scalars which solves a deterministic linear PDE. Using the techniques of parabolic PDE together with functional inequalities, Rowan obtained the decay of $L^2$-norms by evaluating the two-point correlation function on the diagonal. However, as mentioned in \cite[Section 3.1]{Rowan_2025_Obukhovcorrsin}, this approach does not work for establishing anomalous regularity of passive scalar driven by Kraichnan noise, see in particular the disscussions below (3.2) in \cite{Rowan_2025_Obukhovcorrsin}. 

On the other hand, together with G.H. Zhao, we proposed in \cite{Luo_2024_Elementary} to establish the mixing property of solutions to stochastic transport equations 
  \begin{equation}\label{transp-eq}
  \d \rho + \circ \d W\cdot \nabla \rho =0, \quad \rho_0\in L^2(\T^d) 
  \end{equation} 
by investigating the energy spectrum $X_k(t):= \E|\widehat \rho_k(t)|^2$, $k\in \Z_0^d$, $t\ge 0$. It is known that $\{X_k \}_k$ satisfies an infinite dimensional ODE system, see e.g. \cite[Lemma 2, p.859]{Galeati_2020_Convergence}. A key observation in \cite{Luo_2024_Elementary} is that the negative Sobolev norms of $\rho$ can be estimated by suitable $\ell^p(\Z_0^d)$-norms of the energy spectrum, e.g.
  $$\E \|\rho(t) \|_{H^{-1}}^2 = \sum_k |k|^{-2} X_k(t) \le \|X(t) \|_{\ell^p} \bigg(\sum_k |k|^{-\frac{2p}{p-1}} \bigg)^{\frac{p-1}{p}} ,$$
where the last series is finite for suitable $p>1$. Therefore, to obtain decay rate of $\E \|\rho(t) \|_{H^{-1}}^2$, one can turn to studying the time evolution of $\|X(t) \|_{\ell^p}$. To this end, we showed that 
  $$\aligned
    \frac{\d}{\d t} \sum_k X_k^p = - 8 \pi^2 p \nu \sum_{k,l}\theta_l^2 |\Pi_{l^\perp} k|^2 \big(X_{k+l}^{p-1} -X_k^{p-1}\big) (X_{k+l} -X_k) ,
  \endaligned $$
where $\nu>0$ and $\{\theta_l \}_{l\in \Z_0^d}$ are the intensity and coefficients of noise, and $\Pi_{l^\perp} k$ is the projection of $k$ on $l^\perp= \{y\in \R^d: l\cdot y=0 \}$. We proved a discrete Poincar\'e inequality by combining its one-dimensional version (see \cite[Lemma 2.6]{Luo_2024_Elementary}) with a decomposition of $\Z_0^d$ into one-dimensional paths, and deduced from it the exponential decay of $\|X(t) \|_{\ell^p}$ and hence the mixing of $\rho$. 

Anomalous regularity, however, requires a discrete Poincar\'e inequality in weighted form, with the weight measuring the coercivity of the noise at each frequency. Partly motivated by our work, Rowan revisited in \cite{Rowan_2025_Obukhovcorrsin} the Kraichnan passive scalar model and successfully proved a weighted discrete Poincar\'e inequality (see \cite[Proposition 3.1]{Rowan_2025_Obukhovcorrsin} and Proposition \ref{prop: weighted Poincare} below for a slightly different version), which allowed him to establish anomalous regularity of solutions to stochastic transport-diffusion equation, and to provide a proof of the so-called Obukhov--Corrsin spectrum prediction. In the current paper, we extend the methods of \cite{Luo_2024_Elementary, Rowan_2025_Obukhovcorrsin} to the active scalar model \eqref{stoch-2D-Euler} with Kraichnan-type noise, proving anomalous regularity and dissipation of analytical weak solutions, as well as their strong well-posedness in the probabiliatical sense. Our extension gives a weighted lower bound uniform over the whole noise class $\cH_{\alpha}$ (see Section \ref{subsec: noise} for definition), hence uniform anomalous regularity. This uniformity enables us to strengthen the previous scaling limit results in \cite{Flandoli_2024_Quantitative}, showing convergence in $L^2([0,T]; H^{1-\alpha-}) \cap L^{\infty}([0,T]; L^2)$ of solutions to those of the deterministic 2D Navier--Stokes/mSQG equations.

\subsection{Main results} \label{subsec-main-results}

In this section, we present the main results obtained in our paper. Using the notation of noise explained in Section \ref{subsec: noise}, we can write equation \eqref{stoch-2D-Euler} more precisely  as follows:
  \begin{equation*}
  \d\omega+ u\cdot \nabla\omega \d t + \sqrt{2\nu} \sum_k \theta_k \sigma_k \cdot \nabla\omega \circ \d B^k_t=0,
  \end{equation*}
where $u= \nabla^\perp (-\Delta)^{-\beta/2}\omega$ with $\beta\in (1,2]$, thus we treat the stochastic 2D Euler and mSQG equations at the same time. Thanks to the choices of noise coefficients $\{\theta_k \}_k$ and divergence-free fields $\{\sigma_k \}_k$, the above equation admits the It\^o form 
  \begin{equation}\label{stoch-2D-Euler-1}
  \d\omega+ u\cdot \nabla\omega \d t + \sqrt{2\nu} \sum_k \theta_k \sigma_k \cdot \nabla\omega \d  B^k_t= \nu \Delta \omega \d t, 
  \end{equation}
where $\nu \Delta$ comes from the It\^o-Stratonovich corrector. It is well known that, in formal $L^2$-estimate of vorticity, the corrector is cancelled by the quadratic variation of martingale part, thus the above equation is of hyperbolic type. A probabilistical and analytical weak solution to \eqref{stoch-2D-Euler-1} is defined as usual in the distributional sense, see Definition \ref{def-stoch-Euler/mSQG} for the exact meaning. Our purpose is to show that, when the noise is of Kraichnan type and rough in space, with coefficients in the class $\cH_\alpha$ introduced in Section \ref{subsec: noise}, weak solutions of \eqref{stoch-2D-Euler-1} enjoy an $H^{1-\alpha-}$ anomalous regularity, which in turn implies that they are pathwise unique.

Here is the precise statement of first result, see Section \ref{subsec: setting} for the definition of space $\mathcal{X}_{\alpha,p}$. 

\begin{theorem}[Existence and pathwise uniqueness] \label{thm-Euler/mSQG}
Assume $\beta\in (1,2]$, $\alpha\in (0, \frac{\beta-1}2)$ and $\omega_0 \in L^2(\T^2)$. For every fixed noise coefficients $\{\theta_k \}_k \in \cH_{\alpha}$, the stochastic 2D Euler/mSQG equation \eqref{stoch-2D-Euler-1} admits a weak solution $\omega$ in the sense of Definition \ref{def-stoch-Euler/mSQG}, satisfying the anomalous regularity: for any $p \in (1,\frac{\beta-2\alpha+1}{2})$ and $\alpha^{\prime}\in (\alpha, 1)$,
  \begin{equation}\label{eq-anomal-regular}
    \int_{0}^T \| \omega_t \|_{\mathcal{X}_{\alpha,p}}^{2p} \, \d t<+\infty, \quad 
    \E \int_0^T \|\omega_t \|_{H^{1-\alpha^{\prime}}}^2\, \d t< \infty,
  \end{equation}
and the map $t \mapsto \|\omega_t\|_{L^2}$ is $\P$-a.s. nonincreasing. Furthermore, if 
  \begin{equation}\label{eq-param-alpha-p}
  0<\alpha<\frac{\beta+1-\sqrt{(\beta-1)^2+4}}{4},\quad \frac{1-\alpha}{1-2\alpha} <p<\frac{\beta-2\alpha+1}{2},
  \end{equation}
then equation \eqref{stoch-2D-Euler-1} enjoys pathwise uniqueness of weak solutions in the space
  $$L^{2p}([0,T];\mathcal{X}_{\alpha,p}) \cap L^{\infty}\big(\Omega;L^{\infty}([0,T];L^2(\T^2)) \big). $$
\end{theorem}

We briefly describe the strategy of proof which is based on vanishing viscosity arguments. We introduce a small viscosity term $\kappa \Delta$ in \eqref{stoch-2D-Euler-1}:
  \begin{equation}\label{stoch-2D-NS-0}
  \d\omega+ u\cdot \nabla\omega \d t + \sqrt{2\nu} \sum_k \theta_k \sigma_k \cdot \nabla\omega \d  B^k_t= (\kappa+ \nu) \Delta \omega \d t. 
  \end{equation}
Setting $Y_k(t)= \E|\widehat \omega_k(t)|^2$, we follow \cite{Luo_2024_Elementary} to study the time evolution of the energy spectrum in Section \ref{sec: anomalous regularity}, where we also derive an identity for the $\ell^p$-norm of $\{Y_k(t) \}_k$ (Proposition \ref{prop-Y_k^p}). Compared with the linear transport equation \eqref{transp-eq} discussed above, the main difficulty now lies in the estimate of the nonlinear term, which is carried out in Proposition \ref{prop-nonlinearity} by a careful analysis of the Fourier expansion of $u\cdot \nabla \omega$. Combining these estimates with the discrete Poincar\'e inequality (cf. Proposition \ref{prop: weighted Poincare}), we show that the solutions to \eqref{stoch-2D-NS-0} enjoy anomalous regularity and dissipation which are independent of $\kappa>0$ and uniform over $\cH_\alpha$. A classical vanishing viscosity argument then yields weak solutions to \eqref{stoch-2D-Euler-1} with the desired regularity. Pathwise uniqueness is proved by a similar procedure, namely by analysing the energy spectrum of the difference of two weak solutions corresponding to the same initial value and the same noise; here the derivation of the identity for the $\ell^p$-norm of the spectrum requires both the anomalous regularity established above and the exchangeability of summation order proved in Lemma \ref{lem-sum-cutoff}, since the difference lacks the $H^1$ regularity available in the viscous case; see Section \ref{sec: pathwise uniqueness} for details. 

We remark that the $H^{1-\alpha-}$ anomalous regularity of Theorem \ref{thm-Euler/mSQG} is optimal in the sense that the borderline exponent $1-\alpha$ cannot be attained. For small initial data, Proposition \ref{prop-optimal} shows that the weak solutions constructed by the vanishing viscosity method satisfy $\E \int_0^T \|\omega_s \|_{H^{1-\alpha}}^2 \d s=+\infty$ for all $T$ sufficiently large. The proof is by contradiction: if that integral were finite, then the $L^2$-norm of solution would be exactly conserved, in contradiction with the exponential decay established in Theorem \ref{thm: anomalous regularity} for small initial data. We refer to \cite{Drivas_2025_Anomalous, GGP26} for more discussions on the sharpness of anomalous regularity due to Kraichnan-type noises.


\begin{remark}
We mention that here we only consider $L^2$-initial data $\omega_0$, while in the previous papers (see e.g. \cite{BGalMau25, CogMau26, JiaoLuo26}) working on the full space, more general initial data belonging to $L^p$-spaces with $p\neq 2$ are allowed. Note that $p>1$ is used in our paper as an auxiliary parameter for summability of Fourier modes, cf. the discussions below \eqref{transp-eq}. It is possible to extend our results to more general initial data $\omega_0 \in L^q$ for some $q>1$, but this will greatly increase technical complexties of proofs. Therefore, we decide to treat only $L^2$-initial data in this paper. 
\end{remark}

As an application of the uniform anomalous regularity proved in Theorem \ref{thm-Euler/mSQG}, we provide a result which strengthens the scaling limits obtained in \cite{Flandoli_2024_Quantitative}. To this end, consider the following stochastic 2D Euler/mSQG equations
\begin{equation} \label{2D Euler/mSQG-n}
  \d \omega^n + u^n \cdot\nabla\omega^n \d t + \sqrt{2\nu}\sum_k\theta_{n,k}\sigma_k\cdot\nabla\omega^n \d B_t^k = \nu \Delta\omega^n \d t, \quad u^n = \nabla^\perp(-\Delta)^{-\beta/2}\omega^n,
\end{equation}
with a sequence of noise coefficients $\{\theta_{n,k}\}_{k\in\mathbb Z_0^2}$ defined as in \eqref{def: theta-n}. By Lemma \ref{lem: power tail}, one has
\begin{equation*}
  \|\theta_n \|_{\ell^\infty} :=\sup_{k\in\mathbb Z_0^2}|\theta_{n,k}| = Z_{n,\alpha} n^{-1-\alpha} \sim \sqrt{\alpha}\, n^{-1} \to 0, \quad n\to\infty.
\end{equation*}
Flandoli et al. \cite{Flandoli_2024_Quantitative} proved quantitative bounds, with respect to the norm of $C([0,T]; H^{-\epsilon})$, between the solutions to \eqref{2D Euler/mSQG-n} and that of the deterministic 2D Navier--Stokes/mSQG equation:
\begin{equation} \label{2D determined NSE/mSQG}
  \partial_t \widetilde{\omega} + \widetilde{u} \cdot\nabla \widetilde{\omega} = \nu \Delta\widetilde{\omega}, \quad \widetilde{u} = \nabla^\perp(-\Delta)^{-\beta/2} \widetilde{\omega}.
\end{equation}
Thanks to the anomalous regularity of solutions to \eqref{stoch-2D-Euler-1}, we can upgrade the limit results to stronger spaces $L^2([0,T], H^{1-\alpha-}) \cap L^{\infty}([0,T]; L^2)$. 

\begin{theorem}[Scaling limits] \label{thm: scal_Sobolev}
Fix $\nu > 0$, $\beta \in (1,2]$ and $\alpha\in \big(0, \frac{1}{4} (\beta+1-\sqrt{(\beta-1)^2+4}\, ) \big)$. Let $\omega^n$ and $\widetilde{\omega}$ be respectively the unique solution of \eqref{2D Euler/mSQG-n} and \eqref{2D determined NSE/mSQG} with initial data $\omega_0 \in L^2(\T^2)$. 
\begin{itemize}
\item[\rm (a)] For any $T >0$, $\alpha^{\prime} \in (\alpha,1]$, $0<s_0<1-\frac{1}{\beta}$ and any subinterval $I \subset [0,T]$, there exists a constant $C=C(\alpha,\alpha^{\prime},s_0,\beta,\nu,T)>0$ such that
  \begin{equation} \label{eq: scal_sobolev}
    \E \int_I \|\omega_t^n - \widetilde{\omega}_t\|_{H^{1-\alpha^{\prime}}}^2 \, \d t
    \le C \exp \left( C \|\omega_0\|_{L^2}^{\frac{4(1-\alpha)}{\beta-2\alpha-1}} \right)  \|\omega_0\|_{L^2}^2  \frac{|I|^{(\frac{3}{8}+\frac{s_0}{16})(\alpha^{\prime}-\alpha)}}{n^{\frac{3}{2} s_0 (\alpha^{\prime}-\alpha)}}.
  \end{equation}
\item[\rm (b)] For any $T > 0$, it holds
  \begin{equation*}
    \lim_{n \to + \infty} \E \Big[ \sup_{t \in [0, T]} \|\omega_t^n - \widetilde{\omega}_t\|_{L^2}^2 \Big]=0.
  \end{equation*}
\end{itemize}
\end{theorem}

Note that claim (a) provides us with a quantitative convergence rate, while (b) is only a qualitative statement. Theorem \ref{thm-quanti-C-L-2} provides a quantitative version of (b), but it is not uniform in the $L^2$-norm of $\omega_0$; indeed, convergence in $L^\infty([0,T];L^2)$ cannot be uniform in $\|\omega_0\|_{L^2}$, as shown in Proposition \ref{prop-nonuniform-C-L-2}. We also mention that, as the constant in the weighted discrete Poincar\'e inequality \eqref{def-C-theta} does not depend on $\alpha$, one may go further and choose accordingly $\alpha= \alpha_n\to0$; in this way, we directly obtain the scaling limit in $H^{1-}$. Such an adaptation is routine thanks to the uniformity of anomalous regularity established above; however, we do not pursue this here for convenience of presentation.

We finish this section with the organization of the paper. In Section \ref{sec: preliminaries}, we make some preparations by introducing notation and functional spaces, together with the exact choice of the space-time noise $W(t,x)$. Then we study in Section \ref{sec: anomalous regularity} the stochastic Navier--Stokes/mSQG equation \eqref{stoch-2D-NS-0}, trying to show anomalous regularity and dissipation which are uniform in $\kappa>0$, see Theorem \ref{thm: anomalous regularity}. Sections \ref{sec: pathwise uniqueness} and \ref{sec: scaling limit} are devoted to the proofs of Theorems \ref{thm-Euler/mSQG} and \ref{thm: scal_Sobolev}, respectively. Finally, in the appendices, we collect in Section \ref{append-tools} some analytical results used in the paper, prove in Section \ref{sec: poincare inequality} a slightly improved version of the weighted discrete Poincar\'e inequality, and present in Section \ref{sec: summatation reordering} a technical result on the exchange of summation order.

\section{Preliminaries}\label{sec: preliminaries}

This section contains some preparatory work. We introduce in Section \ref{subsec: setting} some notations and conventions used throughout the paper, and describe in Section \ref{subsec: noise} the precise choice of the noise $W(t,x)$ and some properties on it. 

\subsection{Notions and functional setting} \label{subsec: setting}

Let $\T^d=\R^d/\Z^d$ be the $d$--dimensional torus and $\Z_0^d = \Z^d \setminus \{0\}$ be the set of nonzero lattice points. Let $\Pi_{l^\perp} k$ denote the projection of $k$ on $l^\perp= \{y\in \R^d: l\cdot y=0 \}$. For $d=2$ and $l=(l_1,l_2)$, we also write $l^{\perp}:=(-l_2,l_1)$. In the main part of the paper, we restrict to $d=2$, while the higher-dimensional cases are treated only in the appendix. Given positive numbers $a,b$, we write $a \lesssim b$, if there exists constant $C>0$ such that $a \leq C b$, $a \sim b$ if $a \lesssim b$ and $b \lesssim a$; we write $a \lesssim_{\gamma} b$ to stress that the constant $C=C(\gamma)$ depends on some parameter $\gamma$. As we assume the noise is spatially divergence-free (see Section \ref{subsec: noise}), the spatial average of all solutions considered in this paper are preserved. Hence, we shall assume for simplicity that the function spaces in this paper consist of functions on $\T^2$ with zero average.

Let $e_k(x):=\exp \{ 2 \pi {\rm i} k \cdot x\}$ for any $k \in \Z_0^2$ and $x \in \T^2$. Then $\{e_k \}_k$ is the usual complex basis of $L^2(\T^2; \mathbb C)$. Let $\widehat{\omega}_{k}=\langle \omega,e_{-k}\rangle$ be the Fourier coefficient of $\omega$ and $\overline{\widehat{\omega}_k}$ its complex conjugate. For any $a \in \R$, we define the Sobolev space $H^{a}(\T^2)$ as
\begin{equation*}
    H^{a}(\T^2)=\Big\{ \omega = \sum_{k \in \Z_0^2} \widehat{\omega}_{k} \, e_{-k} ;\ 
    \|\omega \|_{H^a}^2:= \sum_{k \in \Z_0^2} |2\pi k|^{2a} 
    |\widehat{\omega}_k|^2 < +\infty\, \mbox{ and } \,
    \overline{\widehat{\omega}_k} = \widehat{\omega}_{-k} \Big\}.
\end{equation*}
Following standard notation, we denote the space $H^{0}(\T^2)$ by $L^2(\T^2)$ with norm $\|\omega\|_{L^2}=\|\omega\|_{H^0}$ and inner product $\langle \omega, \phi \rangle=\sum_{k \in \Z_0^2} \widehat{\omega}_{k}\, \widehat{\phi}_{-k} $. Let $C([0,T]; L^2_w(\T^2))$ be the space of functions $f: [0,T] \rightarrow L^2(\T^2)$ which are continuous with respect to the weak topology of $L^2(\T^2)$. For any $\sigma \geq 0$ and $q \geq 2$, we define the Fourier--Lebesgue space consisting of mean-zero functions taking values in $L^2(\Omega)$ by
\begin{equation*}
  \cF L^{\sigma,q}(\T^2;L^2(\Omega)) := \left\{ f \in L^2(\T^2;L^2(\Omega)) :\, 
  \P \mbox{-a.s.}\, \int_{\T^2} f(x) \, \d x=0, \;  
  \|f\|_{\cF L^{\sigma,q}} < \infty \right\},
\end{equation*}
where
\begin{equation*}
  \|f\|_{\cF L^{\sigma,q}} := \big\| |k|^\sigma \|\widehat{f}(k)\|_{L^2(\Omega)} 
  \big\|_{\ell^q(\Z_0^2)} = \bigg( \sum_{k \in \Z_0^2} |k|^{\sigma q}  
  \big( \E|\widehat{f}(k)|^2 \big)^{q/2} \bigg)^{1/q}.
\end{equation*}
Note that $\|f\|_{\cF L^{\sigma,2}}= \big(\E \|f\|_{H^\sigma}^2 \big)^{1/2}$. For convenience, we denote 
$ \cX_{\alpha,p}:=\cF L^{\frac{1-\alpha}{p},2p}(\T^2;L^2(\Omega))$ for 
$\alpha \in [0,1]$ and $p \geq 1$; then
\begin{equation*}
  \|f\|_{\cX_{\alpha,p}}
  :=\|f\|_{\cF L^{\frac{1-\alpha}{p},2p}(\T^2;L^2(\Omega))}
  =\bigg( \sum_{k \in \Z_0^2} |k|^{2(1-\alpha)}  
  \big( \E|\widehat{f}(k)|^2 \big)^{p} \bigg)^{\frac{1}{2p}}.
\end{equation*}
The relations between $\cX_{\alpha,p}$ and Sobolev space $H^{1-\alpha}$ are summarized in the following proposition.

\begin{proposition} \label{prop:X-Sobolev}
Let $p>1$, $\alpha \in (0,1)$,
  \begin{enumerate}
    \item[(1)] If $s\geq (1-\alpha)/p$, then
    \[
    L^2(\Omega; H^s) \hookrightarrow \mathcal X_{\alpha,p},
    \quad
    \|u\|_{\mathcal X_{\alpha,p}} \leq \|u\|_{L^2(\Omega; H^s)}.
    \]
    \item[(2)] For every $\alpha^{\prime}>2-\frac{2-\alpha}{p}$, then 
    \[
    \mathcal X_{\alpha,p} \hookrightarrow L^2(\Omega; H^{1-\alpha^{\prime}}), \quad
    \|u\|_{L^2(\Omega; H^{1-\alpha^{\prime}})}  \lesssim_{\alpha,\alpha^{\prime},p}\|u\|_{\mathcal X_{\alpha,p}}.
    \]
  \end{enumerate}
\end{proposition}

\begin{proof}
If $s\geq(1-\alpha)/p$, then the first assertion is given by
  \[
  \begin{aligned}
  \|u\|_{\mathcal X_{\alpha,p}}^{2p}
  & = \sum_{k\neq0} |k|^{2(1-\alpha)-2sp}
  \bigl(|k|^{2s}\mathbb E|\widehat u_k|^2\bigr)^p \\
  & \leq \sum_{k\neq0}\bigl(|k|^{2s}\mathbb E|\widehat u_k|^2\bigr)^p
  \leq \bigg(\sum_{k\neq0}|k|^{2s}\mathbb E|\widehat u_k|^2\bigg)^p.
  \end{aligned}
  \]

Next, for $\alpha^{\prime}>2-\frac{2-\alpha}{p}$, Hölder's inequality gives
  \[
  \begin{aligned}
  \mathbb E\|u\|_{H^{1-\alpha^{\prime}}}^2
  & = \sum_{k\neq0} \left(|k|^{2(1-\alpha)}Y_k^p\right)^{1/p}
  |k|^{2(1-\alpha^{\prime})-\frac{2(1-\alpha)}p} \\
  & \leq \|u\|_{\mathcal X_{\alpha,p}}^2 \bigg( \sum_{k\neq0}
  |k|^{\frac{2p}{p-1} \left( 1-\alpha^{\prime}- \frac{1-\alpha}{p} \right)} \bigg)^{\frac{p-1}{p}}.
  \end{aligned}
  \]
  The lattice sum is finite precisely when
  $\alpha^{\prime}>2-(2-\alpha)/p$, proving the second assertion.
\end{proof}

\subsection{Choice of noise} \label{subsec: noise}

In the main body, we restrict to $d=2$ and consider the following space-time noise:
\begin{equation} \label{def:noise}
W(t,x)= \sqrt{2\nu} \sum_{k\in \Z_0^2} \theta_k \sigma_k(x) B^k_t,
\end{equation}
where $\nu>0$ is the noise intensity, $\theta= \{\theta_k\}_k \in \ell^2(\Z_0^2)$ is the noise coefficient with $\|\theta \|_{\ell^2}=1$, and $\{\sigma_k\}_k$ is the classical complex basis of square-integrable, divergence-free vector fields on $\T^2$. More precisely, let $\Z_0^2= \Z_+^2 \cup \Z_-^2$ be a decomposition of $\Z_0^2$ with $\Z_+^2\cap \Z_-^2= \emptyset$ and $\Z_+^2= - \Z_-^2$, and define $a_k= \frac{k^\perp}{|k|}$ for $k\in \Z_+^2$ and $a_k= a_{-k}$ for $k\in \Z_-^2$. Then
\begin{equation*}
\sigma_k(x)= a_k e_k(x), \quad x\in \T^2, \, k\in \Z_0^2.
\end{equation*}
Finally, $\{B^k\}_k$ is a family of standard complex Brownian motions satisfying
\begin{equation}\label{covar-BM}
  [B^k, B^l]_t= 2\, \delta_{k,-l}\, t, \quad k,l\in \Z_0^2.
\end{equation}
The noise coefficients used in the main body are chosen as follows. For any fixed $\alpha \in (0,1/2)$ and $n \in \N$, define
\begin{equation}\label{def: theta-n}
  \theta_{n,k} :=Z_{n,\alpha}|k|^{-1-\alpha}\mathbf1_{\{|k|\ge n\}},  \quad k\in\mathbb Z_0^2,
\end{equation}
where $Z_{n,\alpha} := \big( \sum_{|k|\ge n} |k|^{-2-2\alpha} \big)^{-1/2}$ is a normalizing constant, and let
\begin{equation*}
\cH_{\alpha}:=\big\{\theta_n= \{\theta_{n,k}\}_k \in \ell^2(\Z_0^2) : n \in \N \big\}
\end{equation*}
be the corresponding set of noise coefficients. Since $\theta_{n,k}$ depends only on $|k|$ and $\sum_k \theta_{n,k}^2=1$, we have the covariance identity
\begin{equation}\label{covar-identity}
  \sum_k \theta_{n,k}^2\, \sigma_k(x)\otimes \sigma_{-k}(x)= \frac12 I_2,
\end{equation}
where $I_2$ is the $2\times 2$ unit matrix, see e.g. \cite[Section 5.1]{Fang_2018_Constantin} or \cite[Section 2.3]{Galeati_2020_Convergence} for a proof. For every $\theta \in \cH_{\alpha}$, Proposition \ref{prop: weighted Poincare} and Example \ref{exm: power law} imply that there exists a uniform constant $C_0>0$, independent of $\theta\in\cH_{\alpha}$, such that for every $p>1$ and every nonnegative sequence $\{ a_{k}\}_{k}$ satisfying $\sum_{k \in \Z_0^2} |k|^{2(1-\alpha)} a_{k}^p <+\infty$, the following discrete Poincar\'e inequality holds:
\begin{equation} \label{def-C-theta}
  8 \pi^2 p \sum_{k,l} \theta_{l}^2  |\Pi_{l^\perp} k|^2  \big(a_{k+l}^{p-1} -a_k^{p-1}\big) (a_{k+l} -a_k) \geq C_0 \frac{p-1}{p}  \sum_{k} |k|^{2(1-\alpha)} a_{k}^p.
\end{equation}
We stress that, by Lemma \ref{lem: power tail}, the constant $C_0$ is independent of $\alpha>0$.

For comparison, the general $d$-dimensional divergence-free noise takes the form (cf. \cite{Galeati_2020_Convergence, FlanLuo21})
\begin{equation*}
  W(t, x)= \sqrt{C_{d} \nu}\, \sum_{k \in \Z_0^d} \sum_{i=1}^{d-1} \theta_{k} \sigma_{k, i}(x) B^{k, i}_{t},
\end{equation*}
where $C_d= \frac{d}{d-1}$, $\sigma_{k,i}(x)= a_{k,i} e_k(x)$, $\{a_{k,i}\}_{i=1}^{d-1}$ being an orthonormal basis of $k^\perp= \{ y\in \R^d: k\cdot y=0\}$, and the noise coefficients $\theta= \{\theta_k\}_k$ are specified in Appendix \ref{sec: poincare inequality}.

\section{Anomalous regularization} \label{sec: anomalous regularity}

In this section, we consider the stochastic 2D Navier-Stokes/mSQG equations with Stratonovich transport noise:
\begin{equation*}
  \d\omega+ u\cdot \nabla\omega\, \d t+ \circ \d W_t \cdot \nabla\omega= \kappa \Delta\omega\, \d t, \quad \omega_0\in L^2(\T^2),
\end{equation*}
where $u$ and $\omega$ denote the velocity and vorticity of the fluid, respectively, and are related by $u= \nabla^\perp (-\Delta)^{-\beta/2} \omega$ for some $\beta\in (1,2]$. Here $\kappa>0$ is the molecular diffusivity, and the space-time noise $W(t,x)$ is defined in Section \ref{subsec: noise} with a noise coefficient $\theta \in \cH_{\alpha}$. As usual, we assume that $\omega_0$ has zero average on $\T^2$, a property preserved by the dynamics.

Based on the above preparations (in particular \eqref{covar-identity}), it can be checked that the above equation admits the It\^o form \eqref{stoch-2D-NS-0}, which is repeated here for the reader's convenience:
\begin{equation}\label{SNSE-Ito}
  \d\omega+ u\cdot \nabla\omega\, \d t+ \sqrt{2 \nu} \sum_k \theta_k \sigma_k \cdot \nabla\omega\, \d  B^k_t= (\kappa + \nu) \Delta\omega\, \d t, \quad \omega_0\in L^2(\T^2).
\end{equation}
A weak solution to the above equation is understood in the classical sense of weak formulation. 

\begin{definition}\label{def-stoch-NS/mSQG}
  We say that \eqref{SNSE-Ito} has a weak solution if there exist a filtered probability space $(\Omega,\mathcal F, \mathcal F_t, \P)$, a sequence of $\mathcal F_t$-complex Brownian motions $\{B^k\}_k$ satisfying \eqref{covar-BM}, and a $\mathcal F_t$-progressively measure process $\omega$ with trajectories in $C([0,T], L^2) \cap L^2([0,T], H^1)$ such that for any $\phi\in C^\infty(\T^2)$, $\P$-a.s. for all $t\in [0,T]$, it holds
  \begin{equation}\label{SNSE-Ito-weak}
    \aligned
    \<\omega_t, \phi\>
    &= \<\omega_0, \phi\>+ \int_0^t\<\omega_s, u_s\cdot \nabla \phi\>\, \d s + (\kappa + \nu)\int_0^t \<\omega_s, \Delta\phi\>\, \d s \\
    &\quad + \sqrt{2 \nu} \sum_k \theta_k \int_0^t\<\omega_s, \sigma_k \cdot \nabla\phi\> \, \d  B^k_s,
    \endaligned
  \end{equation}
  where $u_s= \nabla^\perp (-\Delta)^{-\beta/2} \omega_s$, $s\in [0,T]$.
\end{definition}

Due to the presence of the molecular diffusivity $\kappa>0$, we have the following classical well-posedness result.

\begin{proposition}\label{prop-well-posed}
  For any $\omega_0\in L^2(\T^2)$ and every fixed $\theta \in \cH_{\alpha}$, Eq. \eqref{SNSE-Ito} admits a pathwise unique weak solution $\omega$ in the sense of Definition \ref{def-stoch-NS/mSQG}, satisfying
  \begin{equation}\label{eq: well-posed-viscosity}
    \P\mbox{-a.s.}, \quad \|\omega_t \|_{L^2}^2 + 2\kappa\int_0^t \|\nabla \omega_s\|_{L^2}^2\, \d s = \|\omega_0 \|_{L^2}^2.
  \end{equation}
\end{proposition}

Our aim is to establish the following noise-induced anomalous regularity for the solution $\omega$: there exists a constant $C>0$, independent of $\kappa>0$, such that for some $\alpha^{\prime} \in (\alpha,1)$,
\begin{equation*}
\E \int_0^T \|\omega_s\|_{H^{1-\alpha^{\prime}}}^2 \, \d s \le C \|\omega_0 \|_{L^2}^2,
\end{equation*}
where $\alpha$ is the parameter appearing in the definition of $\cH_{\alpha}$. This will be proved in the following sections. Henceforth, we shall not mention explicitly that the noise coefficient lies in $\cH_{\alpha}$.

\subsection{Energy spectrum of solution}

Recall $\widehat{\omega}_k(t)= \<\omega_t, e_{-k}\>$ as the Fourier coefficients of $\omega_t$, and $\overline{\widehat{\omega}_k (t)}$ its complex conjugate. For simplicity of notation, we will often omit the time variable. We first prove the following result.

\begin{lemma}\label{lem-Fourier-coeff}
  It holds that, for any $k\in \Z_0^2$,
  $$\aligned
    \d|\widehat{\omega}_k|^2 &= 2\,{\rm Re} \big(\, \overline{\widehat{\omega}_k} \<\omega, u\cdot\nabla e_{-k}\> \big)\, \d t+ 2\,{\rm Re}\big(\widehat{\omega}_k\, \d\overline{M^k_t}\, \big) \\
    &\quad - 8\pi^2|k|^2 (\kappa+\nu) |\widehat{\omega}_k|^2\, \d t+ 16 \pi^2 \nu \sum_l  \theta_l^2 |\Pi_{l^\perp} k|^2 |\widehat{\omega}_{k-l}|^2\, \d t,
    \endaligned $$
  where ${\rm Re}$ means the real part, $\{M^k\}_k$ are martingales and $\Pi_{l^\perp} k$ is the orthogonal projection of $k$ in the direction $l^\perp$. Moreover, letting $Y_k= \E|\widehat{\omega}_k|^2$, then one has
  $$ \aligned
    \frac{\d}{\d t} Y_k &= 2\,{\rm Re} \E \big(\,\overline{\widehat{\omega}_k} \<\omega, u\cdot\nabla e_{-k}\> \big) - 8\pi^2|k|^2 (\kappa+\nu) Y_k + 16\pi^2 \nu \sum_l \theta_{l}^2 |\Pi_{l^\perp} k|^2 Y_{k-l} \\
    &= 2\,{\rm Re} \E \big(\, \overline{\widehat{\omega}_k} \<\omega, u\cdot\nabla e_{-k}\> \big) - 8\pi^2|k|^2 \kappa Y_k + 16\pi^2 \nu \sum_l \theta_{l}^2 |\Pi_{l^\perp} k|^2 (Y_{k-l} -Y_k).
    \endaligned $$
\end{lemma}

\begin{proof}
Taking $\phi =e_{-k}$ in \eqref{SNSE-Ito-weak} and writing the equation in differential form, we obtain
  $$\d\widehat{\omega}_k =\<\omega , u\cdot \nabla e_{-k}\>\, \d t- 4\pi^2|k|^2(\kappa + \nu) \widehat\omega_k\, \d t+ \sqrt{2 \nu} \sum_l \theta_l \<\omega, \, \sigma_l \cdot \nabla e_{-k}\> \, \d  B^l_t. $$
Denote the last martingale by $\d M^k_t$; we deduce from \eqref{covar-BM} the quadratic variation
  $$\aligned
    \d\big[M^k, \overline{M^k}\,\big]_t &= 4\nu \sum_l \theta_l^2 |\<\omega, \sigma_l \cdot \nabla e_{-k}\>|^2\, \d t
    = 16 \pi^2 \nu \sum_l \theta_l^2 \frac{(l^\perp\cdot k)^2}{|l|^2} |\<\omega, e_{l-k}\>|^2\, \d t\\
    &= 16\pi^2 \nu \sum_l \theta_l^2 |\Pi_{l^\perp} k|^2 |\widehat{\omega}_{k-l}|^2\, \d t.
    \endaligned $$
  Hence, by the It\^o formula,
  $$\aligned
    \d|\widehat{\omega}_k|^2
    &= \overline{\widehat{\omega}_k}\, \d\widehat{\omega}_k + \widehat{\omega}_k\, \d\overline{\widehat{\omega}_k} + \d\big[\widehat{\omega}_k, \overline{\widehat{\omega}_k}\, \big] \\
    &= 2\,{\rm Re} \big(\overline{\widehat{\omega}_k} \<\omega, u\cdot\nabla e_{-k}\> \big)\, \d t+ 2\,{\rm Re} \big( \widehat{\omega}_k\, \d \overline{M^k_t} \big) \\
    &\quad - 8\pi^2|k|^2 (\kappa+\nu) |\widehat{\omega}_k|^2\, \d t+ 16\pi^2 \nu \sum_l \theta_l^2 |\Pi_{l^\perp} k|^2 |\widehat{\omega}_{k-l}|^2\, \d t.
    \endaligned $$
  Taking expectation gives us the first equation for $Y_k$. Finally, using \eqref{covar-identity} and $|\Pi_{l^\perp} k|^2= \frac{(k\cdot l^\perp)^2}{|l|^2}$, we can derive the second equation of $Y_k$.
\end{proof}

Following \cite{Luo_2024_Elementary}, we study the time evolution of $\{Y_k\}_k$ as a sequence in $\ell^p(\Z_0^2)$, $p>1$.

\begin{proposition}\label{prop-Y_k^p}
For any $p > 1$, it holds
  $$\aligned
    & \frac{\d}{\d t} \sum_k Y_k^p + 8 \pi^2 p \nu \sum_{k,l}\theta_l^2 |\Pi_{l^\perp} k|^2 \big(Y_{k+l}^{p-1} -Y_k^{p-1}\big) (Y_{k+l} -Y_k) \\
    &\qquad = 2p \sum_k Y_k^{p-1} {\rm Re} \E \big(\,\overline{\widehat{\omega}_k} \<\omega, u\cdot\nabla e_{-k}\> \big)- 8\pi^2 p \kappa \sum_k |k|^2 Y_k^p.
    \endaligned $$
\end{proposition}

\begin{proof}
  We follow the proof of \cite[Lemma 3.1]{Luo_2024_Elementary}. For any $p \geq 1$, by Lemma \ref{lem-Fourier-coeff}, we have the following integral equation:
  \begin{equation} \label{eq: Fourier-coefficient-norm-2}
    \aligned
    Y_{k}^p(t)-Y_{k}^p(0) & = 2p\int_0^t Y_k^{p-1}{\rm Re} \E \big(\overline{\widehat{\omega}_k} \<\omega, u\cdot\nabla e_{-k}\> \big)\, \d s - 8 \pi^2 \kappa \, p|k|^2 \int_0^t  Y_{k}^p\, \d s\\
    &\quad + 16 \pi^2 \nu \, p \sum_l \theta_l^2 |\Pi_{l^\perp} k|^2 \int_{0}^{t} Y_k^{p-1} (Y_{k-l} -Y_k) \, \d  s.
    \endaligned
  \end{equation}
  By Proposition \ref{prop-well-posed}, we have $\sum_{k} Y_{k}(t) = \E \|\omega_t\|_{L^2}^2 \leq \|\omega(0)\|_{L^2}^2 < +\infty$, hence $\sup_{k} Y_{k}(t) \leq \|\omega(0)\|_{L^2}^2 < + \infty$ and
  $$\sum_{k} Y_{k}^{p}(t) \leq \sup\limits_{k} Y_{k}^{p-1}(t) \sum_{k} Y_{k}(t) \leq \|\omega(0)\|_{L^2}^{2p}  < +\infty, \quad \forall \, t \geq 0.$$
  Since $\sum_{k} Y_{k}^p (t)<+\infty$, summing equation \eqref{eq: Fourier-coefficient-norm-2} over $k \in \Z_0^2$ yields
  \begin{equation} \label{eq: Fourier sum}
    \aligned
    \sum_{k} Y_{k}^p(t) - \sum_{k} Y_{k}^p (0)
    &= 2p \sum_k \int_0^t Y_k^{p-1}{\rm Re} \E \big(\, \overline{\widehat{\omega}_k} \<\omega, u\cdot\nabla e_{-k}\> \big)\, \d s - 8 \pi^2 \kappa\, p \sum_k |k|^2 \int_0^t  Y_{k}^p\, \d s\\
    &\quad + 16 \pi^2 \nu \, p \sum_k \sum_l \theta_l^2 |\Pi_{l^\perp} k|^2 \int_{0}^{t} Y_k^{p-1} (Y_{k-l} -Y_k) \, \d  s.
    \endaligned
  \end{equation}

  Denote the last series on the right-hand side of equation \eqref{eq: Fourier sum} by
  $$S_0:= 16 \pi^2 \nu \, p \sum_k \sum_l \theta_l^2 |\Pi_{l^\perp} k|^2 \int_{0}^{t} Y_k^{p-1} (Y_{k-l} -Y_k) \, \d  s.$$
  According to Proposition \ref{prop-well-posed}, the process $\omega$ belongs to $ L^2(\Omega; L^2([0,T]; H^1(\T^2) ))$, therefore, the series $S_0$ converges absolutely. By Fubini's theorem, we have
  \begin{equation}\label{eq: Fourier sum-2}
    S_0= 16 \pi^2 \nu \, p \sum_l \sum_k \theta_l^2 |\Pi_{l^\perp} k|^2 \int_{0}^{t} Y_k^{p-1} (Y_{k-l} -Y_k) \, \d  s.
  \end{equation}
  For any fixed $l \in \Z_0^2$, since $\Pi_{l^\perp} l=0$, substituting $k$ for $\tilde{k}+l$ on the right-hand side of the above equality yields
  \begin{equation*}
    S_0= 16 \pi^2 \nu \, p \sum_l \sum_{\tilde k} \theta_l^2 |\Pi_{l^\perp} \tilde k|^2 \int_{0}^{t} Y_{\tilde{k}+l}^{p-1} \big( Y_{\tilde k} -Y_{\tilde{k}+l} \big) \, \d  s.
  \end{equation*}
  Substituting $l$ for $-\tilde{l}$ in the identity \eqref{eq: Fourier sum-2} and using the radial symmetry of $\theta_k$ in $k$, we get
  \begin{equation*}
    S_0= 16 \pi^2 \nu \, p \sum_{\tilde l} \sum_k \theta_{\tilde l}^2 \big|\Pi_{{\tilde l}^\perp} k \big|^2 \int_{0}^{t} Y_k^{p-1} \big(Y_{k+ \tilde l} -Y_k \big) \, \d  s
  \end{equation*}
  Changing $\tilde{k}, \tilde{l}$ back to $k,l$ in the above two identities and summing them up, we obtain
  $$S_0 = -8\pi^2 \nu \, p \sum_l \sum_k \theta_l^2 \big|\Pi_{l^\perp} k \big|^2 \int_{0}^{t}\big(Y_{k+ l}^{p-1}- Y_k^{p-1}\big) \big(Y_{k+ l} -Y_k \big) \, \d  s. $$
  Replacing the last series of \eqref{eq: Fourier sum} by the above expression, we arrive at
  \begin{equation*}
    \aligned
    \sum_{k} Y_{k}^p(t) - \sum_{k} Y_{k}^p (0)
    &= 2p \sum_k \int_0^t Y_k^{p-1}{\rm Re} \E \big(\, \overline{\widehat{\omega}_k} \<\omega, u\cdot\nabla e_{-k}\> \big)\, \d s - 8 \pi^2 \kappa\, p \sum_k |k|^2 \int_0^t  Y_{k}^p\, \d s\\
    &\quad -8\pi^2 \nu \, p \sum_l \sum_k \theta_l^2 \big|\Pi_{l^\perp} k \big|^2 \int_{0}^{t}\big(Y_{k+ l}^{p-1}- Y_k^{p-1}\big) \big(Y_{k+ l} -Y_k \big) \, \d  s.
    \endaligned
  \end{equation*}
  We obtain the desired result by differentiating the above equation.
\end{proof}

Next, we provide some key estimates on the nonlinearity
\begin{equation}\label{def-J-beta-p}
  J_{\beta,p}:= 2p \sum_k Y_k^{p-1} {\rm Re} \E \big(\,\overline{\widehat{\omega}_k} \<\omega, u\cdot\nabla e_{-k}\> \big),
\end{equation}
where $u= \nabla^\perp(-\Delta)^{-\beta/2} \omega = \frac{\rm i}{(2\pi)^{ \beta -1}} \sum_m \frac{m^\perp}{|m|^{ \beta}} \widehat{\omega}_m e_m$ with some $\beta \in (1,2]$.

\begin{proposition} \label{prop-nonlinearity}
  For any $\alpha\in \big(0,\frac{\beta-1}{2} \big)$, $1<p \le  \frac{\beta-2\alpha+1}{2}$ and $\lambda \geq \frac{4(1-\alpha)}{\beta-2\alpha-1}$, there exist constants $C_{\alpha,\beta}, C_{\alpha,\beta,p,\lambda} >0$ such that for any $\delta>0$, it holds
  \begin{align}
    J_{\beta,p} & \le C_{\alpha,\beta}  \|\omega_0\|_{L^2} \sum_k |k|^{2(1-\alpha)}\, Y_k^p, \label{eq-nonlinearity 1}                                                                                        \\
    J_{\beta,p} & \le \delta (p-1) \sum_k |k|^{2(1-\alpha)} Y_k^p +  \frac{C_{\alpha,\beta,p,\lambda}}{((p-1) \delta)^{\lambda-1}}  \|\omega_0\|_{L^2}^{\lambda}  \sum_{k} Y_{k}^p. \label{eq-nonlinearity 2}
  \end{align}
  Moreover, estimate \eqref{eq-nonlinearity 1} holds with $p=1$ and $\lim\limits_{p \to 1} C_{\alpha,\beta,p,\lambda} <\infty$.
\end{proposition}

\begin{proof}
  We first consider the case $p>1$.
  As $u=  \frac{\rm i}{(2\pi)^{\beta -1}} \sum_{m\ne 0} \frac{m^\perp}{|m|^{\beta}} \widehat{\omega}_m e_m$, it is clear that
  $$ |\<\omega, u\cdot\nabla e_{-k}\> |
    = (2\pi)^{2-\beta} \Big|\sum_{m\ne 0} \frac{m^\perp\cdot k}{|m|^{\beta}} \widehat{\omega}_{k-m} \, \widehat{\omega}_m \Big|
    \le (2\pi)^{2-\beta} \sum_{m\ne 0} \frac{|m^\perp\cdot k|}{|m|^{\beta}} |\widehat{\omega}_{k-m}|\, | \widehat{\omega}_m|. $$
  We have the following simple facts:
  $$|m^\perp\cdot k| \le \min\{|m|\, |k|, |m|\, |k-m|, |k-m|\, |k|\}; $$
  taking $0 \leq \alpha_1, \alpha_2\le 1$ and small enough $ \epsilon >0$ to be determined later, such that
  \begin{equation}\label{alpha-1-2}
    3 - \beta +\epsilon \leq \alpha_1+\alpha_2 \leq 2,
  \end{equation}
  one easily deduces
  $$|m^\perp\cdot k| \le |k|^{\alpha_1} |k-m|^{\alpha_2} |m|^{2-\alpha_1 -\alpha_2}. $$
  Therefore,
  $$|\<\omega, u\cdot\nabla e_{-k}\> |\le  (2\pi)^{2-\beta} |k|^{\alpha_1} \sum_{m\ne 0,k} \frac{|k-m|^{\alpha_2} }{|m|^{\beta+\alpha_1 + \alpha_2-2}}|\widehat{\omega}_{k-m}| \, |\widehat{\omega}_m|. $$
  By the definition \eqref{def-J-beta-p} of $J_{\beta,p}$ and Cauchy--Schwarz inequality,
  \begin{equation} \label{eq-J}
    \aligned
    J_{\beta,p}
    & \leq 2p (2\pi)^{2-\beta} \E \bigg[ \sum_{k\ne 0}  |k|^{\alpha_1} Y_k^{p-1}  |\widehat{\omega}_k| \sum_{m\ne 0,k}  \frac{ |\widehat{\omega}_m|  }{|m|^{\beta+\alpha_1 + \alpha_2-2}} |k-m|^{\alpha_2} |\widehat{\omega}_{k-m}|  \bigg] \\
    & \leq 2p (2\pi)^{2-\beta} \bigg[\E \sum_{k\ne 0}  |k|^{2 \alpha_1} Y_k^{2(p-1)}  |\widehat{\omega}_k|^2 \bigg]^{1/2} \\
    &\quad \times \bigg[ \E \sum_{k\ne 0} \Big( \sum_{m\ne 0,k}  \frac{ |\widehat{\omega}_m|  }{|m|^{\beta+\alpha_1 + \alpha_2-2}} |k-m|^{\alpha_2} |\widehat{\omega}_{k-m}| \Big)^2 \bigg]^{1/2}.
    \endaligned
  \end{equation}

  Firstly, we have
  \begin{equation} \label{eq-J-1}
    \aligned
    & \E \sum_k  |k|^{2 \alpha_1} Y_k^{2(p-1)} |\widehat{\omega}_k|^2= \sum_k |k|^{2 \alpha_1}  Y_k^{2p-1} \\
    & \leq \sup_{k} \Big( |k|^{\frac{2 p \alpha_1}{2p-1}} Y_{k}^p \Big)^{\!\frac{p-1}{p}}  \sum_k  |k|^{\frac{2 p \alpha_1}{2p-1}} Y_k^p \leq  \Big( \sum_k |k|^{\frac{2 p \alpha_1}{2p-1}} Y_k^p \Big)^{\frac{2p-1}{p}}.
    \endaligned
  \end{equation}

  Next, we estimate the last term in \eqref{eq-J}.  Let $C_\epsilon := \sum_m |m|^{-2(1+\epsilon)} < + \infty$. By Cauchy--Schwarz inequality and \eqref{alpha-1-2}, $\P$-a.s.,
  \begin{align*}
    \sum_m \frac{ |\widehat{\omega}_m| }{|m|^{\beta+\alpha_1 + \alpha_2-2}} & \leq \Big( \sum_{m} \frac{ 1  }{|m|^{2(\beta+\alpha_1 + \alpha_2-2)}} \Big)^{1/2} \Big( \sum_{m} |\widehat{\omega}_m|^2 \Big)^{1/2}
    \leq C_\epsilon^{1/2} \| \omega_0 \|_{L^2},
  \end{align*}
  where in the last step we have used the pathwise bound \eqref{eq: well-posed-viscosity}. Thus, we can apply the Cauchy--Schwarz inequality to the sum over $m$ and obtain
  \begin{equation}\label{eq-nonlin-2} 
  \aligned
     & \E \sum_k \Big( \sum_m \frac{ |\widehat{\omega}_m| }{|m|^{\beta+\alpha_1 + \alpha_2-2}} |k-m|^{\alpha_2} |\widehat{\omega}_{k-m}|  \Big)^2 \\
     & \leq C_\epsilon^{1/2} \| \omega_0 \|_{L^2} \, \E \sum_k \sum_m  \frac{ |\widehat{\omega}_m|  }{|m|^{\beta+\alpha_1 + \alpha_2-2}} |k-m|^{2 \alpha_2} \, |\widehat{\omega}_{k-m}|^2 \\
     & = C_\epsilon^{1/2} \| \omega_0 \|_{L^2} \, \E \sum_m  \frac{ |\widehat{\omega}_m|  }{|m|^{\beta+\alpha_1 + \alpha_2-2}} \sum_k |k|^{2 \alpha_2} \, |\widehat{\omega}_{k}|^2 \\
     & \leq C_{\epsilon} \| \omega_0 \|_{L^2}^2 \sum_{k} |k|^{2 \alpha_2} Y_{k}.
  \endaligned
  \end{equation}
  Using the H\"older inequality, we get
  \begin{align*}
    \sum_{k} |k|^{2 \alpha_2} Y_{k} 
    & =  \sum_{k} |k|^{-\frac{2(1+\epsilon)(p-1)}{p}} |k|^{2 \alpha_2 + \frac{2 (1+\epsilon)(p-1)}{p}} Y_{k} \\
    & \leq  \bigg( \sum_{k} \frac{1}{|k|^{2(1+\epsilon)}}  \bigg)^{\!\frac{p-1}{p}} \bigg( \sum_{k} |k|^{2 ( p \alpha_2 + (1+\epsilon)(p-1) )  } Y_{k}^p \bigg)^{\!\frac1p}.
  \end{align*}
  Substituting this estimate into the above inequality yields
  \begin{equation} \label{eq-J-2}
    \E \sum_k \Big( \sum_m  \frac{ |\widehat{\omega}_m|  }{|m|^{\beta+\alpha_1 + \alpha_2-2}} |k-m|^{\alpha_2} |\widehat{\omega}_{k-m}|  \Big)^2 \leq C_{\epsilon}^{\frac{2p-1}{p}} \! \| \omega_0 \|_{L^2}^2 \Big( \sum_{k} |k|^{2 ( p \alpha_2 + (1+\epsilon)(p-1) )} Y_{k}^p \Big)^{1/p}.
  \end{equation}

  Combining the estimates \eqref{eq-J}, \eqref{eq-J-1} and \eqref{eq-J-2}, we have
  \begin{equation} \label{eq-J-3}
    J_{\beta,p}  \leq  2p (2\pi)^{2-\beta} C_{\epsilon}^{\frac{2p-1}{2p}} \| \omega_0 \|_{L^2} \Big( \sum_k |k|^{q_1} Y_{k}^p \Big)^{\frac{2p-1}{2p}} \Big( \sum_k |k|^{q_2} Y_{k}^p \Big)^{\frac{1}{2p}},
  \end{equation}
  where
  $$q_1= \frac{2 p\alpha_1}{2p-1}, \quad q_2=2 p \alpha_2 + 2(1+\epsilon)(p-1).$$
  Taking $\epsilon=\frac{(\beta-2\alpha-1)^2}{8}<\frac{\beta-2\alpha-1}{8}$ and
  \begin{equation*}
    \alpha_1=1-\alpha, \quad \alpha_2=2-\beta+\alpha+\epsilon,
  \end{equation*}
  then $\alpha_1, \alpha_2 \in (0,1)$ and $\alpha_1+\alpha_2=3-\beta+\epsilon$, thus \eqref{alpha-1-2} holds. Now for all $1<p\le \frac{\beta-2\alpha+1}{2}$, one has $q_1 < 2 \alpha_1= 2(1-\alpha)$ and
  \begin{align*}
    q_2 & \le (\beta-2\alpha+1) (2-\beta+\alpha+\epsilon) + (\beta-2\alpha-1)(1+\epsilon)         \\
        & = 2(\beta-2\alpha) \epsilon + (\beta-2\alpha+1) (2-\beta+\alpha ) + (\beta-2\alpha-1) \\
        & = 2(\beta-2\alpha) \epsilon + 2 (1-\alpha)- (\beta-2\alpha-1)(\beta-\alpha-1)         \\
        & < 2(1-\alpha),
  \end{align*}
  where in the last step we have used $\epsilon=\frac{(\beta-2\alpha-1)^2}{8}  < \frac{(\beta-2\alpha-1)^2}{4(\beta-2\alpha)}<\frac{(\beta-2\alpha-1)(\beta-\alpha-1)}{2(\beta-2\alpha)}$, which follows from $\beta\in (1,2]$ and $\alpha \in (0,\frac{\beta-1}{2})$.
  Thus, by \eqref{eq-J-3}, we obtain the estimate \eqref{eq-nonlinearity 1}  with constant
  $$ C_{\alpha,\beta}:= \sup_{p \in [1,2]} \bigg\{ 2p (2\pi)^{2-\beta} C_{\frac{(\beta-2\alpha-1)^2}{8}}^{\frac{2p-1}{2p}} \bigg\}.$$

  For the case $p=1$, one can deduce from \eqref{eq-J} and \eqref{eq-nonlin-2} that
  $$J_{\beta,1}\le 2(2\pi)^{2-\beta} C_{\epsilon}^{\frac{1}{2}} \| \omega_0 \|_{L^2} \Big( \sum_k |k|^{2\alpha_1} Y_{k} \Big)^{\frac{1}{2}} \Big( \sum_k |k|^{2\alpha_2} Y_{k} \Big)^{\frac{1}{2}} .$$
  Choosing $\alpha_1$ and $\alpha_2$ as above, it is easy to show that $\alpha_2<1-\alpha$, hence we obtain the desired estimate in this case.

  Finally, we turn to deriving the estimate \eqref{eq-nonlinearity 2} from \eqref{eq-J-3} in the case $p>1$. By the H\"older inequality, for $i=1,2$,
  \begin{equation*}
    \sum_{k} |k|^{q_{i}} Y_{k}^{p} \leq \Big(\sum_{k} |k|^{2(1-\alpha)} Y_{k}^p \Big)^{\frac{q_i}{2(1-\alpha)}} \Big(\sum_{k}  Y_{k}^p \Big)^{1-\frac{q_i}{2(1-\alpha)}}.
  \end{equation*}
  Then \eqref{eq-J-3} becomes
  \begin{equation} \label{eq-J-4}
    J_{\beta,p}  \leq  2p (2\pi)^{2-\beta} C_{\epsilon}^{\frac{2p-1}{2p}} \|\omega_0 \|_{L^2} \Big(\sum_{k} |k|^{2(1-\alpha)} Y_{k}^p \Big)^{\vartheta_{p}} \Big(\sum_{k}  Y_{k}^p \Big)^{1-\vartheta_{p}},
  \end{equation}
  where $\vartheta_{p}:=\frac{2p-1}{2p} \frac{q_1}{2(1-\alpha)}  + \frac{1}{2p}  \frac{q_2}{2(1-\alpha)} \in (0,1)$ since $0< q_i <2(1-\alpha)$. Recalling the definitions of $q_1$ and $q_2$, we have
  \begin{equation*}
    \aligned
    \vartheta_{p} &= \frac{1}{2(1-\alpha)} \Big(\alpha_1+\alpha_2 + (1+\epsilon)\frac{p-1}{p} \Big) = \frac{1}{2(1-\alpha)} \Big(3-\beta+ \epsilon+ (1+\epsilon)\frac{p-1}{p} \Big) \\
    &\le \frac{1}{2(1-\alpha)} \Big(3-\beta+ \epsilon+ (1+\epsilon)\frac{\beta-2\alpha-1}{\beta-2\alpha+1} \Big)
    \endaligned
  \end{equation*}
  for all $1<p \leq \frac{\beta-2\alpha+1}{2}$. Then,
  \begin{align*}
    1-\vartheta_{p}
     & > 1- \frac{1}{2(1-\alpha)} \Big(3-\beta+ \epsilon+ (1+\epsilon)\frac{\beta-2\alpha-1}{\beta-2\alpha+1} \Big) \\
     & = \frac{(\beta-2\alpha)(\beta-2\alpha-1) -2(\beta-2\alpha)\epsilon}{2(1-\alpha)(\beta-2\alpha+1)}            \\
     & >\frac{\beta-2\alpha-1}{4(1-\alpha)},
  \end{align*}
  where the last step is due to $\epsilon=\frac{(\beta-2\alpha-1)^2}{8} < \frac{(\beta-2\alpha-1)^2}{4(\beta-2\alpha)}$. Finally, applying Young's inequality to \eqref{eq-J-4}, for any $\lambda \geq \frac{1}{1-\vartheta_{p}}$, there is constant $C_{\alpha,\beta,p,\lambda} >0$ such that
  \begin{align*}
    J_{\beta,p} & \leq \delta (p-1) \bigg( \Big( \sum_k |k|^{2(1-\alpha)} Y_{k}^p \Big)^{\vartheta_{p}} \Big(  \sum_{k} Y_{k}^p \Big)^{(1-\vartheta_{p}-\frac{1}{\lambda})} \bigg)^{\frac{\lambda}{\lambda-1}} + \frac{C_{\alpha,\beta,p,\lambda}}{((p-1) \delta)^{\lambda-1}}  \|\omega_0\|_{L^2}^{\lambda}  \sum_{k} Y_{k}^p \\
                & \le \delta (p-1) \Big( \sum_k |k|^{2(1-\alpha)} Y_{k}^p \Big)^{\vartheta_{p} \frac{\lambda}{\lambda-1} + (1-\vartheta_{p}-\frac{1}{\lambda}) \frac{\lambda}{\lambda-1}}  + \frac{C_{\alpha,\beta,p,\lambda}}{((p-1) \delta)^{\lambda-1}}  \|\omega_0\|_{L^2}^{\lambda}  \sum_{k} Y_{k}^p                     \\
                & = \delta (p-1) \sum_k |k|^{2(1-\alpha)} Y_{k}^p + \frac{C_{\alpha,\beta,p,\lambda}}{((p-1) \delta)^{\lambda-1}}  \|\omega_0\|_{L^2}^{\lambda}  \sum_{k} Y_{k}^p,
  \end{align*}
  and $\lim_{p \rightarrow 1} C_{\alpha,\beta,p,\lambda} <+\infty$. Then for any $\lambda \geq \frac{4(1-\alpha)}{\beta-2\alpha-1}$, we obtain \eqref{eq-nonlinearity 2}.
\end{proof}

\subsection{Anomalous regularity and dissipation} \label{subsec: anomalous regularity and dissipation}

Having Proposition \ref{prop-nonlinearity} at hand, now we are ready to show the anomalous regularity and dissipation of solutions to \eqref{SNSE-Ito}. Due to the presence of the molecular diffusivity, we have the energy equality \eqref{eq: well-posed-viscosity} which implies $\sum_{k \in \Z_0^2} |k|^2 Y_{k} <+\infty$ for a.e. $t \in [0,T]$. Hence, we can apply the discrete Poincar\'e inequality \eqref{def-C-theta} in Section \ref{subsec: noise}: there exists a uniform constant $C_0>0$ such that, for every $p>1$,
\begin{equation}\label{discrete-Poincare-1}
  8 \pi^2 p \sum_{k,l} \theta_{l}^2  |\Pi_{l^\perp} k|^2  \big(Y_{k+l}^{p-1} -Y_k^{p-1}\big) (Y_{k+l} -Y_k) \geq C_0 \frac{p-1}{p} \sum_{k} |k|^{2(1-\alpha)} Y_{k}^p.
\end{equation}

\begin{theorem}[Anomalous regularity and dissipation] \label{thm: anomalous regularity}
For every $\alpha\in(0,\frac{\beta-1}{2})$ and every $\{\theta_k\}_k\in\cH_{\alpha}$, the following assertions hold uniformly over $\cH_{\alpha}$.
  \begin{itemize}
    \item[\rm (i)] {\rm (Small initial data)}
    If $\|\omega_0\|_{L^2} < \frac{\beta-2\alpha-1}{\beta-2\alpha+1} \frac{C_0 \nu}{C_{\alpha,\beta}}$, where the constants $C_0$ and $C_{\alpha,\beta}$ are given in \eqref{discrete-Poincare-1} and \eqref{eq-nonlinearity 1} respectively, then there exists $1<p<\frac{\beta-2\alpha+1}{2}$ such that
    \begin{equation}\label{def-mu-p}
    \mu_{p}:= \frac{p-1 }{p^2} C_0 \nu- \frac{1}{p} C_{\alpha,\beta} \|\omega_0\|_{L^2}>0.
    \end{equation}
    For any $\alpha^{\prime} \in (2-\frac{2-\alpha}{p}, \frac{1}{p})$ and for all $t \geq s \geq 0 $, it holds that
    \begin{equation*}
    \E  \int_{s}^{t}  \| \omega(\tau) \|_{H^{1-\alpha^{\prime}}}^2  \, \d \tau \lesssim_{\alpha,\alpha^{\prime},p} \mu_{p}^{-1/p} (t-s)^{\frac{p-1}{p}}  e^{- \mu_{p} s} \|\omega_0\|_{L^2}^2.
    \end{equation*}
    Furthermore, there is a constant $C_{\alpha,p,\mu_{p}}>0$ such that
    \begin{equation}\label{eq-anomal-dissip}
    \E \| \omega_t \|_{L^2}^2 \leq C_{\alpha,p,\mu_{p}} e^{- \mu_{p} t} \|\omega_0\|_{L^2}^2.
    \end{equation}
    \item[\rm (ii)] {\rm (General initial data)} For any $1<p \leq \frac{\beta-2\alpha+1}{2}$, $\lambda \geq \frac{4(1-\alpha)}{\beta-2\alpha-1}$ and $T\geq 0$, one has
    \begin{equation*}
    \int_{0}^T  \| \omega_t \|_{\mathcal{X}_{\alpha,p}}^{2p} \, \d t\leq  \bigg( \frac{2 p}{(p-1)C_0 \nu} + C_{\alpha,\beta,p,\lambda} T \Big( \frac{2p}{(p-1)C_0 \nu} \Big)^{\lambda} \|\omega_0\|_{L^2}^{\lambda} \bigg) \|\omega_0\|_{L^2}^{2p},
    \end{equation*}
    where the constant $C_{\alpha,\beta,p,\lambda}>0$ is given in \eqref{eq-nonlinearity 2}. Furthermore, for any $\alpha^{\prime} \in (\alpha,1)$, $1<p<\min\{\frac{\beta-2\alpha+1}{2}, \frac{2-\alpha}{2-\alpha^{\prime}}  \}$ and for all $T \geq 0$, it holds
    \begin{equation*}
    \E \int_{0}^T \|\omega_t \|_{H^{1-\alpha^{\prime}}}^2 \,\d t \lesssim_{\alpha,\alpha^{\prime},\beta,\lambda,p} T^{\frac{p-1}{p}} \Big( \frac{1}{\nu^{1/p}}+ \frac{\|\omega_0\|_{L^2}^{\lambda/p}}{\nu^{\lambda/p}} \Big) \|\omega_0\|_{L^2}^2.
    \end{equation*}
  \end{itemize}
\end{theorem}

\begin{remark}
When $\|\omega_0\|_{L^2}$ is small or the noise intensity $\nu$ is large, estimate \eqref{eq-anomal-dissip} implies anomalous dissipation, and one has an estimate on the Lyapunov exponent:
  $$\limsup_{t\to +\infty} \frac1t \log \E\|\omega_t \|_{L^2}^2 \le -\mu_p. $$
This estimate is extrinsic since it depends on the parameter $p$. In the general setting, for any fixed $\kappa>0$, since $\|\omega_t\|_{L^2}\to 0$ as $t\to\infty$, there exists $t_{\kappa}>0$ such that $\|\omega(t_{\kappa})\|_{L^2}$ is sufficiently small. Hence, for $t>t_{\kappa}$, repeating the proof of \eqref{eq-anomal-dissip} with $p$ close to $\frac{\beta-2\alpha+1}{2}$ shows that the Lyapunov exponent is independent of $\kappa$ and $p$:
\begin{equation*}
  \limsup_{\kappa \rightarrow 0} \lim_{t \rightarrow +\infty} \frac1t \log \E \| \omega_t \|_{L^2}^2 \le - \frac{\beta-2\alpha-1}{(\beta-2\alpha+1)^2} C_0 \nu.
\end{equation*}
\end{remark}

\begin{proof}[Proof of Theorem \ref{thm: anomalous regularity}]
  We treat the small-data and general-data cases separately.

  (i) Recall the nonlinear term $J_{\beta,p}$ defined in \eqref{def-J-beta-p}. For any $1<p < \frac{\beta-2\alpha+1}{2}$, combining Proposition \ref{prop-Y_k^p} and the estimates \eqref{eq-nonlinearity 1}, \eqref{discrete-Poincare-1}, we get
  \begin{align*}
     & \frac{\d}{\d t} \sum_k Y_k^p  + 8\pi^2 p \kappa \sum_k |k|^2 Y_k^p                                                           \\
     & = J_{\beta,p} - 8\pi^2 p \nu \sum_{k,l} \theta_l^2 |\Pi_{l^\perp} k|^2 \big(Y_{k+l}^{p-1} -Y_k^{p-1}\big) (Y_{k+l} -Y_k)  \\
     & \leq   \Big( C_{\alpha,\beta} \| \omega_0 \|_{L^2}-  \frac{p-1 }{p} C_0 \nu \Big)  \sum_k |k|^{2(1-\alpha)} Y_k^p.
  \end{align*}
  We denote by $R= \| \omega_0 \|_{L^2}$ and let $ \mu_{p}$ be defined as in \eqref{def-mu-p}. Due to $\|\omega_0\|_{L^2} < \frac{\beta-2\alpha-1}{\beta-2\alpha+1} \frac{C_0 \nu}{C_{\alpha,\beta}}$, there exists $p>1$ slightly less than $1+\frac{\beta-2\alpha-1}{2}$ such that $\mu_{p}>0$; then we obtain
  \begin{equation}\label{eq-anomal-1}
    \frac{\d}{\d t} \sum_k Y_k^{p} \le -p\, \mu_{p} \sum_k |k|^{2(1-\alpha)} Y_k^{p} \le -p\, \mu_{p} \sum_k Y_k^{p}.
  \end{equation}
  Using the Gr\"onwall inequality and fact $\sum_{k} Y_{k}^{p}(0)\le \big(\sup_k Y_{k}^{p-1}(0)\big) \sum_{k} Y_{k}(0) \leq R^{2p}$, we have
  \begin{equation} \label{eq-anomalous-a-diss-p}
    \sum_k Y_k^{p}(t) \leq e^{- p \, \mu_{p} t} \sum_k Y_k^{p}(0) \leq e^{- p \, \mu_{p} t} R^{2p}, \quad \forall\, t \geq 0 .
  \end{equation}
  On the other hand, integrating the first inequality in \eqref{eq-anomal-1} in time yields that
  \begin{equation} \label{eq-anomalous-regul-a-p}
    \sum_k Y_k^{p}(t) + p \, \mu_{p} \int_{s}^t \sum_k |k|^{2(1-\alpha)} Y_k^{p} (\tau)  \, \d \tau \leq   \sum_k Y_k^{p}(s) \leq e^{- p \, \mu_{p} s} R^{2p}
  \end{equation}
  for any $t \geq s \geq 0$.

  Since $1<p< 1+\frac{\beta-2\alpha-1}{2}<\frac{3-\alpha}{2}$, it is easy to check that
  \begin{equation*}
    \alpha<2-\frac{2-\alpha}{p} < \frac{1}{p}.
  \end{equation*}
  Now for any $\alpha^{\prime} \in (2-\frac{2-\alpha}{p}, \frac{1}{p})$, we have
  \begin{equation*}
    (1-\alpha)-p(1-\alpha^{\prime})>p-1\quad \mbox{and}\quad \frac{p}{p-1}(1-\alpha^{\prime}) >1,
  \end{equation*}
  which will be utilized in the subsequent proofs of the anomalous regularity and anomalous dissipation, respectively.

  For convenience, we denote
  \begin{equation*}
    H_{\alpha^{\prime}}(t):= \sum_{k} |k|^{2(1-\alpha^{\prime})} Y_{k}(t),\quad t\ge 0.
  \end{equation*}
  Using H\"older's inequality, we get that for any $t>s \geq 0$,
  \begin{align*}
    \int_{s}^{t} H_{\alpha^{\prime}}(\tau) \,\d \tau
     & = \int_{s}^t \sum_k |k|^{2(1-\alpha^{\prime})-2(1-\alpha)/p} |k|^{2(1-\alpha)/p} Y_k(\tau) \,\d \tau           \\
     & \leq  \bigg( \sum_{k} \frac{t-s}{|k|^{2 \frac{(1-\alpha)-(1-\alpha^{\prime})p}{p-1}}}  \bigg)^{\frac{p-1}{p}}
    \bigg(\int_{s}^t  \sum_k |k|^{2(1-\alpha)} Y_k^{p}(\tau)\,\d \tau \bigg)^\frac{1}{p} .
  \end{align*}
Since $(1-\alpha)-(1-\alpha^{\prime})p >p-1$, one has $ \sum_{k} |k|^{-2 \frac{(1-\alpha)- (1-\alpha^{\prime}) p}{p-1}} <+\infty$. Thus, by the estimate \eqref{eq-anomalous-regul-a-p}, we obtain the anomalous regularity for small initial data:
  \begin{equation} \label{eq-anomalous-regul-a}
    \int_{s}^{t} H_{\alpha^{\prime}}(\tau) \,\d \tau \lesssim_{\alpha,\alpha^{\prime},p} \mu_{p}^{-1/p} (t-s)^{\frac{p-1}{p}} e^{- \mu_{p} s} R^2.
  \end{equation}

Next, we turn to the anomalous dissipation for small initial data. Using the non-increasing property of $t\mapsto \E \| \omega_t \|_{L^2}^2$ and H\"older's inequality, it holds that
  \begin{align*}
    \E\| \omega_t \|_{L^2}^2
     & \le \int_{t-1}^{t} \sum_{k} Y_{k}(s) \,\d s = \int_{t-1}^{t} \sum_k Y_{k}^{\frac12}(s) |k|^{-(1-\alpha^{\prime})} \big( |k|^{1-\alpha^{\prime}} Y_{k}^{\frac12} (s) \big) \,\d s \\
     & \leq \bigg[ \int_{t-1}^t \sum_{k} Y_{k}^{p} (s) \,\d s \bigg]^{\frac{1}{2p}} \bigg[\sum_{k} |k|^{-\frac{2p(1-\alpha^{\prime})}{p-1}} \bigg]^{\frac{p-1}{2p}} \bigg[\int_{t-1}^t \sum_{k} |k|^{2(1-\alpha^{\prime})} Y_{k}(s) \,\d s \bigg]^{\frac{1}{2}}
  \end{align*}
for any $t \geq 1$. As $\frac{p (1-\alpha^{\prime})}{p-1}>1$, we have $C_{p,\alpha^{\prime}}:=\sum_{k} |k|^{-\frac{2p(1-\alpha^{\prime})}{p-1}} < +\infty$. Using the anomalous regularity \eqref{eq-anomalous-regul-a} for small initial data and estimate \eqref{eq-anomalous-a-diss-p}, we get
  \begin{align*}
    \E\| \omega_t \|_{L^2}^2
     & \leq \bigg[ \sup_{s \in [t-1,t]} \sum_{k} Y_{k}^{p} (s)  \bigg]^{\frac{1}{2p}} C_{p,\alpha^{\prime}} \bigg[\int_{t-1}^t H_{\alpha^{\prime}}(s) \,\d s \bigg]^{\frac{1}{2}} \\
     & \lesssim_{\alpha,\alpha^{\prime},p} C_{p,\alpha^{\prime}}\, e^{- \frac{\mu_{p}}{2} (t-1)} R \, \mu_{p}^{-1/2 p} e^{- \frac{\mu_{p}}{2} (t-1)} R               \\
     & \leq C_{\alpha,\alpha^{\prime},p,\mu_{p}} e^{- \mu_{p} t} R^2,
  \end{align*}
where the constant $C_{\alpha,\alpha^{\prime},p,\mu_{p}}<+\infty$ depends only on $\alpha$, $\alpha^{\prime}$, $p$ and $\mu_{p}$. Taking 
  $$\alpha^{\prime} = \frac12 \Big(2-\frac{2-\alpha}{p}+ \frac{1}{p} \Big) =1-\frac{1-\alpha}{2p},$$ 
we obtain the anomalous dissipation with the new constant $C_{\alpha,p,\mu_{p}}:= C_{\alpha,1-\frac{1-\alpha}{2p}, p,\mu_{p}}$.

(ii) Finally, we aim to prove the anomalous regularity for general initial data. For any $1<p \leq 1+ \frac{\beta- 2\alpha-1}{2}$, combining Proposition \ref{prop-Y_k^p} and estimates \eqref{eq-nonlinearity 2}, \eqref{discrete-Poincare-1}, we get
  \begin{align*}
     & \frac{\d}{\d t} \sum_k Y_k^p  + 8\pi^2 p \kappa \sum_k |k|^2 Y_k^p  \\
     & = J_{\beta,p} - 8\pi^2 p \nu \sum_{k,l} \theta_l^2 |\Pi_{l^\perp} k|^2 \big(Y_{k+l}^{p-1} -Y_k^{p-1}\big) (Y_{k+l} -Y_k) \\
     & \leq  (p-1) \Big( \delta-  \frac{C_0 \nu}{p} \Big) \sum_k |k|^{2(1-\alpha)} Y_k^p +  \frac{C_{\alpha,\beta,p,\lambda}}{ ((p-1)\delta)^{\lambda-1}} \, \|\omega_0\|_{L^2}^{\lambda} \sum_{k} Y_{k}^p
  \end{align*}
which holds for any $\delta>0$ and $\lambda \geq \frac{4(1-\alpha)}{\beta-2\alpha-1}$. Taking $\delta = \frac{C_0 \nu}{2 p}$ leads to
  $$\frac{\d}{\d t} \sum_k Y_k^p + \frac{p-1}{ 2 p} C_0 \nu \sum_k |k|^{2(1-\alpha)} Y_k^p \le C_{\alpha,\beta,p,\lambda} \Big(\frac{ 2 p }{  (p-1) C_0 \nu } \Big)^{\lambda-1} \|\omega_0\|_{L^2}^{\lambda} \sum_{k} Y_{k}^p. $$
Next, using the fact $\sum_{k} Y_{k}^p \leq \| \omega_0 \|_{L^2}^{2p} = R^{2p}$ and integrating from $0$ to $T$, we obtain
  \begin{equation}\label{eq-anomalous-b-diss-p}
    \int_{0}^T  \sum_k |k|^{2(1-\alpha)} Y_k^p(t)\, \d t\leq  \bigg( \frac{2 p}{(p-1)C_0 \nu} + C_{\alpha,\beta,p,\lambda} T \Big( \frac{2p}{(p-1)C_0 \nu} \Big)^{\lambda} R^{\lambda} \bigg) R^{2p}.
  \end{equation}
This provides the $L^{2p}([0,T];\mathcal{X}_{\alpha,p})$ regularity of $\omega$, which is independent of $\kappa>0$.

For any fixed $\alpha^{\prime} \in (\alpha,1)$ and $1<p<\min\{\frac{\beta-2\alpha+1}{2}, \frac{2-\alpha}{2-\alpha^{\prime}}  \}$, we have $(1-\alpha)-(1-\alpha^{\prime})p>p-1$ and estimate \eqref{eq-anomalous-b-diss-p} holds. Moreover, similarly to the proof of \eqref{eq-anomalous-regul-a}, 
  \begin{equation*}
  \aligned
  \E \int_{0}^T \|\omega_t \|_{H^{1-\alpha^{\prime}}}^2 \,\d t
  &= \int_{0}^T \sum_k |k|^{2(1-\alpha^{\prime})-2(1-\alpha)/p} |k|^{2(1-\alpha)/p} Y_k(t) \,\d t \\
  & \leq  \bigg( \sum_{k} \frac{T}{|k|^{2 \frac{(1-\alpha)-(1-\alpha^{\prime})p}{p-1}}}  \bigg)^{\frac{p-1}{p}}
    \bigg(\int_{0}^T  \sum_k |k|^{2(1-\alpha)} Y_k^p(t)\,\d t\bigg)^\frac{1}{p} \\
  &\lesssim_{\alpha,\alpha^{\prime},\beta,\lambda,p} T^{\frac{p-1}{p}} \Big( \frac{1}{\nu^{1/p}}+ \frac{R^{\lambda/p}}{\nu^{\lambda/p}} \Big) R^{2} 
  \endaligned
  \end{equation*}
for all $T\ge 0$. Thus, we get the anomalous regularity.
\end{proof}

\section{Existence and Pathwise uniqueness} \label{sec: pathwise uniqueness}

This section is devoted to proving Theorem \ref{thm-Euler/mSQG}. The proof relies on the anomalous regularity obtained in Theorem~\ref{thm: anomalous regularity} and the energy spectrum analysis.

Throughout this section, we work on $\T^2$ with noise coefficients $\{\theta_k \}_k \in \cH_{\alpha}$ for  $\alpha\in(0,\frac{\beta-1}{2})$. Recall that the It\^o form of the stochastic 2D Euler/mSQG equation \eqref{stoch-2D-Euler-1} reads as
\begin{equation} \label{2D Euler/mSQG}
  \d \omega + u\cdot\nabla\omega\, \d t + \sqrt{2\nu}\sum_k\theta_k\sigma_k\cdot\nabla\omega\, \d B_t^k = \nu \Delta\omega\, \d t,
\end{equation}
where $u = \nabla^\perp(-\Delta)^{-\beta/2}\omega$ with $\beta\in(1,2]$. Before moving forward, we first give the meaning of weak solutions to the above equation \eqref{2D Euler/mSQG}. 

\begin{definition} \label{def-stoch-Euler/mSQG}
  We say that \eqref{2D Euler/mSQG} has a weak solution if there exist a filtered probability space $(\Omega,\mathcal F, \mathcal F_t, \P)$, a sequence of $\mathcal F_t$-complex Brownian motions $\{B^k\}_k$ satisfying \eqref{covar-BM}, and a $\mathcal F_t$-progressively measurable process $\omega_t$ with trajectories in $C([0,T], L_w^2)$ such that for any $\phi\in C^\infty(\T^2)$, $\P$-a.s. for all $t\in [0,T]$, it holds
  \begin{equation} \label{SEE-Ito-weak}
    \aligned
    \<\omega_t, \phi\>
    &= \<\omega_0, \phi\>+ \int_0^t\<\omega_s, u_s\cdot \nabla \phi\>\, \d s + \nu \int_0^t \<\omega_s, \Delta\phi\>\, \d s \\
    &\quad + \sqrt{2 \nu} \sum_k \theta_k \int_0^t\<\omega_s, \sigma_k \cdot \nabla\phi\> \, \d  B^k_s,
    \endaligned
  \end{equation}
  where $u_s= \nabla^\perp (-\Delta)^{-\beta/2} \omega_s$, $s\in [0,T]$.
\end{definition}

The existence and anomalous regularity of solutions to \eqref{2D Euler/mSQG} in the space $L^{2p}([0,T]; \mathcal{X}_{\alpha,p}) \cap L^{\infty} \big(\Omega;L^{\infty}([0,T]; L^2(\T^2)) \big)$ follow from Theorem \ref{thm: anomalous regularity} and a standard compactness argument combined with the vanishing viscosity method. In the sequel, we concentrate on the pathwise uniqueness, and consider two different weak solutions $\omega^1$ and $\omega^2$ of \eqref{2D Euler/mSQG} in $L^{2p}([0,T];\mathcal{X}_{\alpha,p}) \cap L^{\infty} \big(\Omega; L^{\infty}([0,T]; L^2(\T^2)) \big)$ corresponding to the same initial condition $\omega_{0} \in L^2(\T^2)$ and the same sequence of $\mathcal F_t$-complex Brownian motions $\{B^k_t \}_k$. Denote by $\bar{\omega}=\omega^1-\omega^2$ the difference of the two solutions. Then, in the weak sense, it satisfies, 
\begin{equation} \label{2D Euler/mSQG diff}
  \d \, \bar{\omega} + u^1 \cdot \nabla \bar{\omega} \,\d t +  \bar{u}  \cdot \nabla \omega^2 \, \d t + \sqrt{2\nu}\sum_k\theta_k\sigma_k\cdot\nabla \bar{\omega} \,\d B_t^k = \nu \Delta \bar{\omega} \,\d t,
\end{equation}
where $u^i= \nabla^\perp (-\Delta)^{-\beta/2} \omega^{i}$ for $i=1,2$ and $\bar{u} = u^1-u^2$.

\subsection{Energy spectrum}

To obtain the pathwise uniqueness, analogously to Proposition \ref{prop-Y_k^p}, we will study the time evolution of energy spectrum of $\bar{\omega}$. Aside from the nonlinear term being different, the summation interchange of the It\^o-Stratonovich correction term needs to be reconsidered as well, owing to the lack of $L^{2} \big(\Omega; L^2([0,T];H^1(\T^2)) \big)$ regularity in the inviscid setting (i.e., 2D Euler/mSQG). Fortunately, the anomalous regularity in $L^{2p}([0,T];\mathcal{X}_{\alpha,p})$ and a cutoff technique ensure that the exchange of summation order is valid, see Lemma \ref{lem-sum-cutoff} for the details. Thanks to this fact, we can obtain the time evolution of energy spectrum for the difference $\bar{\omega}$ of two solutions to \eqref{2D Euler/mSQG}.

Let $\chi\in C_c^\infty(\mathbb R_+;[0,1])$ be a decreasing cutoff function, with $\chi(r)=1$ if $r\in [0,1]$ and $\chi(r)=0$ if $r \geq 2$; for any $N\ge 1$, we define $\chi_N(|k|):= \chi(|k|/N)$ for all $k \in \Z_0^2$.

\begin{proposition}\label{prop-Y_k^p-diff}
  Denote $X_{k}(t)=\E \big| \widehat{\bar{\omega}}_{k}(t) \big|^2$ with $k \in \Z_0^2$. For any $p \in (1,\frac{\beta-2\alpha+1}{2})$, it holds
  \begin{align*}
     & \sum_k X_k^p (t) +  \frac{p-1}{p} C_0 \nu  \int_{0}^{t}  \sum_{k} |k|^{2(1-\alpha)} X_{k}^p \,\d s\\
     & \leq 2p \lim_{N \rightarrow +\infty} \int_{0}^t \sum_k  \chi_{N}(|k|)  X_k^{p-1} {\rm Re} \E \Big[\, \overline{\widehat{ \bar{\omega}}_k} \big( \<\bar{\omega}, u^1\cdot\nabla e_{-k}\> + \< \omega^2, \bar{u} \cdot\nabla e_{-k}\> \big) \Big] \, \d s 
  \end{align*}
for any $t \in [0,T]$, where the constant $C_0>0$ is given in inequality \eqref{discrete-Poincare-1}.
\end{proposition}

\begin{proof}
Recall the equation \eqref{2D Euler/mSQG diff} for the difference $\bar{\omega}$ between two solutions; the Fourier coefficient $\widehat{\bar{\omega}}_{k}=\<\bar{\omega},e_{-k} \>$ satisfies
  \begin{align*}
    \d \widehat{ \bar{\omega}}_k & =\< \bar{\omega} , u^1 \cdot \nabla e_{-k}\>\, \d t + \< \omega^2 , \bar{u} \cdot \nabla e_{-k}\>\, \d t- 4\pi^2|k|^2 \nu \, \widehat{\bar{\omega}}_k\, \d t \\
    & \quad + \sqrt{2 \nu} \sum_l \theta_l \< \bar{\omega}, \sigma_l \cdot \nabla e_{-k}\> \, \d B^l_t.
  \end{align*}
Let $\d M^k_t$ be the last martingale; by \eqref{covar-BM}, it has the quadratic variation
  $$\aligned
    \d\big[M^k, \overline{M^k}\,\big]_t &= 4\nu \sum_l \theta_l^2 |\<\bar{\omega}, \sigma_l \cdot \nabla e_{-k}\>|^2\, \d t
    = 16 \pi^2 \nu \sum_l \theta_l^2 \frac{(l^\perp\cdot k)^2}{|l|^2} |\<\bar{\omega}, e_{l-k}\>|^2\, \d t\\
    &= 16\pi^2 \nu \sum_l \theta_l^2 |\Pi_{l^\perp} k|^2 |\widehat{ \bar{\omega}}_{k-l}|^2\, \d t.
    \endaligned $$
Hence, by the It\^o formula,
  $$\aligned
    \d|\widehat{\bar{\omega}}_k|^2
    &= \overline{ \widehat{ \bar{\omega}}_k}\, \d\widehat{\bar{\omega}}_k + \widehat{\bar{\omega}}_k\, \d\overline{\widehat{\bar{\omega}}_k} + \d\big[\widehat{\bar{\omega}}_k, \overline{\widehat{\bar{\omega}}_k}\, \big] \\
    &= 2\,{\rm Re}  \Big[\, \overline{\widehat{ \bar{\omega}}_k} \big( \<\bar{\omega}, u^1\cdot\nabla e_{-k}\> + \< \omega^2, \bar{u} \cdot\nabla e_{-k}\> \big) \Big] \d t+ 2\,{\rm Re} \big( \widehat{\bar{\omega}}_k\, \d \overline{M^k_t} \big) \\
    &\quad - 8\pi^2|k|^2 \nu |\widehat{\bar{\omega}}_k|^2\, \d t+ 16\pi^2 \nu \sum_l \theta_l^2 |\Pi_{l^\perp} k|^2 |\widehat{\bar{\omega}}_{k-l}|^2\, \d t.
    \endaligned $$
  Taking expectation gives us
  \begin{align*}
    \d X_{k} & = 2\,{\rm Re} \E \Big[\, \overline{\widehat{ \bar{\omega}}_k} \big( \<\bar{\omega}, u^1\cdot\nabla e_{-k}\> + \< \omega^2, \bar{u} \cdot\nabla e_{-k}\> \big) \Big]\d t + 16 \pi^2 \nu \sum_l \theta_l^2 |\Pi_{l^\perp} k|^2  (X_{k-l} -X_k) \d t.
  \end{align*}
  Then for all $t \in [0,T]$ and $p>1$, it is easy to show that $\sum _{k}\chi_{N}(|k|) X_{k}^p$ satisfies
  \begin{align*}
    \Big( \sum_k \chi_{N}(|k|) X_k^p \Big) (t) & =   2p  \int_{0}^t \sum_k  \chi_{N}(|k|)  X_k^{p-1} {\rm Re} \E \Big[\, \overline{\widehat{ \bar{\omega}}_k} \big( \<\bar{\omega}, u^1\cdot\nabla e_{-k}\> + \< \omega^2,\bar{u} \cdot\nabla e_{-k}\> \big) \Big] \d s \\
    & \quad + 16 \pi^2 p \, \nu \int_{0}^{t} \sum_k \sum_l \theta_l^2 |\Pi_{l^\perp} k|^2 \chi_{N}(|k|) X_k^{p-1} (X_{k-l} -X_k) \, \d  s.
  \end{align*}
  Due to $\omega^1,\omega^2 \in L^{2p}([0,T]; \mathcal{X}_{\alpha,p})$, one has $\int_0^T \sum_{k} |k|^{2(1-\alpha)} X_{k}^p \, \d s <+\infty$. By the monotone convergence theorem and Lemma \ref{lem-sum-cutoff}, as $N \rightarrow +\infty$, we obtain, for any $t \in [0,T]$,
  \begin{align*}
     & \sum_k X_k^p (t) +  8 \pi^2 p \, \nu \int_{0}^{t} \sum_{k,l} \theta_l^2 |\Pi_{l^\perp} k|^2 \big(X_{k+l}^{p-1} -X_k^{p-1}\big) (X_{k+l} -X_k) \, \d s  \\
     & = 2p \lim_{N \rightarrow +\infty} \int_{0}^t  \sum_k  \chi_{N}(|k|)  X_k^{p-1} {\rm Re} \E \Big[ \, \overline{\widehat{ \bar{\omega}}_k} \big( \<\bar{\omega}, u^1\cdot\nabla e_{-k}\> + \< \omega^2, \bar{u} \cdot\nabla e_{-k}\> \big) \Big] \d s.
  \end{align*}
  Using Proposition \ref{prop: weighted Poincare}, there is a constant $C_0>0$ such that
  \begin{equation*}
    8 \pi^2 p\nu \sum_{k,l} \theta_{l}^2  |\Pi_{l^\perp} k|^2  \big(X_{k+l}^{p-1} -X_k^{p-1}\big) (X_{k+l} -X_k) \geq C_0 \nu  \frac{p-1}{p}  \sum_{k} |k|^{2(1-\alpha)} X_{k}^p
  \end{equation*}
  for any $p>1$. Thus, we get the desired result.
\end{proof}

\subsection{Estimates of nonlinear terms}

Denote the nonlinear terms as follows:
\begin{align*}
  J_{\beta,p}^{(1)} & := 2p \sum_k   X_k^{p-1} \Big| {\rm Re} \E \Big[\, \overline{\widehat{ \bar{\omega}}_k}\, \<\bar{\omega}, u^1\cdot\nabla e_{-k}\> \Big] \Big| ;     \\
  J_{\beta,p}^{(2)} & := 2p \sum_k   X_k^{p-1} \Big| {\rm Re} \E \Big[\, \overline{\widehat{ \bar{\omega}}_k} \, \< \omega^2, \bar{u} \cdot\nabla e_{-k}\> \Big] \Big|.
\end{align*}
If they are integrable in $[0,T]$, then the dominated convergence theorem yields
\begin{align*}
  \lim_{N \rightarrow +\infty} 2p & \int_{0}^t  \sum_k  \chi_{N}(|k|)  X_k^{p-1} {\rm Re} \E \Big[ \overline{\widehat{ \bar{\omega}}_k} \big( \<\bar{\omega}, u^1\cdot\nabla e_{-k}\> + \< \omega^2, \bar{u} \cdot\nabla e_{-k}\> \big) \Big] \, \d s \\
  & \leq \int_{0}^t \big(J_{\beta,p}^{(1)} + J_{\beta,p}^{(2)} \big) \, \d s<+\infty, \quad \quad \forall\, t \in [0,T].
\end{align*}
If in addition the nonlinear terms $J_{\beta,p}^{(1)} $ and $ J_{\beta,p}^{(2)}$ verify the following bound analogous to \eqref{eq-nonlinearity 2}: for $i=1,2$,
\begin{equation}\label{eq-nonlin-3}
  J_{\beta,p}^{(i)}  \le \delta (p-1) \sum_k |k|^{2(1-\alpha)} X_k^p + \frac{C_{\alpha,\beta,p,\lambda}}{((p-1) \delta)^{\lambda-1}}  \|\omega_0\|_{L^2}^{\lambda} \sum_{k} X_{k}^p ,
\end{equation}
then Proposition \ref{prop-Y_k^p-diff} and Gr\"onwall's inequality yield the pathwise uniqueness of solutions to 2D stochastic Euler/mSQG equation \eqref{2D Euler/mSQG}. Consequently, the remainder of this subsection is devoted to establishing this nonlinear estimate.

\begin{lemma} \label{lem-nonlinear-1}
  For any $\alpha\in \big(0,\frac{\beta-1}{2} \big)$, $1<p <  \frac{\beta-2\alpha+1}{2}$, $\lambda \geq \frac{4(1-\alpha)}{\beta-2\alpha-1}$, and $\delta>0$, the estimate \eqref{eq-nonlin-3} holds for $i=1$ with the same constant $C_{\alpha,\beta,p,\lambda}>0$ as in \eqref{eq-nonlinearity 2}.
\end{lemma}

\begin{proof}
  Note that $u^1 = \nabla^\perp(-\Delta)^{-\beta/2} \omega^1 = \frac{\rm i}{(2\pi)^{\beta -1}} \sum_{m\ne 0} \frac{m^\perp}{|m|^{\beta}} \widehat{\omega^1}_m\, e_m$, thus
  \begin{equation*}
    \big| \<\bar{\omega}, u^1\cdot\nabla e_{-k}\> \big|
    = (2\pi)^{2-\beta} \Big|\sum_{m\ne 0} \frac{m^\perp\cdot k}{|m|^{\beta}} \widehat{ \bar{\omega}}_{k-m} \, \widehat{\omega^1}_m \Big| \le (2\pi)^{2-\beta} \sum_{m\ne 0} \frac{|m^\perp\cdot k|}{|m|^{\beta}} |\widehat{\bar{\omega}}_{k-m}|\, | \widehat{\omega^1}_m|.
  \end{equation*}
  Due to $\esssup_{\Omega} \sum_{m} | \widehat{\omega^1}_m|^2 \leq \| \omega_0 \|_{L^2}^2$, for any $\beta+\alpha_1+\alpha_2-2>1$, we have
  \begin{equation*}
    \esssup_{\Omega}  \sum_{m} \frac{1}{|m|^{\beta+\alpha_1+\alpha_2-2}} | \widehat{\omega^1}_m| \lesssim \Big( \esssup_{\Omega} \sum_{m} | \widehat{\omega^1}_m|^2 \Big)^{1/2} \leq \| \omega_0 \|_{L^2}.
  \end{equation*}
  Thus, the desired result follows by replacing $Y_{k}$, $\widehat{\omega}_{k}$, $\widehat{\omega}_{k-m}$ and $\widehat{\omega}_{m}$ with $X_{k}$, $\widehat{\bar{\omega}}_{k}$, $\widehat{\bar{\omega}}_{k-m}$ and $\widehat{\omega^1}_{m}$ respectively, in the proof of Proposition \ref{prop-nonlinearity}.
\end{proof}

For the second nonlinear term $J_{\beta,p}^{(2)}$,  we will adopt a method different from the proof of Proposition \ref{prop-nonlinearity} to obtain the desired estimate.

\begin{lemma} \label{lem-nonlinear-2}
  Let $1<\beta\leq2$ and $ 0<\alpha<\frac{\beta+1-\sqrt{(\beta-1)^2+4}}{4}<\frac{\beta-1}{4}$. For any
  \begin{equation}\label{eq-parameter-p}
    \frac{1-\alpha}{1-2\alpha} < p<  \frac{\beta-2\alpha+1}{2}, 
  \end{equation}
   $\lambda \geq \frac{4(1-\alpha)}{\beta-2\alpha-1}$ and $\delta>0$, estimate \eqref{eq-nonlin-3} holds for $i=2$ with a constant $C_{\alpha,\beta,p,\lambda}^{\prime}>0$.
\end{lemma}

\begin{proof}
  We first check that the assumptions on $\alpha$ and $p$ are consistent. Since $0<\alpha<1/2$, the interval $(\frac{1-\alpha}{1-2\alpha},    \frac{\beta-2\alpha+1}{2})$ is nonempty if and only if
  \[
    2(1-\alpha)
    <
    (1-2\alpha)(\beta-2\alpha+1),
  \]
  equivalently,
  \[
    4\alpha^2-2(\beta+1)\alpha+\beta-1>0.
  \]
  The roots of the corresponding quadratic polynomial are
  \[
    \alpha_\pm
    =
    \frac{\beta+1\pm\sqrt{(\beta-1)^2+4}}{4},
  \]
  thus the solution of the above inequality is $\alpha>\alpha_+$ or $\alpha<\alpha_-$. We omit the first range since $\alpha_+ >1$, hence we obtain $\alpha<\alpha_- \in (0, \frac{\beta-1}4)$.

  Note that $\bar{u} = \nabla^\perp(-\Delta)^{-\beta/2} \bar{\omega} = \frac{\rm i}{(2\pi)^{\beta -1}} \sum_{m\ne 0} \frac{m^\perp}{|m|^{\beta}} \widehat{\bar{\omega}}_m e_m$, thus
  \begin{equation*}
    \big| \<\omega^2, \bar{u} \cdot\nabla e_{-k}\> \big|
    = (2\pi)^{2-\beta} \Big|\sum_{m\ne 0} \frac{m^\perp\cdot k}{|m|^{\beta}} \widehat{\omega^2}_{k-m} \, \widehat{\bar{\omega}}_m \Big| \leq (2\pi)^{2-\beta} \sum_{m\ne 0} \frac{|k|}{|m|^{\beta-1}} |\widehat{\omega^2}_{k-m}|\, | \widehat{\bar{\omega}}_m|.
  \end{equation*}
  By the definition of $J_{\beta,p}^{(2)}$, we have
  \begin{align*}
    J_{\beta,p}^{(2)} 
    & \le 2p\sum_k X_k^{p-1} \E \Big( \big|\overline{\widehat{ \bar{\omega}}_k} \big|\, \big| \<\omega^2, \bar{u} \cdot\nabla e_{-k}\> \big|\Big) \\
    & \leq  2p (2\pi)^{2-\beta} \E \Big( \sum_k X_k^{p-1} |\widehat{ \bar{\omega}}_k| \sum_{m\ne 0} \frac{|k|}{|m|^{\beta-1}} |\widehat{\omega^2}_{k-m}|\, | \widehat{\bar{\omega}}_m| \Big) \\
    & = 2p (2\pi)^{2-\beta} \E \Big( \sum_{m\ne 0} \frac{| \widehat{\bar{\omega}}_m|}{|m|^{\beta-1}} \sum_k |k|  X_k^{p-1} |\widehat{ \bar{\omega}}_k|\, |\widehat{\omega^2}_{k-m}| \Big).
  \end{align*}
  For fixed $m \in \Z_0^2$, the Cauchy--Schwarz inequality implies, $\P$-a.s.,
  \begin{equation*}
    \aligned
    \sum_k |k|  X_k^{p-1} |\widehat{ \bar{\omega}}_k|\, |\widehat{\omega^2}_{k-m}|
    &\leq \Big( \sum_{k} |\widehat{\omega^2}_{k-m}|^2 \Big)^{1/2} \Big(  \sum_k |k|^2  X_k^{2(p-1)} |\widehat{ \bar{\omega}}_k|^2\Big)^{1/2} \\
    &\le \|\omega_0 \|_{L^2} \Big(  \sum_k |k|^2  X_k^{2(p-1)} |\widehat{ \bar{\omega}}_k|^2\Big)^{1/2},
    \endaligned
  \end{equation*}
  where we have used $\esssup_{\Omega} \sum_{k} |\widehat{\omega^2}_{k}|^2 \leq \|\omega_0\|_{L^2}^2$. Again by the Cauchy--Schwarz inequality, we obtain
  \begin{align}
    J_{\beta,p}^{(2)} & \leq 2p(2\pi)^{2-\beta} \| \omega_0\|_{L^2}  \bigg[ \E \Big( \sum_{m\ne 0} \frac{| \widehat{\bar{\omega}}_m|}{|m|^{\beta-1}}  \Big)^2 \bigg]^{1/2} \bigg[\E  \sum_k |k|^2  X_k^{2(p-1)} |\widehat{ \bar{\omega}}_k|^2 \bigg]^{1/2} \nonumber \\
                      & = 2p(2\pi)^{2-\beta} \| \omega_0\|_{L^2} \bigg[ \E \Big( \sum_{m\ne 0} \frac{| \widehat{\bar{\omega}}_m|}{|m|^{\beta-1}}  \Big)^2 \bigg]^{1/2} \bigg[ \sum_k |k|^2  X_k^{2p-1} \bigg]^{1/2}.
    \label{eq-J-(2)}
  \end{align}

  Next, we will estimate the last two terms on the right side of the above inequality one by one. For any $\epsilon>0$, since $|m|^{1-\beta}  = |m|^{-1-\epsilon} |m|^{2-\beta+\epsilon}$,  the Cauchy--Schwarz inequality gives
  \begin{align*}
    \E \Big( \sum_{m\ne 0} \frac{| \widehat{\bar{\omega}}_m|}{|m|^{\beta-1}}  \Big)^2 \leq \Big( \sum_{m \in \Z_0^2} |m|^{-2-2\epsilon}\Big) \E \Big(\sum_{m} |m|^{2(2-\beta+\epsilon)} | \widehat{\bar{\omega}}_m|^2 \Big) = C_{\epsilon} \, \sum_{m} |m|^{2(2-\beta+\epsilon)} X_{m},
  \end{align*}
  where $C_{\epsilon}= \sum_{m\in\mathbb Z_0^2}|m|^{-2(1+\epsilon)}<+\infty$. One has $p<\frac{\beta-2\alpha+1}{2} < \frac{2-\alpha}{3-\beta}$, where the second inequality is due to
  \[
    2(2-\alpha) -(\beta-2\alpha+1)(3-\beta) = (\beta-1)^2+2\alpha(2-\beta)>0.
  \]
  Since $ (3-\beta)p<2-\alpha$, for sufficiently small $\epsilon>0$, we have $p< \frac{2-\alpha+ \epsilon}{3-\beta +2\epsilon}$ which implies
  \[
    q_1 :=2p(2-\beta+\epsilon)+2(1+\epsilon)(p-1) <2(1-\alpha).
  \]
  Obviously $q_1>0$. The H\"older  inequality yields
  \begin{align*}
    \E \Big( \sum_{m\ne 0} \frac{| \widehat{\bar{\omega}}_m|}{|m|^{\beta-1}}  \Big)^2
     & \leq C_{\epsilon} \, \sum_{m} |m|^{2(2-\beta+\epsilon) - \frac{q_1}{p}} |m|^{\frac{q_1}{p}} X_{m}                                                                   \\
     & \leq  C_{\epsilon} \, \Big( \sum_{m} |m|^{ \frac{1}{p-1} (2p(2-\beta+\epsilon) - q_1 )} \Big)^{\frac{p-1}{p}} \Big( \sum_{m} |m|^{q_1} X_{m}^p \Big)^{\frac{1}{p}}.
  \end{align*}
  Due to $2p(2-\beta+\epsilon)-q_1= -2(p-1)(1+\epsilon)$ and $C_{\epsilon}= \sum_{m\in\mathbb Z_0^2}|m|^{-2(1+\epsilon)}<+\infty$, we have
  \begin{equation} \label{eq-J-(2)-1}
    \E \Big( \sum_{m\ne 0} \frac{| \widehat{\bar{\omega}}_m|}{|m|^{\beta-1}}  \Big)^2 \leq C_{\epsilon}^{\frac{2p-1}{p}} \Big( \sum_{k} |k|^{q_1} X_{k}^p \Big)^{\frac{1}{p}}.
  \end{equation}

  Now we deal with the last quantity in \eqref{eq-J-(2)}. Set $q_2:=\frac{2p}{2p-1}$; by \eqref{eq-parameter-p}, $p>\frac{1-\alpha}{1-2\alpha}$, thus
  \[
    q_2 <\frac{2\frac{1-\alpha}{1-2\alpha}}{2\frac{1-\alpha}{1-2\alpha}-1}=2(1-\alpha).
  \]
  Since $(2p-1)/p>1$, the Minkowski inequality yields
  \[
    \sum_k|k|^2X_k^{2p-1} = \sum_k \left(|k|^{q_2}X_k^p\right)^{\frac{2p-1}{p}} \leq \Big( \sum_k|k|^{q_2}X_k^p \Big)^{\frac{2p-1}{p}}.
  \]
  Combining this estimate with \eqref{eq-J-(2)} and \eqref{eq-J-(2)-1} gives
  \begin{equation*}
    J_{\beta,p}^{(2)} \leq 2p (2\pi)^{2-\beta} C_{\epsilon}^{\frac{2p-1}{2p}} \| \omega_0\|_{L^2} \Big( \sum_{k} |k|^{q_1} X_{k}^p \Big)^{\frac{1}{2p}} \Big( \sum_k|k|^{q_2}X_k^p \Big)^{\frac{2p-1}{2p}}
  \end{equation*}
  with $q_1,q_2 \in (0, 2(1-\alpha))$. Note that
  \begin{align*}
      \vartheta_p^{\prime}
& =\frac{1}{2p} \frac{q_1}{2(1-\alpha)}+\frac{2p-1}{2p}\frac{q_2}{2(1-\alpha)}
=\frac{1}{2(1-\alpha)}\big((2-\beta+\epsilon)+\frac{(1+\epsilon)(p-1)}{p}+1\big) \\
&=\frac{1}{2(1-\alpha)} \Big( 3-\beta+\epsilon+\frac{(1+\epsilon)(p-1)}{p} \Big)
  \end{align*}
  is the same as $\vartheta_{p}$ given in the proof of Proposition \ref{prop-nonlinearity}. Then, by an argument analogous to the derivation of \eqref{eq-nonlinearity 2} from \eqref{eq-J-3}, which relies on the H\"older and Young inequalities, we obtain that for any $\lambda\geq \frac{4(1-\alpha)}{\beta-2\alpha-1}$ and $\delta>0$, there exists constant $C_{\alpha,\beta,p,\lambda}^{\prime}>0$ such that
  \begin{equation*}
    J_{\beta,p}^{(2)}  \leq  \delta(p-1)  \sum_{k\in\mathbb Z_0^2} |k|^{2(1-\alpha)}X_k^p + \frac{C_{\alpha,\beta,p,\lambda}^{\prime}} {\bigl((p-1)\delta\bigr)^{\lambda-1}} \|\omega_0\|_{L^2}^\lambda \sum_{k\in\mathbb Z_0^2}X_k^p.
  \end{equation*}
  The proof is complete.
\end{proof}

\subsection{Well-posedness of \eqref{2D Euler/mSQG} and sharpness of anomalous regularity}

Thanks to the anomalous regularity due to Kraichnan noise for the solutions to stochastic Navier--Stokes/mSQG equations \eqref{SNSE-Ito}, we can now provide Theorem \ref{thm-Euler/mSQG}.  

\begin{proof}[Proof of Theorem \ref{thm-Euler/mSQG}]
Let $\omega^{(\kappa)}$ be the pathwise unique solution to \eqref{stoch-2D-NS-0} with diffusivity $\kappa>0$. Then by the vanishing-viscosity argument, there exist a sequence $\kappa_n\downarrow0$ and a weak solution $\omega$
  to \eqref{2D Euler/mSQG} such that, $\mathbb P$-a.s.,
  \begin{equation*}
    \< \omega^{(\kappa_{n})}_t, \varphi \> \to \< \omega_t, \varphi \> \quad \text{as } n \rightarrow +\infty,
  \end{equation*}
for all $t\in[0,T]$ and $\varphi\in C^\infty(\T^2)$. Taking $\varphi=e_{k}$, we conclude that all the Fourier coefficients converge. The uniform anomalous regularity (Theorem \ref{thm: anomalous regularity}) and Fatou lemma imply that for any $p \in (1,\frac{\beta-2\alpha+1}{2})$ and $\alpha^{\prime}\in (\alpha, 1)$, the solutions $\omega$ obtained by the vanishing-viscosity method belong to $ L^{2p}([0,T]; \mathcal{X}_{\alpha,p})$ and $L^2\big( \Omega;L^2([0,T];H^{1-\alpha^{\prime}}(\T^2)) \big)$, that is to say, the estimates in \eqref{eq-anomal-regular} hold. The nonincreasing property of the $L^2$-norm $\|\omega_t \|_{L^2}$ follows from the same property of the approximations $\omega^{(\kappa)}$.

  Finally, we prove the pathwise uniqueness. Assume the parameters $\alpha$ and $p$ are given as in \eqref{eq-param-alpha-p}.  Let $\omega^1$ and $\omega^2$ be two solutions to \eqref{2D Euler/mSQG} with the same initial datum $\omega_0\in L^2(\mathbb T^2)$ and satisfying
  \[
    \omega^i\in L^{2p}(0,T;\mathcal X_{\alpha,p}) \cap L^\infty\bigl(\Omega;L^\infty(0,T;L^2(\mathbb T^2))\bigr),\qquad i=1,2.
  \]
  By Proposition \ref{prop-Y_k^p-diff}, Lemma \ref{lem-nonlinear-1} and Lemma \ref{lem-nonlinear-2} with $\delta=\frac{C_0}{2p} \nu$, we get
  \begin{equation*}
    \sum_k X_k^p (t) \leq   \frac{C_{\alpha,\beta,p,\lambda} + C_{\alpha,\beta,p,\lambda}^{\prime}}{\bigl((p-1)\delta\bigr)^{\lambda-1}} \|\omega_0\|_{L^2}^\lambda \int_0^t \sum_{k} X_k^p(s) \, \d s 
  \end{equation*}
  for any $t \in [0,T]$, where $X_k(t)=\E | \langle\omega^1(t)-\omega^2(t),e_{-k}\rangle |^2$. Since $X_k(0)=0$, the Gr\"onwall inequality yields $\sum_k X_k^p(t) \equiv 0$. Hence $X_k(t) =0$ for every $k\in\mathbb Z_0^2$ and $t\in[0,T]$, and therefore $\omega^1=\omega^2$ almost surely. This proves pathwise uniqueness.
\end{proof}

The following proposition shows that, for small initial data, the $L^2\big(\Omega;L^2(\R_+;H^{1-\alpha^{\prime}}) \big)$ regularity of the solution to \eqref{2D Euler/mSQG} with $\alpha^{\prime}>\alpha$ is optimal.

\begin{proposition}\label{prop-optimal}
  Assume $\|\omega_0\|_{L^2} < \frac{\beta-2\alpha-1}{\beta-2\alpha+1} \frac{C_0 \nu}{C_{\alpha,\beta}}$, where the constants $C_0$ and $C_{\alpha,\beta}$ are given in \eqref{def-C-theta} and \eqref{eq-nonlinearity 1} respectively. Then there exists $T_0>0$ such that for all $T>T_0$, the solution $\omega$ to Eq. \eqref{2D Euler/mSQG} satisfies
  \begin{equation*}
    \E \int_{s}^{s+T}  \| \omega_t \|_{H^{1-\alpha}}^2  \, \d t = \infty, \quad \forall\, s\ge 0.
  \end{equation*}
\end{proposition}

\begin{proof}
  Assume $s=0$ for simplicity. Let $\chi_{N}$ be the smooth cutoff function defined above Proposition \ref{prop-Y_k^p-diff}, and set $\widehat{ \omega}_k(t)= \langle\omega_t,e_{-k}\rangle$, $Y_k(t)=\E | \widehat{ \omega}_k(t) |^2$. 
  
  If $\E \int_{0}^{T}  \| \omega_t \|_{H^{1-\alpha}}^2  \, \d t < \infty$, then by applying Lemma \ref{lem-sum-cutoff} with $p=1$ and following the computations analogous to those in the proofs of Propositions \ref{prop-Y_k^p} and \ref{prop-Y_k^p-diff}, we have
  \begin{equation}\label{prop-optimal-1}
    \sum_k Y_k (T)-\sum_k Y_k (0) =  2 \lim_{N \rightarrow \infty} \int_{0}^T  \sum_k  \chi_{N}(|k|) {\rm Re} \E \Big[ \overline{\widehat{ \omega}_k} \< \omega, u \cdot\nabla e_{-k}\>  \Big] \, \d s,
  \end{equation}
  where $u = \nabla^\perp(-\Delta)^{-\beta/2}\omega$. Note that the noise terms do not appear in the case $p=1$. Thus, using $\overline{\widehat{\omega}_{k}}=\widehat{\omega}_{-k}$, it holds
  \begin{equation*}
    \sum_k \chi_{N}(|k|) \overline{\widehat{ \omega}_k} \<\omega,  u \cdot\nabla e_{-k}\>
    = (2\pi)^{2-\beta} \, \i \,  \sum_k \chi_{N}(|k|)  \sum_{m\ne 0} \frac{m^\perp\cdot k}{|m|^{\beta}} \widehat{\omega}_{k-m} \, \widehat{ \omega}_m \, \widehat{ \omega}_{-k}.
  \end{equation*}
  Making the change of variable $k\to m-k$ and using $m^{\perp} \cdot k= - m^{\perp} \cdot (m-k)$, we get
  \begin{equation*}
    \sum_k \chi_{N}(|k|) \overline{\widehat{ \omega}_k} \<\omega,  u \cdot\nabla e_{-k}\>
    = \frac{(2\pi)^{2-\beta}}{2} \, \i \,  \sum_k \sum_{m\ne 0} \big(\chi_{N}(|k|) - \chi_{N}(|m-k|) \big)  \frac{m^\perp\cdot k}{|m|^{\beta}} \widehat{\omega}_{k-m} \, \widehat{ \omega}_m \, \widehat{ \omega}_{-k}.
  \end{equation*}
  Proposition \ref{prop-nonlinearity} (taking $p=1$ in estimate \eqref{eq-nonlinearity 1}) shows that
  \begin{align*}
     & \sup_{N} \Big| \sum_k  \chi_{N}(|k|) {\rm Re} \E \Big[ \overline{\widehat{ \omega}_k} \< \omega, u \cdot\nabla e_{-k}\>  \Big] \Big|  \leq   C_{\alpha,\beta} \|\omega_0\|_{L^2}\, \E \| \omega_t \|_{H^{1-\alpha}}^2,
  \end{align*}
  the latter being integrable on $[0,T]$ under our assumption. Thus, the dominated convergence theorem implies
  \begin{equation*}
    \aligned
    &\lim_{N \to \infty} \int_{0}^T  \sum_k  \chi_{N}(|k|) {\rm Re} \E \Big[ \overline{\widehat{ \omega}_k} \< \omega, u \cdot\nabla e_{-k}\>  \Big] \d s \\
    &= \lim_{N \to \infty} \frac{(2\pi)^{2-\beta}}{2} \i  \int_{0}^T \, \sum_k \sum_{m\ne 0} \big(\chi_{N}(|k|) - \chi_{N}(|m-k|) \big) {\rm Re} \E \Big[ \frac{m^\perp\cdot k}{|m|^{\beta}} \widehat{\omega}_{k-m} \, \widehat{ \omega}_m \, \widehat{ \omega}_{-k} \Big] \d s=0.
    \endaligned
  \end{equation*}
  Consequently, we have $\E \| \omega_t \|_{L^2}^2= \sum_k Y_k (T)=\sum_k Y_k (0)=\| \omega_0 \|_{L^2}^2$.

  However, due to $\|\omega_0\|_{L^2} < \frac{\beta-2\alpha-1}{\beta-2\alpha+1} \frac{C_0 \nu}{C_{\alpha,\beta}}$, Theorem \ref{thm: anomalous regularity}-(i) shows anomalous dissipation for solutions to \eqref{stoch-2D-Euler}. Since this property is independent of $\kappa$, similarly to the proof of Theorem \ref{thm-Euler/mSQG}, the vanishing-viscosity argument implies that there is a constant $C_{\alpha,p,\mu_{p}}>0$ such that
  \begin{equation*}
    \E \| \omega_t \|_{L^2}^2 \leq C_{\alpha,p,\mu_{p}} e^{- \mu_{p} T} \|\omega_0\|_{L^2}^2.
  \end{equation*}
  There exists a $T_0>0$ such that $C_{\alpha,p,\mu_{p}} e^{- \mu_{p} T_0}=1$, thus for all $T>T_0$, we get a contradiction. As a result, we conclude that $\E \int_{0}^{T}  \| \omega_t \|_{H^{1-\alpha}}^2  \, \d t = \infty$.
\end{proof}

\section{Scaling limits} \label{sec: scaling limit}

The purpose of this section is to prove Theorem \ref{thm: scal_Sobolev}. We consider a sequence of the stochastic 2D Euler/mSQG equations \eqref{2D Euler/mSQG-n} driven by Kraichnan noises without low frequencies, and aim to show that they converge in 
  $$L^2\big(\Omega; L^2([0,T]; H^{1-\alpha-}(\T^2)) \cap L^{\infty}([0,T]; L^2(\T^2)) \big) $$ 
to the deterministic viscous 2D Navier--Stokes/mSQG equation \eqref{2D determined NSE/mSQG}. The proof relies on the scaling limit in $L^2\big(\Omega; C([0,T]; H^{-\epsilon}(\T^2)) \big)$ (c.f. \cite[Section 3.3]{Flandoli_2024_Quantitative}) and on the uniform anomalous regularity over $\cH_{\alpha}$ with $\sup_{k}|\theta_{n,k}| \rightarrow 0$, established in Theorem \ref{thm-Euler/mSQG}.

Let $Y_{n,k}=  \E | \widehat{\omega}^n_{k} |^2$. Due to the anomalous regularity proved in Theorem \ref{thm-Euler/mSQG}, we have 
  $$ \int_0^T \sum_{k \in \Z_0^2} |k|^{2(1-\alpha)} Y_{n,k}^p (t)\,\d t =\int_{0}^T  \| \omega^n_t \|_{\mathcal{X}_{\alpha,p}}^{2p} \, \d t<+\infty.$$ 
Hence, we can apply the discrete Poincar\'e inequality \eqref{def-C-theta} in Section \ref{subsec: noise}: there exists a uniform constant $C_0>0$ such that, for every $n\ge 1$ and $p>1$,
\begin{equation} \label{eq: weighted-poincare-n}
  8 \pi^2 p \sum_{k,l} \theta_{n,l}^2  |\Pi_{l^\perp} k|^2  \big(Y_{n,k+l}^{p-1} -Y_{n,k}^{p-1}\big) (Y_{n,k+l} -Y_{n,k}) \geq C_0 \frac{p-1}{p}  \sum_{k} |k|^{2(1-\alpha)} Y_{n,k}^p .
\end{equation}

We first recall the known scaling limit result in the weaker space $L^2\big(\Omega; C([0,T]; H^{-\epsilon}(\T^2)) \big)$ from \cite[Section 3.3]{Flandoli_2024_Quantitative}. 
Since our noise coefficients satisfy all the assumptions of \cite{Flandoli_2024_Quantitative} and $\|\theta_n \|_{\ell^\infty} \sim \sqrt{\alpha}\, n^{-1}$, we have the following convergence rate of $\omega^n$ to $\widetilde{\omega}$ in $L^2 \big(\Omega; C([0,T]; H^{-\epsilon}(\T^2)) \big)$.

\begin{proposition} \label{prop: scal_lim_negative}
Fix $\nu > 0$, $\beta \in (1,2]$, $q_0 > \frac{2}{\beta-1}$ and $s \in (0,\beta-1)$. Let $\omega^n$ and $\widetilde{\omega}$ be the unique solution of \eqref{2D Euler/mSQG-n} and of \eqref{2D determined NSE/mSQG}, with common initial datum $\omega_0$. Then for any $T < +\infty$, $p \geq 1$ and $\epsilon \in (0,s)$, there exists a constant $C = C(T, \nu, s, \beta,p,q_0,\epsilon)>0$ such that
  \begin{equation*}
    \E \Big[ \sup_{t \in [0, T]} \|\omega_t^n - \tilde{\omega}_t\|_{H^{-s}}^p \Big]^{1/p}
    \leq C  \exp \left( C \|\omega_0\|_{L^2}^{q_0} \right) \frac{\|\omega_0\|_{L^2}}{n^{s-\epsilon}}
    .
  \end{equation*}
\end{proposition}

The remainder of this section is devoted to two tasks. We first prove the uniform anomalous regularity with respect to the noise coefficients $\{\theta_{n,k} \}_{k\in\mathbb Z_0^2}$, and use interpolation to establish the scaling limit in $L^2 \big(\Omega; L^2([0,T]; H^{1-\alpha-}(\T^2)) \big)$ from the stochastic 2D Euler/mSQG equation \eqref{2D Euler/mSQG-n} to the limit equation \eqref{2D determined NSE/mSQG}. Next, we prove the scaling limit in $L^2\big(\Omega; L^{\infty}([0,T]; L^2(\T^2)) \big)$ and provide a quantitative estimate that depends on the initial data.

\subsection{Scaling limit in $L^2([0,T]; H^{1-\alpha-}(\T^2))$}

To show the scaling limit in $L^2([0,T]; H^{1-\alpha-}(\T^2))$, we need the uniform anomalous regularity for the solutions $\omega^n$ of \eqref{2D Euler/mSQG-n} which can be deduced from \eqref{eq: weighted-poincare-n} and Theorem \ref{thm-Euler/mSQG}.

\begin{lemma} \label{lem: uniform-anomalous}
Assume $\beta\in (1,2]$ and $\alpha\in \big(0, \frac{1}{4} (\beta+1-\sqrt{(\beta-1)^2+4}\, ) \big)$. Let $\omega^n$ be the unique solution of the stochastic 2D Euler/mSQG equation \eqref{2D Euler/mSQG-n} with initial value $\omega_0$. Then $\omega^n$ satisfies the uniform anomalous regularity: for any $p \in \big(1,\frac{\beta-2\alpha+1}{2} \big)$ and $\tilde{\alpha}\in (\alpha, 1)$, there is $C=C(\alpha,\beta,p,\nu)>0$ such that
  \begin{equation} \label{eq: uniform-Fourier}
    \sup_{n} \int_{0}^T  \| \omega^n_t \|_{\mathcal{X}_{\alpha,p}}^{2p} \, \d t \leq  C \Big(1+ T \|\omega_0\|_{L^2}^{\frac{4(1-\alpha)}{\beta-2\alpha-1}} \Big) \|\omega_0\|_{L^2}^{2p},
  \end{equation}
and the map $t \mapsto \|\omega^n_t\|_{L^2}$ is $\P$-a.s. nonincreasing. Consequently, for any $\tilde{\alpha}\in (\alpha, 1)$ and subinterval $I \subset [0,T]$, there is a constant $C=C(\alpha,\tilde{\alpha},\beta,\nu,T)>0$ such that
  \begin{equation} \label{eq: uniform-Sobolev}
    \sup_{n} \E \int_I \|\omega^n_t\|_{H^{1-\tilde{\alpha}}}^2\, \d t \leq C |I|^{\eta(\alpha,\tilde{\alpha})} \Big(1+ \|\omega_0 \|_{L^2}^{\frac{4(1-\alpha)}{\beta-2\alpha-1}} \Big)  \|\omega_0\|_{L^2}^{2},
  \end{equation}
where $\eta(\alpha,\tilde{\alpha}):=\min \big\{\frac{\tilde{\alpha} -\alpha}{4-\alpha-\tilde{\alpha}}, \frac{\beta-2\alpha-1}{\beta-2\alpha+3}  \big\}$.
\end{lemma}

\begin{proof}
  Using the uniform weighted Poincar\'e inequality \eqref{eq: weighted-poincare-n} and following the proof of Theorem \ref{thm-Euler/mSQG}, we apply the intermediate result Theorem \ref{thm: anomalous regularity}-(ii) with $\lambda=\frac{4(1-\alpha)}{\beta-2\alpha-1}$ to obtain the uniform anomalous regularity \eqref{eq: uniform-Fourier}. By Proposition \ref{prop:X-Sobolev}-(2) and H\"older's inequality, we have
  \begin{align*}
    \sup_{n} \E \int_I \|\omega^n_t\|_{H^{1-\tilde{\alpha}}}^2\, \d t
     & \lesssim_{\alpha,\tilde{\alpha},p} |I|^{1-\frac{1}{p}} \Big( \sup_{n} \int_I  \| \omega^n_t \|_{\mathcal{X}_{\alpha,p}}^{2p} \, \d t \Big)^{\frac{1}{p}}
    \\ & \lesssim_{\alpha,\tilde{\alpha},p} |I|^{\eta(\alpha,\tilde{\alpha})} \Big(1+\|\omega_0\|_{L^2}^{\frac{4(1-\alpha)}{\beta-2\alpha-1}} \Big)  \|\omega_0\|_{L^2}^{2},
  \end{align*}
where $p =\min \big\{ \frac{4-\alpha-\tilde{\alpha}}{2(2-\tilde{\alpha})}, \frac{\beta-2\alpha+3}{4} \big\}< \min \big\{ \frac{2-\alpha}{2-\tilde{\alpha}}, \frac{\beta-2\alpha+1}{2} \big\}$ and $\lambda=\frac{4(1-\alpha)}{\beta-2\alpha-1} $. It gives the uniform estimate \eqref{eq: uniform-Sobolev}. The nonincreasing property of the $L^2$-norm has been proved in Theorem \ref{thm-Euler/mSQG}.
\end{proof}

Combining the estimate in Proposition \ref{prop: scal_lim_negative} and the uniform anomalous regularity \eqref{eq: uniform-Sobolev}, we can establish the scaling limit in $L^2 \big(\Omega; L^2([0,T]; H^{1-\alpha-}(\T^2)) \big)$.

\begin{proof}[Proof of Theorem \ref{thm: scal_Sobolev}-(a)]
  Using the interpolation of Sobolev spaces and H\"older's inequality, for any $0<s<\beta-1$ and $\tilde{\alpha} \in (\alpha,\alpha^{\prime}) \subset (\alpha,1)$, we have
  \begin{equation}\label{eq: scal_sobolev.1}
  \aligned
    & \E \int_I \|\omega_t^n - \widetilde{\omega}_t\|_{H^{1-\alpha^{\prime}}}^2 \, \d t
     \le \E \int_I \|\omega_t^n - \widetilde{\omega}_t\|_{H^{-s}}^{\frac{2(\alpha^{\prime}-\tilde{\alpha})}{1-\tilde{\alpha} +s}} \|\omega_t^n - \widetilde{\omega}_t\|_{H^{1-\tilde{\alpha}}}^{\frac{2(1-\alpha^{\prime}+s)}{1-\tilde{\alpha} +s}} \, \d t \\
     &\qquad \le \bigg[\E \int_I \|\omega_t^n - \widetilde{\omega}_t\|_{H^{-s}}^2 \, \d t\bigg]^{\frac{\alpha^{\prime}-\tilde{\alpha}}{1-\tilde{\alpha} +s}} \bigg[\E \int_I \|\omega_t^n - \widetilde{\omega}_t\|_{H^{1-\tilde{\alpha}}}^2 \, \d t\bigg]^{\frac{1-\alpha^{\prime}+s}{1-\tilde{\alpha} +s}}.
  \endaligned
  \end{equation}
Applying Proposition \ref{prop: scal_lim_negative} with $p=2$ and $q_0=\frac{4(1-\alpha)}{\beta-2\alpha-1}>\frac{2}{\beta-1}$, for any $0<\epsilon<s$, we get
  \begin{equation}\label{eq: scal_sobolev.2}
    \aligned
    \E \int_I \|\omega_t^n - \widetilde{\omega}_t\|_{H^{-s}}^2 \, \d t 
    & \le |I|\, \E\Big(\sup_{t\in [0,T]} \|\omega_t^n - \widetilde{\omega}_t\|_{H^{-s}}^2 \Big) \\
    & \le C |I| \exp \Big( C \|\omega_0 \|_{L^2}^{\frac{4(1-\alpha)}{\beta-2\alpha-1}} \Big) \frac{\|\omega_0 \|_{L^2}^2}{n^{2(s-\epsilon)}} .
    \endaligned
  \end{equation}
Next, using the uniform anomalous regularity \eqref{eq: uniform-Sobolev}, we have
  \begin{equation*}
    \sup_{n} \E \int_I \|\omega^n_t\|_{H^{1-\tilde{\alpha}}}^2\, \d t
    \le C |I|^{\eta(\alpha,\tilde{\alpha})} \Big(1+ \|\omega_0 \|_{L^2}^{\frac{4(1-\alpha)}{\beta-2\alpha-1}} \Big)  \|\omega_0\|_{L^2}^{2} .
  \end{equation*}
  By the Sobolev interpolation,  H\"older's inequality, and energy estimate of \eqref{2D determined NSE/mSQG}, we obtain
  \begin{align*}
    \int_I \|\widetilde{\omega}_t \|_{H^{1-\tilde{\alpha}}}^2\,\d t
     & \le \int_I \|\widetilde{\omega}_t \|_{L^2}^{2 \tilde{\alpha}} \|\widetilde{\omega}_t \|_{H^1}^{2(1-\tilde{\alpha} )} \, \d t                                            \\
     & \le \bigg[\int_I \|\widetilde{\omega}_t \|_{L^2}^{2} \,\d t\bigg]^{\tilde{\alpha}} \bigg[ \int_I \|\widetilde{\omega}_t \|_{H^1}^{2} \, \d t \bigg]^{1-\tilde{\alpha} } \\
     & \lesssim_\nu |I|^{\tilde{\alpha}} \|\omega_0 \|_{L^2}^{2}.
  \end{align*}
Using the triangle inequality and the fact $\tilde{\alpha}>\frac{\tilde{\alpha}-\alpha}{4-\alpha-\tilde{\alpha}} \geq \eta(\alpha,\tilde{\alpha})$, we get
  \begin{equation} \label{eq: scal_sobolev.3}
  \aligned
    \E \int_I \|\omega_t^n - \widetilde{\omega}_t\|_{H^{1-\tilde{\alpha}}}^2 \, \d t
     & \le 2\, \E \int_I \|\omega_t^n \|_{H^{1-\tilde{\alpha}}}^2 \, \d t + 2\,\E \int_I \|\widetilde{\omega}_t\|_{H^{1-\tilde{\alpha}}}^2 \, \d t       \\
     & \le C |I|^{\eta(\alpha,\tilde{\alpha})} \Big( 1+\|\omega_0\|_{L^2}^{\frac{4(1-\alpha)}{\beta-2\alpha-1}} \Big) \|\omega_0 \|_{L^2}^2          \\
     & \le C |I|^{\eta(\alpha,\tilde{\alpha})} \exp \left( C \|\omega_0\|_{L^2}^{\frac{4(1-\alpha)}{\beta-2\alpha-1}} \right)  \|\omega_0\|_{L^2}^2.
  \endaligned
  \end{equation}

Due to $s<\beta-1\le1$ and $\tilde{\alpha}<\alpha^{\prime} \leq 1$, we have $\frac{1}{1-\tilde{\alpha}+s} \geq \frac{1}{1+s} \geq \frac{1}{2}$ and $\frac{1-\alpha^{\prime}+s}{1-\alpha+s} \geq \frac{s}{1+s}$. Substituting the above estimates \eqref{eq: scal_sobolev.2} and \eqref{eq: scal_sobolev.3} into \eqref{eq: scal_sobolev.1}, and taking $\tilde{\alpha}=\frac{\alpha^{\prime}+3\alpha}{4}$ (so $\alpha^{\prime}-\tilde{\alpha}=\frac{3(\alpha^{\prime}-\alpha)}{4}$), then we have
  \begin{align*}
    \E \int_I \|\omega_t^n - \widetilde{\omega}_t\|_{H^{1-\alpha^{\prime}}}^2 \, \d t
     & \le C \exp \left( C \|\omega_0\|_{L^2}^{\frac{4(1-\alpha)}{\beta-2\alpha-1}} \right)  \|\omega_0\|_{L^2}^2  \frac{|I|^{\frac{\alpha^{\prime}-\tilde{\alpha}}{1-\tilde{\alpha} +s}+\eta(\alpha,\tilde{\alpha})\frac{1-\alpha^{\prime}+s}{1-\tilde{\alpha} +s}}}{n^{2(s-\epsilon)\frac{\alpha^{\prime}-\tilde{\alpha}}{1-\tilde{\alpha} +s}}} \\
     & \le C \exp \left( C \|\omega_0\|_{L^2}^{\frac{4(1-\alpha)}{\beta-2\alpha-1}} \right)  \|\omega_0\|_{L^2}^2  \frac{|I|^{(\frac{3}{8}+\frac{1}{16}\frac{s}{1+s})(\alpha^{\prime}-\alpha)}}{n^{\frac{3}{2} \frac{s-\epsilon}{1+s} (\alpha^{\prime}-\alpha)}}
  \end{align*}
  Taking $\epsilon>0$ and $s = \frac{s_0+\epsilon}{1-s_0} \in (0,\beta-1)$ (i.e. $s_0 =\frac{s-\epsilon}{1
      +s} \in (0,1-\frac{1}{\beta})$), we obtain \eqref{eq: scal_sobolev}.
\end{proof}

\subsection{Scaling limit in $L^{\infty}([0,T]; L^2(\T^2))$} \label{subsec-C-L^2}

To establish the scaling limit in $L^2 \big(\Omega;L^{\infty}([0,T]; L^2(\T^2)) \big)$, we will exploit the nonincreasing property of the $L^2$-norm for the solution $\omega^n$ to \eqref{2D Euler/mSQG-n}. This allows us to upgrade the known scaling limit in $L^2(\Omega; C([0,T]; H^{-\epsilon}(\T^2)))$ to the stronger limit in $L^2\big(\Omega; L^{\infty}([0,T]; L^2(\T^2)) \big)$.

\begin{proof}[Proof of Theorem \ref{thm: scal_Sobolev}-(b)]
  Let $E_n(t)=\|\omega^n_t\|_{L^2}^2$ and $E(t)=\|\widetilde{\omega}_t\|_{L^2}^2$, then we have
  \begin{equation*}
    \| \omega^n_t - \widetilde{\omega}_t \|_{L^2}^2 = E_n(t) - E(t) -2 \< \omega^n_t-\widetilde{\omega}_t, \widetilde{\omega}_t \>.
  \end{equation*}
  Let $R>0$ and $P_{\le R}$ be the projection onto the subspace of functions with frequency less than or equal to $R$. Then, for any $t \in [0,T]$ and $s \in (0, \beta-1)$, we have
  \begin{align*}
    | \< \omega^n_t-\widetilde{\omega}_t, \widetilde{\omega}_t \>|
     & \le  | \< \omega^n_t-\widetilde{\omega}_t, P_{\le R} \widetilde{\omega}_t \>| + | \< \omega^n_t-\widetilde{\omega}_t, (I-P_{\le R})  \widetilde{\omega}_t \>|                            \\
     & \le R^{s}  \| \widetilde{\omega}_t \|_{L^2} \| \omega^n_t-\widetilde{\omega}_t \|_{H^{-s}} + \| \omega^n_t-\widetilde{\omega}_t \|_{L^2} \|  (I-P_{\le R}) \widetilde{\omega}_t \|_{L^2} \\
     & \le R^{s} \| \omega_0 \|_{L^2} \sup_{t \in [0,T]} \| \omega^n_t-\widetilde{\omega}_t \|_{H^{-s}}  + 2 \| \omega_0 \|_{L^2} \|  (I-P_{\le R}) \widetilde{\omega}_t \|_{L^2},
  \end{align*}
  where in the last step we used $\| \omega^n_t \|_{L^2} \le \| \omega_0 \|_{L^2}$ and $\| \widetilde{\omega}_t \|_{L^2} \le \| \omega_0 \|_{L^2}$ for all $t\in[0,T]$.

Using the nonincreasing property of the $L^2$-norm of $\omega^n$, for any $h>0$ and $t \in [h,T]$ we have
  \begin{equation*}
    E_n(t) - E(t) \le \frac{1}{h} \int_{t-h}^{t} \big( E_n(\tau) - E(\tau) \big) \d \tau + \sup_{\tau\in[t-h,t]} |E(\tau) - E(t)|.
  \end{equation*}
By the Cauchy inequality, we have
  \begin{align*}
     & \frac{1}{h} \int_{t-h}^{t} \big( E_n(\tau) - E(\tau) \big) \d \tau = \frac{1}{h} \int_{t-h}^{t} \<\omega^n_\tau +\widetilde{\omega}_\tau, \omega^n_\tau -\widetilde{\omega}_\tau \> \d \tau \\
     & \le \frac{1}{h} \int_{t-h}^{t} \big( \| \omega^n_{\tau} \|_{H^s} + \| \widetilde{\omega}_{\tau} \|_{H^s} \big) \| \omega^n_\tau - \widetilde{\omega}_\tau \|_{H^{-s}} \d \tau \\
     & \le \frac{1}{h^{1/2}} \bigg[ \Big(\int_{t-h}^{t} \|\omega^n_\tau\|_{H^{s}}^2 \d \tau \Big)^{1/2} +\Big(\int_{t-h}^{t} \|\widetilde{\omega}_{\tau} \|_{H^{s}}^2 \d \tau \Big)^{1/2} \bigg] \sup_{t \in [0,T]} \| \omega^n_t-\widetilde{\omega}_t \|_{H^{-s}} .
  \end{align*}
For $s \in (0,1-\alpha)$, applying \eqref{eq: uniform-Sobolev} with $\tilde{\alpha}=1-s$ and using the energy estimate of $\widetilde{\omega}_\tau$, we obtain
  \begin{equation*}
    \E \bigg[\sup_{t\in [h,T]} \frac{1}{h} \int_{t-h}^{t} \big( E_n(\tau) - E(\tau) \big) \d \tau \bigg] \lesssim_{\nu,T} \Big(1+ \|\omega_0 \|_{L^2}^{\frac{2(1-\alpha)}{\beta-2\alpha-1}} \Big)  \frac{\|\omega_0\|_{L^2}}{h^{\frac{1-\eta(\alpha,1-s)}{2}}}  \E\Big[ \sup_{t \in [0,T]} \| \omega^n_t-\widetilde{\omega}_t \|_{H^{-s}}^2  \Big]^{\frac12},
  \end{equation*}
where $\eta(\alpha,1-s)$ is given by Lemma \ref{lem: uniform-anomalous}. For $t \in [0,h]$, using the nonincreasing property of $t \mapsto \| \omega_t \|_{L^2}^2$, we have
  \begin{equation*}
    E_n(t) - E(t) =  E_n(t) -\| \omega_0\|_{L^2}^2 + \| \omega_0\|_{L^2}^2 - E(t) \le E(0)-E(t).
  \end{equation*}

Combining the above estimates, we obtain that for any $s \in (0,\beta-1) \cap (0,1-\alpha)$,
  \begin{align}
     & \E \Big[ \sup_{t \in [0,T]} \| \omega^n_t - \widetilde{\omega}_t \|_{L^2}^2 \Big] \nonumber \\
     & \lesssim 2 R^{s}  \| \omega_0 \|_{L^2}  \E\Big[ \sup_{t \in [0,T]} \| \omega^n_t-\widetilde{\omega}_t \|_{H^{-s}}  \Big] + 4  \|\omega_0 \|_{L^2} \sup_{t\in[0,T]}  \|  (I-P_{\le R})\widetilde{\omega}_t\|_{L^2} \nonumber  \\
     & \quad +  \frac{\|\omega_0\|_{L^2}}{h^{\frac{1-\eta(\alpha,1-s)}{2}}} \Big(1+ \|\omega_0 \|_{L^2}^{\frac{2(1-\alpha)}{\beta-2\alpha-1}} \Big)  \E\Big[ \sup_{t \in [0,T]} \| \omega^n_t-\widetilde{\omega}_t \|_{H^{-s}}^2  \Big]^{1/2} + \sup_{ \substack{t,s  \in [0,T] \\ |t-s| \le h } } |E(t) - E(s)|. \label{est: convergence_L2}
  \end{align}
Applying Proposition \ref{prop: scal_lim_negative}, the first and third terms on the right-hand side converge to zero as $n\to+\infty$. The second term vanishes as $R\to+\infty$, and the last term converges to zero as $h\to0$ due to the continuity of $E(t)$. Thus, we have
  \begin{equation*}
    \lim_{n\to+\infty} \E\Big[ \sup_{t\in[0,T]}\|\omega_t^n-\widetilde{\omega}_t\|_{L^2}^2 \Big] =0,
  \end{equation*}
which is the scaling limit in $L^2 \big(\Omega; L^{\infty}([0,T]; L^2(\T^2)) \big)$.
\end{proof}


As mentioned below Theorem \ref{thm: scal_Sobolev}, assertion (b) only gives us a qualitative convergence result. The remainder part of this section is devoted to deriving a certain quantitative rate of convergence in $L^2 \big(\Omega; L^{\infty}([0,T]; L^2(\T^2)) \big)$. To this end, we will invoke Lemma \ref{lem: convergence rate pre} to estimate the terms in \eqref{est: convergence_L2}. First we need the following preparation which provides an estimate on the convolution of the nonlinear term.


\begin{lemma} \label{lem: determined nonlinear}
  For any $0\le a<b \le T$, $0 \leq \gamma<\beta-1$ and $0<\epsilon < \beta-1-\gamma$, the solution $\widetilde{\omega}_t$ of 2D NS/mSQG Eq. \eqref{2D determined NSE/mSQG} with initial datum $\omega_0 \in L^2(\T^2)$ satisfies
  \begin{equation*}
    \Big\| \int_a^b e^{\nu(b-r)\Delta} ( \widetilde{u}_r \cdot \nabla\widetilde{\omega}_r) \d r \Big\|_{H^{\gamma}}^2  \lesssim_{\beta,\epsilon} \frac{\| \omega_0 \|_{L^2}^{4}}{\nu^{3-\beta +\gamma+\epsilon}} (b-a)^{\beta-1-\gamma-\epsilon} ,
  \end{equation*}
  where $\widetilde{u} = \nabla^\perp(-\Delta)^{-\beta/2} \widetilde{\omega}$.
\end{lemma}

\begin{proof}
Using the classical semigroup property in Lemma \ref{lem-heat-semigroup-2}, it holds
  \begin{equation} \label{lem: determined nonlinear.1}
    \Big\| \int_a^b e^{\nu(b-r)\Delta} (\widetilde{u}_r \cdot\nabla\widetilde{\omega}_r) \d r \Big\|_{H^{\gamma}}^2 \lesssim\frac{1}{\nu} \int_a^b \| \widetilde{u}_r \cdot\nabla \widetilde{\omega}_r \|_{H^{\gamma-1}}^2 \d r.
  \end{equation}
By the classical estimate on the Sobolev norm of products (see Lemma \ref{lem-product-Sobolev}), we have
  \begin{equation} \label{lem: determined nonlinear.2}
  \| \widetilde{u}_r\cdot \nabla\widetilde{\omega}_r \|_{H^{\gamma-1}} 
  \lesssim \| \widetilde{u}_r\, \widetilde{\omega}_r \|_{H^\gamma} 
  \lesssim \| \widetilde{u}_r \|_{H^{\beta-1 - \epsilon}} \| \widetilde{\omega}_r \|_{H^{2-\beta+ \gamma+ \epsilon}}. 
  \end{equation}  
Note that $\| \widetilde{u}_r \|_{H^{\beta-1 - \epsilon}} \lesssim \| \widetilde{\omega}_r \|_{H^{- \epsilon}} \lesssim \| \widetilde{\omega}_r \|_{L^2}$. Due to our choice of parameters, it holds $2-\beta+ \gamma+ \epsilon < 1$, thus by interpolation,
  $$\| \widetilde{u}_r\cdot \nabla\widetilde{\omega}_r \|_{H^{\gamma-1}} 
  \le \| \widetilde{\omega}_r \|_{L^2} \| \widetilde{\omega}_r \|_{L^2}^{\beta-1-\gamma -\epsilon} \| \widetilde{\omega}_r \|_{H^1}^{2-\beta+ \gamma+ \epsilon} 
  \le \|\omega_0 \|_{L^2}^{\beta -\gamma -\epsilon} \| \widetilde{\omega}_r \|_{H^1}^{2-\beta+ \gamma+ \epsilon}, $$
where in the last step we have used the energy estimate of $\widetilde{\omega}_r$. Plugging this estimate into \eqref{lem: determined nonlinear.1}, we arrive at
  $$\aligned
  \Big\| \int_a^b e^{\nu(b-r)\Delta}  \widetilde{u}_r\cdot\nabla\widetilde{\omega}_r\d r \Big\|_{H^{\gamma}}^2 
  & \lesssim \frac{\|\omega_0 \|_{L^2}^{2(\beta -\gamma -\epsilon)}}{\nu} \int_a^b \| \widetilde{\omega}_r \|_{H^1}^{2(2-\beta+ \gamma+ \epsilon)} \,\d r \\
  &\le \frac{\|\omega_0 \|_{L^2}^{2(\beta -\gamma -\epsilon)}}{\nu} (b-a)^{\beta-1 -\gamma -\epsilon} \Big( \int_a^b \| \widetilde{\omega}_r \|_{H^1}^{2} \,\d r \Big)^{2-\beta+ \gamma+ \epsilon}.
  \endaligned $$
Using again the energy estimate for the solution $\widetilde{\omega}_r$ of \eqref{2D determined NSE/mSQG}, we obtain the desired result. 
\end{proof}

Before moving to Lemma \ref{lem: convergence rate pre}, we denote the spectral modulus of the initial datum by
\begin{equation}\label{def: initial-modulus}
  M_{\omega_0}(R) := \bigg( \sum_{k\ne 0}\min\Big\{1,\frac{|k|^2}{R^2}\Big\} |\widehat{\omega}_{0,k}|^2 \bigg)^{1/2}.
\end{equation}
As $R \rightarrow +\infty$, the dominated convergence theorem yields $M_{\omega_0}(R) \rightarrow 0$, but it does not provide a uniform rate over all $\omega_0 \in  L^2(\T^2) $ with $\|\omega_0\|_{L^2} \le M$. For $s>0$, if the initial datum $\omega_0$ satisfies $\|\omega_0\|_{H^s} \le M$, then the convergence rate is uniform and satisfies $M_{\omega_0}(R) \le \frac{M}{R^{s}}$ for all $R \geq 1$.

\begin{lemma} \label{lem: convergence rate pre}
  For any $h>0$, $R \geq 1$ and $0 \leq \gamma<\beta-1$, the solution $\widetilde{\omega}_t$ of 2D NS/mSQG Eq. \eqref{2D determined NSE/mSQG} with initial datum $\omega_0 \in L^2(\T^2)$ satisfies
  \begin{align}
    \sup_{\substack{t,s\in[0,T] \\|t-s|\le h}}
    \big| \|\widetilde{\omega}_t\|_{L^2}^2 -\|\widetilde{\omega}_s\|_{L^2}^2 \big|
    & \lesssim_{\beta,\gamma} \frac{(\nu h)^{\gamma/2}}{\nu^{\frac{5-\beta + 3\gamma}4}} \| \omega_0 \|_{L^2}^{3}+M_{\omega_0}\big((\nu h)^{-1/2}\big) \| \omega_0 \|_{L^2},   \label{eq:limit-time-modulus} \\
    \sup_{t\in[0,T]} \| (I-P_{\le R}) \widetilde{\omega}_t\|_{L^2}
    & \lesssim_{\beta,\gamma} M_{\omega_0}(R)+\frac{T^{\frac{\beta-1-\gamma}{4}}}{R^{\gamma} \nu^{\frac{5-\beta + \gamma}4}} \| \omega_0 \|_{L^2}^2. \label{eq:limit-tail}
  \end{align}
\end{lemma}

\begin{proof}
Write the mild solution of \eqref{2D determined NSE/mSQG} as
  \begin{equation}\label{eq:mild-split}
    \widetilde{\omega}_t=e^{\nu t\Delta}\omega_0+N(t),
    \qquad
    N(t)=-\int_0^t e^{\nu(t-r)\Delta}  (\widetilde{u}_r \cdot\nabla \widetilde{\omega}_r) \d r.
  \end{equation}
For any $b \geq a \geq 0$ and $b-a \leq h$, we have
  \begin{align*}
    \widetilde{\omega}_b-\widetilde{\omega}_a
     & =(e^{\nu (b-a)\Delta}-I)e^{\nu a\Delta}\omega_0 +(e^{\nu (b-a)\Delta}-I)N(a) \\
     & \quad- \int_a^b e^{\nu(b-r)\Delta} (\widetilde{u}_r \cdot\nabla \widetilde\omega_r) \d r.
  \end{align*}
For the first term, Parseval's identity and $|1-e^{-x}|^2 \le \min\{1,x^2\} \le \min\{1,x\}$ imply
  \begin{align*}
    \|(e^{\nu (b-a)\Delta}-I)e^{\nu a\Delta}\omega_0\|_{L^2}^2 
    & \le \sum_{k\ne0} \big| 1-e^{- 4 \pi^2 \nu (b-a)|k|^2}\big|^2 |\widehat{\omega}_{0,k}|^2 \\
    & \le \sum_{k\ne0} \min\{1,4 \pi^2 \nu h|k|^2\}    |\widehat{\omega}_{0,k}|^2  \leq 4\pi^2 M_{\omega_0}((\nu h)^{-1/2})^2,
  \end{align*}
for all $b-a \leq h$. For the second term, by the classical heat semigroup estimate (cf. Lemma \ref{lem-heat-semigroup-1}), applying Lemma \ref{lem: determined nonlinear} with $\epsilon=\frac12 (\beta-1 -\gamma)$ yields
  \[
    \|(e^{\nu (b-a)\Delta}-I)N(a)\|_{L^2}^2 \lesssim (b-a) ^{\gamma}\|N(a)\|_{H^\gamma}^2 \lesssim \frac{h^{\frac{\beta-1+ \gamma}2}}{\nu^{\frac{5-\beta +\gamma}2}} \| \omega_0 \|_{L^2}^{4} 
  \]
for all $0 \leq \gamma<\beta-1$. For the third term, applying Lemma \ref{lem: determined nonlinear} with $\gamma=0$ and $0 < \epsilon < \beta-1$, we obtain 
  \begin{equation*}
    \Big\| \int_a^b e^{\nu(b-r)\Delta}  (\widetilde{u}_r \cdot\nabla \widetilde{\omega}_r) \d r \Big\|_{L^2}^2  \lesssim_{\beta,\epsilon} \frac{h^{\beta-1-\epsilon}}{\nu^{3-\beta+\epsilon}} \| \omega_0 \|_{L^2}^{4} . 
  \end{equation*}
For any $\gamma\in [0, \beta-1)$, taking $\epsilon= \frac12(\beta-1-\gamma)$ in the above inequality gives us
  \begin{equation*}
    \Big\| \int_a^b e^{\nu(b-r)\Delta}  (\widetilde{u}_r \cdot\nabla \widetilde{\omega}_r) \d r \Big\|_{L^2}^2  \lesssim \frac{h^{\frac{\beta-1+ \gamma}2}}{\nu^{\frac{5-\beta +\gamma}2}} \| \omega_0 \|_{L^2}^{4} 
    \lesssim_T \frac{h^\gamma}{\nu^{\frac{5-\beta +\gamma}2}} \| \omega_0 \|_{L^2}^{4},
  \end{equation*}
which is the same bound for the second term.

Combining the above estimates, we get
  \begin{equation*}
    \| \widetilde{\omega}_b-\widetilde{\omega}_a \|_{L^2}^2 \lesssim_{\beta,\gamma} \frac{h^{\gamma}}{\nu^{\frac{5-\beta +\gamma}2}} \| \omega_0 \|_{L^2}^{4}+M_{\omega_0}((\nu h)^{-1/2})^2 
  \end{equation*}
for all $0 \leq \gamma<\beta-1$. Thus, the Cauchy--Schwarz inequality leads to
  \begin{align*}
  \sup_{\substack{t,s\in[0,T] \\|t-s|\le h}}
  \big| \|\widetilde{\omega}_t\|_{L^2}^2 -\|\widetilde{\omega}_s\|_{L^2}^2 \big|
  & \le \sup_{\substack{t,s\in[0,T] \\|t-s|\le h}}
  (\|\widetilde{\omega}_t\|_{L^2}+ \|\widetilde{\omega}_s\|_{L^2}) \|\widetilde{\omega}_t-\widetilde{\omega}_s\|_{L^2}  \\
  & \lesssim_{\beta,\gamma} \frac{(\nu h)^{\gamma/2}}{\nu^{\frac{5-\beta + 3\gamma}4}} \| \omega_0 \|_{L^2}^{3}+M_{\omega_0}((\nu h)^{-1/2}) \| \omega_0 \|_{L^2}
  \end{align*}
for all $0 \leq \gamma<\beta-1$. We get estimate \eqref{eq:limit-time-modulus}.

  The decomposition~\eqref{eq:mild-split} and Lemma \ref{lem: determined nonlinear} applied with $\epsilon=\frac{\beta-1-\gamma}{2}$ also give \eqref{eq:limit-tail}:
  \begin{align*}
    \sup_{t\in[0,T]}  \|  (I-P_{\le R})\widetilde{\omega}_t\|_{L^2}
    & \le \sup_{t\in[0,T]} \| (I-P_{\le R}) \omega_0\|_{L^2}  +R^{-\gamma} \sup_{t\in[0,T]} \|N(t)\|_{H^\gamma} \\
    & \lesssim_{\beta,\gamma} M_{\omega_0}(R)+\frac{T^{\frac{\beta-1-\gamma}{4}}}{R^{\gamma} \nu^{\frac{5-\beta + \gamma}4}} \| \omega_0 \|_{L^2}^2,
  \end{align*}
for all $0 \leq \gamma<\beta-1$. The proof is finished.
\end{proof}

Using Lemma \ref{lem: convergence rate pre} and the estimate in \eqref{est: convergence_L2}, we give the following quantitative estimate for the scaling limit in $L^2 \big(\Omega; L^{\infty}([0,T]; L^2(\T^2)) \big)$.

\begin{theorem}\label{thm-quanti-C-L-2}
Fix $\nu > 0$, $\beta \in (1,2]$ and $\alpha\in \big(0, \frac{1}{4} (\beta+1-\sqrt{(\beta-1)^2+4}\, ) \big)$. Let $\omega^n$ and $\widetilde{\omega}$ be respectively the unique solution of \eqref{2D Euler/mSQG-n} and of \eqref{2D determined NSE/mSQG}, with common initial datum $\omega_0$. Then for any $T > 0$ and $L\ge 1$, it holds for all $q_0>\frac{2}{\beta-1}$, $s\in (0,\beta-1) \cap (0,1-\alpha)$ and $\gamma \in (0,\beta-1)$,
  \begin{equation*}
    \E \Big[ \sup_{t \in [0, T]} \|\omega_t^n - \widetilde{\omega}_t\|_{L^2}^2 \Big] \lesssim_{\nu,T,q_0,s,\gamma,\epsilon} \exp \left( C \|\omega_0\|_{L^2}^{q_0} \right) \|\omega_0\|_{L^2} \min_{L\geq 1} \Big\{ M_{\omega_0}(n^{1/L}), n^{-\frac{\gamma}{L} }, n^{-s+\epsilon+\frac{1-\eta(\alpha,1-s)}{L}} \Big\},
  \end{equation*}
  where $\eta(\alpha,1-s)$ is given by Lemma \ref{lem: uniform-anomalous}.
\end{theorem}

\begin{proof}
  Recall the estimate \eqref{est: convergence_L2}:
  \begin{align*}
     & \E \Big[ \sup_{t \in [0,T]} \| \omega^n_t - \widetilde{\omega}_t \|_{L^2}^2 \Big] \\
     & \lesssim 2 R^{s}  \| \omega_0 \|_{L^2}  \E\Big[ \sup_{t \in [0,T]} \| \omega^n_t-\widetilde{\omega}_t \|_{H^{-s}}  \Big] + 4  \|\omega_0 \|_{L^2} \sup_{t\in[0,T]}  \|  (I-P_{\le R})\widetilde{\omega}_t\|_{L^2} \nonumber \\
     & \quad + \frac{\|\omega_0\|_{L^2}}{h^{\frac{1-\eta(\alpha,1-s)}{2}}} \Big(1+ \|\omega_0 \|_{L^2}^{\frac{2(1-\alpha)}{\beta-2\alpha-1}} \Big)  \E\Big[ \sup_{t \in [0,T]} \| \omega^n_t-\widetilde{\omega}_t \|_{H^{-s}}^2  \Big]^{1/2}  + \sup_{ \substack{t,s  \in [0,T] \\ |t-s| \le h } } |E(t) - E(s)|.
  \end{align*}
Applying Proposition \ref{prop: scal_lim_negative} with $q_0>\frac{2}{\beta-1}$ and Lemma \ref{lem: convergence rate pre}, we get
  \begin{align*}
     & \E \Big[ \sup_{t \in [0,T]} \| \omega^n_t - \widetilde{\omega}_t \|_{L^2}^2 \Big] \\
     & \lesssim \exp \left( C \|\omega_0\|_{L^2}^{q_0} \right) \frac{R^s \|\omega_0\|_{L^2}^2}{n^{s-\epsilon}} +  M_{\omega_0}(R) \| \omega_0 \|_{L^2} +\frac{T^{\frac{\beta-1-\gamma}{4}}}{R^{\gamma} \nu^{\frac{5-\beta + 3\gamma}4}} \| \omega_0 \|_{L^2}^3 \\
     & \quad +  \exp \left( C \|\omega_0\|_{L^2}^{q_0} \right) \Big(1+ \|\omega_0 \|_{L^2}^{\frac{2(1-\alpha)}{\beta-2\alpha-1}} \Big)  \frac{\|\omega_0\|_{L^2}^2}{h^{\frac{1-\eta(\alpha,1-s)}{2}} n^{s-\epsilon}} + \frac{(\nu h)^{\gamma/2}}{\nu^{\frac{5-\beta + \gamma}4}} \| \omega_0 \|_{L^2}^{3}+M_{\omega_0}((\nu h)^{-1/2}) \| \omega_0 \|_{L^2},
  \end{align*}
  for any $s\in (0,\beta-1) \cap (0,1-\alpha)$ and $\gamma \in (0,\beta-1)$. Taking $R=(\nu h)^{-1/2}=n^{1/L}$, we get
  \begin{align*}
     & \E \Big[ \sup_{t \in [0,T]} \| \omega^n_t - \widetilde{\omega}_t \|_{L^2}^2\Big] \\
     & \lesssim \exp \left( C \|\omega_0\|_{L^2}^{q_0} \right) \|\omega_0\|_{L^2} \bigg( \frac{1 }{n^{ \frac{L-1}{L} s -\epsilon}} +  M_{\omega_0}(n^{1/L})+ \frac{1}{n^{\frac{1}{L} \gamma} } + \frac{1}{n^{s-\epsilon-\frac{1-\eta(\alpha,1-s)}{L}}} \bigg),
  \end{align*}
which holds for all $L \geq 1$. Due to $\eta(\alpha,1-s)\leq \frac{1-s-\alpha}{4-(1-s)-\alpha} <1- s$, then $\frac{L-1}{L} s > s - \frac{1-\eta(\alpha,1-s)}{L}$ and we obtain the desired result.
\end{proof}

The convergence rate in $L^2 \big(\Omega;C([0,T];L^2(\T^2)) \big)$ depends on the Fourier mode(s) of the initial value $\omega_0$, rather than only on the $L^2(\T^2)$-norm of $\omega_0$, as demonstrated below.

\begin{proposition}\label{prop-nonuniform-C-L-2}
There exits a sequence of initial data $\{\omega_0^n \}_{n \geq 1}$ with the same $L^2(\T^2)$-norm such that the corresponding solution $\omega_t^n$ of Eq. \eqref{2D Euler/mSQG-n} with initial data $\omega_0^n$ satisfies
  \begin{equation} \label{eq: non-uniform convergence}
    \liminf_{n \rightarrow +\infty} \E \Big[ \sup_{t \in [0, T]} \|\omega_t^n - \widetilde{\omega}^n_t\|_{L^2}^2 \Big] >0,
  \end{equation}
where $\widetilde{\omega}^n_t$ is the solution of Eq. \eqref{2D determined NSE/mSQG} with initial data $\omega_0^n$.    
\end{proposition}

\begin{proof}
For convenience, we set $\| \omega^n_0\|_{L^2}=1$ for all $n \geq 1$. Write \eqref{2D Euler/mSQG-n} in mild form as 
  \begin{equation}\label{eq:mild-split-n} 
  \omega^n_t =e^{\nu t\Delta}\omega_0^n +N_{n}(t) +M_n(t),
  \end{equation}
where
  \[
    N_{n}(t)=-\int_0^t e^{\nu(t-r)\Delta} (u^n_r\cdot\nabla \omega^n_r) \d r,
    \quad M_n(t)= -\sqrt{2\nu}\sum_k\theta_{n,k} \int_0^t e^{\nu(t-r)\Delta} (\sigma_k\cdot\nabla\omega^n_r) \d B_r^k
  \]
and $u^n_r = \nabla^\perp(-\Delta)^{-\beta/2} \omega^n_r= \frac{\rm i}{(2\pi)^{\beta -1}} \sum_{m\ne 0} \frac{m^\perp}{|m|^{\beta}} \widehat{\omega^n}_m(r) e_m$. Similarly to \eqref{eq:mild-split}, we have
  $$ \widetilde{\omega}^n_t = e^{\nu t\Delta} \omega^n_0 + \widetilde N_n(t), \quad \widetilde N_n(t)= -\int_0^t e^{\nu(t-r)\Delta}  \widetilde{u}^n_r\cdot\nabla\widetilde{\omega}^n_r\d r. $$
Therefore, the difference $\omega^n_t-\widetilde{\omega}^n_t$ satisfies
  \begin{equation*}
    \omega^n_t-\widetilde{\omega}^n_t = N_n(t)- \widetilde N_n (t) +M_n(t). 
  \end{equation*}
  Then, by triangle inequality, we obtain that for any $t_n \in [0,T]$,
  \begin{equation*}
    \E \Big[ \sup_{t \in [0, T]} \|\omega_t^n - \widetilde{\omega}^n_t\|_{L^2}^2 \Big]^{\frac{1}{2}}  \geq \E \Big[ \|M_n(t_n)\|_{L^2}^2 \Big]^{\frac{1}{2}} -\E \Big[  \|N_n(t_n)\|_{L^2}^2 \Big]^{\frac{1}{2}} -\E \Big[  \|\widetilde N_n(t_n)\|_{L^2}^2 \Big]^{\frac{1}{2}}.
  \end{equation*}

Similarly to the estimate \eqref{lem: determined nonlinear.2} with $\gamma=0$, we have for any $0<\epsilon<\beta-1-\alpha$,
  \begin{align*}
     & \E \| u^n_r\cdot\nabla \omega^n_r \|_{H^{-1}}^2 \lesssim_{\beta} \|\omega^n_0 \|_{L^2}^2\, \E \|\omega^n_r \|_{H^{2-\beta +\epsilon}}^2  .
  \end{align*}
Using the uniform anomalous regularity \eqref{eq: uniform-Sobolev} with $\tilde{\alpha}=\beta-1-\epsilon>\alpha$ (i.e. $2-\beta+\epsilon=1-\tilde{\alpha}<1-\alpha$) and arguing as in the proof of Lemma \ref{lem: determined nonlinear} with $a=0$ and $b=t_n$, we get
  \begin{equation*}
    \E \| N_n(t_n)\|_{L^2}^2  \le C |t_n|^{\eta(\alpha,\beta-1-\epsilon)} \Big(1+ \|\omega_0^n \|_{L^2}^{\frac{4(1-\alpha)}{\beta-2\alpha-1}} \Big)  \|\omega_0^n\|_{L^2}^{4} \lesssim |t_n|^{\eta(\alpha,\beta-1-\epsilon)},
  \end{equation*}
  where function $\eta(\alpha,\beta-1-\epsilon)$ is given by Lemma \ref{lem: uniform-anomalous}. Thus, if $t_n \rightarrow 0$ as $n \rightarrow +\infty$, we have $\lim_{n \rightarrow +\infty} \E \| N_n(t_n)\|_{L^2}^2=0$; similarly we also have $\lim_{n \rightarrow +\infty} \E \| \widetilde N_n(t_n)\|_{L^2}^2=0$. Therefore, to prove \eqref{eq: non-uniform convergence}, it suffices to prove
  \begin{equation}\label{limit-marting}
    \liminf_{n \rightarrow +\infty} \E \Big[ \|M_n(t_n)\|_{L^2}^2 \Big]>0.
  \end{equation}

By It\^o's isometry and Jessen's inequality, we have
  \begin{align*}
    \E \Big[ \|M_n(t_n)\|_{L^2}^2 \Big] & =2\nu \int_0^{t_n} \sum_k\theta_{n,k}^2 \E \Big[ \| e^{\nu(t_n-r)\Delta} \big(\sigma_k\cdot\nabla\omega^n_r \big) \|_{L^2}^2 \Big] \, \d r \\
    & \geq 2\nu \int_0^{t_n} \sum_k\theta_{n,k}^2 \Big\| \E \big[  e^{\nu(t_n-r)\Delta} \big(\sigma_k\cdot\nabla\omega^n_r \big) \big]  \Big\|_{L^2}^2 \, \d r.
  \end{align*}
Using \eqref{eq:mild-split-n} and the triangle inequality, we obtain
  \begin{align*}
    \big\| \E \big[  e^{\nu(t_n-r)\Delta} \big(\sigma_k \cdot\nabla \omega^n_r \big) \big] \big\|_{L^2}
    & = \big\|e^{ \nu(t_n-r)\Delta} \big(\sigma_k\cdot\nabla e^{ \nu r\Delta} \omega_0^n  \big)+ e^{ \nu(t_n-r)\Delta} \big(\sigma_k\cdot\nabla\E[N_n(r)] \big)\big\|_{L^2} \\
    & \geq \big\|e^{ \nu(t_n-r)\Delta} \big(\sigma_k\cdot\nabla e^{ \nu r\Delta} \omega_0^n  \big)\big\|_{L^2} -\big\| e^{ \nu(t_n-r)\Delta} \big(\sigma_k\cdot\nabla\E[N_n(r)]  \big)\big\|_{L^2}.
  \end{align*}
A simple manipulation leads to
  \begin{equation}\label{eq-marting-estim-1}
  \aligned
  \big\| \E \big[  e^{\nu(t_n-r)\Delta} \big(\sigma_k \cdot\nabla \omega^n_r \big) \big] \big\|_{L^2}^2  
  & \ge \frac12 \big\|e^{ \nu(t_n-r)\Delta} \big(\sigma_k\cdot\nabla e^{ \nu r\Delta} \omega_0^n  \big)\big\|_{L^2}^2 \\ 
  &\quad -\big\| e^{ \nu(t_n-r)\Delta} \big(\sigma_k\cdot\nabla\E[N_n(r)] \big)\big\|_{L^2}^2. 
  \endaligned 
  \end{equation}
By Jessen's inequality and the classical semigroup estimate (cf. Lemma \ref{lem-heat-semigroup-1}), we get, for small $\epsilon>0$,
  \begin{align*}
     & \int_0^{t_n} \sum_k\theta_{n,k}^2 \big\| e^{ \nu(t_n-r)\Delta} \big(\sigma_k\cdot\nabla\E[N_n(r)]  \big)\big\|_{L^2}^2 \, \d r \\
     & \lesssim   \E \int_0^{t_n} \sum_k\theta_{n,k}^2 |t_n-r|^{\epsilon-1} \big\| \sigma_k\cdot\nabla N_n(r)  \big\|_{H^{-1+\epsilon}}^2 \, \d r \\
     & \lesssim \int_0^{t_n}  |t_n-r|^{\epsilon-1}  \, \d r \sum_k\theta_{n,k}^2 \,  \E \sup_{r \in [0,t_n]} \big\| \sigma_{k} N_n(r)\big\|_{H^{\epsilon}}^2 \\
     & \lesssim  t_n^\epsilon \Big(\sum_k |k|^{2\epsilon} \theta_{n,k}^2 \Big) \,  \E \sup_{r \in [0,t_n]} \big\| N_n(r)\big\|_{H^{\epsilon}}^2 .
  \end{align*}
Analogously to the estimate for $ \E \| N_n(t_n)\|_{L^2}^2$ with taking $\gamma=\epsilon$, we have
  \begin{equation*}
    \E \sup_{r \in [0,t_n]} \big\| N_n(r)\big\|_{H^{\epsilon}}^2 \lesssim  |t_n|^{\frac{\beta-1-2\epsilon-\alpha}{5-\beta-\alpha+2\epsilon}}.
  \end{equation*}
By Lemma \ref{lem: power tail}, we have $\sum_k |k|^{2\epsilon} \theta_{n,k}^2 \lesssim n^{2\epsilon}$; choosing $t_n \leq e^{-n}$, we can easily deduce
  \begin{equation*}
    \lim_{n \rightarrow +\infty}  \int_0^{t_n} \sum_k\theta_{n,k}^2 \big\| e^{ \nu(t_n-r)\Delta} \big(\sigma_k\cdot\nabla\E[N_n(r)]  \big)\big\|_{L^2}^2 \, \d r =0.
  \end{equation*}
Combining this limit with \eqref{eq-marting-estim-1}, we conclude that to prove \eqref{limit-marting}, it suffices to show
  \begin{equation*}
    \liminf_{n \rightarrow +\infty} \int_0^{t_n} \sum_k\theta_{n,k}^2 \big\|e^{ \nu(t_n-r)\Delta} \big(\sigma_k\cdot\nabla e^{ \nu r\Delta} \omega_0^n  \big)\big\|_{L^2}^2 \, \d r>0.
  \end{equation*}

Let $E_1=(1,0)$ and $t_n \leq e^{-n}$ with $\frac{1}{\sqrt{t_n}}$ an integer. Taking $\omega^n_0(x)= \sqrt{2} \cos \{ 2 \pi \frac{E_1}{\sqrt{t_n}} \cdot x\}$. We first show that there exits $0<c<1$ and $L> 1$, such that 
  \begin{equation}\label{eq-lower-bound}
  \sum_{|k| \leq L n} \theta_{n,k}^2 (a_k\cdot E_1)^2 >c .
  \end{equation}
Indeed, by Lemma \ref{lem: power tail}, it is clear that there exits $0<c<1$ and $L>1$, such that $\sum_{|k| \leq L n} \theta_{n,k}^2 >c$. Since $(a_k\cdot E_1)^2 + \big(\frac{k}{|k|} \cdot E_1\big)^2 =1$ for any $k\in \Z_0^2$, by the rotation symmetry of $\theta_{n,k}$ in $k$, one see that \eqref{eq-lower-bound} holds up to replacing $c$ by $c/2$. 
Then, by elementary semigroup property,
  \begin{align*}
    &\int_0^{t_n} \sum_k\theta_{n,k}^2 \big\|e^{ \nu(t_n-r)\Delta} \big(\sigma_k\cdot\nabla e^{ \nu r\Delta} \omega_0^n \big)\big\|_{L^2}^2 \d r \\
    & \gtrsim \sum_k\theta_{n,k}^2 \int_0^{t_n} e^{- 8\pi^2 \nu \big|\frac{E_1}{\sqrt{t_n}}+k \big|^2 (t_n-r) } \frac{(a_k\cdot E_1)^2}{t_n} e^{-8\pi^2 \nu \frac{r}{t_n} } \d r \\
    & \gtrsim \frac{1}{t_n} \sum_{|k| \leq L n} \theta_{n,k}^2 (a_k\cdot E_1)^2 \int_0^{t_n} e^{-8\pi^2 \nu \big|\frac{E_1}{\sqrt{t_n}}+k \big|^2 (t_n-r) } e^{-8\pi^2 \nu \frac{r}{t_n} } \d r.
  \end{align*}
When $n$ is big enough, $ \big|\frac{E_1}{\sqrt{t_n}}+k \big|^2 = \frac{|E_1+ \sqrt{t_n} k|^2}{t_n} \leq \frac{2}{ t_n}$ holds for all $k \leq L n$. Using $\frac{t_n-r}{t_n} \leq 1$ and $\frac{r}{t_n} \leq 1$, we deduce from \eqref{eq-lower-bound} that
  \begin{equation*}
    \int_0^{t_n} \sum_k\theta_{n,k}^2 \big\|e^{ \nu(t_n-r)\Delta} \big(\sigma_k\cdot\nabla e^{ \nu r\Delta} \omega_0^n \big)\big\|_{L^2}^2 \d r 
    \gtrsim e^{-24 \pi^2 \nu} \frac{c}{t_n} \int_{0}^{t_n} \d r 
    \geq e^{-24 \pi^2 \nu} c.
  \end{equation*}
The proof is completed.
\end{proof}

\appendix

\section{Analysis tools}\label{append-tools}

We first collect some useful properties of heat semigroup $\{e^{t\Delta} \}_{t\ge 0}$, which will be stated for higher dimensional cases $\T^d$. We omit their proofs since they are classical.

\begin{lemma}\label{lem-heat-semigroup-1}
Let $f\in H^a(\T^d)$, $a\in \R$. Then,
\begin{itemize}
\item[\rm (i)] for any $\rho\ge 0$, it holds $\|e^{t\Delta} f\|_{H^{a+\rho}} \le C_\rho\, t^{-\rho/2} \|f \|_{H^{a}}$ for some constant increasing in $\rho$;
\item[\rm (ii)] for any $\rho\in [0,2]$, it holds $\|(I-e^{t\Delta}) f\|_{H^{a-\rho}} \lesssim t^{\rho/2} \|f \|_{H^{a}}$ .
\end{itemize}
\end{lemma}

\begin{lemma}\label{lem-heat-semigroup-2}
For any $a\in \R$ and $f\in L^2([s,t]; H^a(\T^d))$, it holds
  $$ \bigg\| \int_s^t e^{\nu (t-r)\Delta} f_r \d r \bigg\|_{H^{a+1}}^2 \lesssim \frac1\nu \int_s^t \|f_r \|_{H^a}^2 \d r. $$
\end{lemma}

We also recall the following result on the product of Sobolev functions, see e.g. \cite[Corollary 2.55, p.90]{BCD11}.

\begin{lemma}\label{lem-product-Sobolev}
Let $a, b<d/2$ with $a + b>0$. For any $f\in H^a(\T^d)$ and $g\in H^b(\T^d)$, the product $fg$ belongs to $H^{a+b-d/2}(\T^d)$ and verifies
  $$\|fg \|_{H^{a+b-d/2}} \lesssim_{d,a,b} \|f \|_{H^a} \|g \|_{H^b} .$$
\end{lemma}

The following inequality is proved in our previous work \cite[Lemma 2.6]{Luo_2024_Elementary}.
\begin{lemma} \label{lem: bound p}
  For $a,b \geq 0$ and $p>1$
  \begin{equation*}
    |a^p - b^p| \leq \frac{p}{\sqrt{p-1}} \sqrt{a^p + b^p} \sqrt{(a^{p-1}- b^{p-1})(a-b)}.
  \end{equation*}
\end{lemma}

\begin{proof}
  Recall that for convex functions $f : [0,\infty) \to \R$, it holds
  \[f(y) \geq f(x) + f'(x) (y-x), \quad x,y \in \R_{+}. \]
 Since $x^p$ and $x^{\frac{p}{p-1}}$ are convex for $p>1$, we obtain 
  \begin{align*}
    p b^{p-1} (a-b) & \leq a^p - b^p \leq p a^{p-1} (a-b) ,
    \\\frac{p}{p-1} b (a^{p-1} - b^{p-1})&\leq a^p- b^p \leq \frac{p}{p-1} a(a^{p-1} - b^{p-1}).
  \end{align*}
  Taking product of the two inequalities and then taking square root, we get the desired result.
\end{proof}

The following result gives us uniform upper and lower bounds on the tail of coefficients in the Kraichnan-type noise.

\begin{lemma} \label{lem: power tail}
  For any $d \geq 2$, \(0<\alpha\le\frac12\) and \(r\ge1\), define
  \[
    T_\alpha(r)  :=\sum_{|l|\ge r} |l|^{-d-2\alpha}.
  \]
  Then there exist constants $C_1>0$ and $C_2>0$, independent of \(r\) and \(\alpha\), such that
  \begin{equation*}
    \frac{C_1}{\alpha} r^{-2\alpha} \leq T_\alpha(r)\leq
    \frac{C_2}{\alpha} r^{-2\alpha}.
  \end{equation*}
\end{lemma}

\begin{proof}
  For any \(j\ge0\), set
  \[
    A_j :=\bigl\{l\in\mathbb Z^d:2^jr\le |l|<2^{j+1}r \bigr\}.
  \]
  Then there exist constants $C_1 \le C_2$ depending only on the dimension $d \geq 2$, such that
  $$ C_1 (2^j r)^d \le \# A_{j} \le C_2 (2^j r)^d, \quad \forall\, j\ge 0. $$
  Therefore, 
  \begin{align*}
    \sum_{l\in A_j}|l|^{-d-2\alpha} & \leq (2^{j} r)^{-d-2\alpha} \# A_{j} \leq  C_2 (2^{j} r)^{-2\alpha}, \\
    \sum_{l\in A_j}|l|^{-d-2\alpha} & \geq (2^{j+1} r)^{-d-2\alpha} \# A_{j} \geq 2^{-d-2\alpha} C_1 (2^j r)^{-2\alpha}.
  \end{align*}
  Summing over \(j\ge0\), we obtain
  \begin{align*}
    T_{\alpha}(r) & =  \sum_{j=0}^{+\infty} \sum_{l\in A_j}|l|^{-d-2\alpha}  \leq  C_2 \sum_{j=0}^{+\infty}  (2^{j} r)^{-2\alpha} = C_2 \frac{r^{-2\alpha}}{1-2^{-2\alpha}}, \\
    T_{\alpha}(r) & =  \sum_{j=0}^{+\infty} \sum_{l\in A_j}|l|^{-d-2\alpha}  \geq  2^{-d-2\alpha} C_1 \sum_{j=0}^{+\infty}  (2^j r)^{-2\alpha} =  2^{-d} C_1 \frac{r^{-2\alpha}}{2^{2\alpha} -1}.
  \end{align*}
  Since \(0<2\alpha\le1\), concavity of \( f(t):=1-2^{-t}\) gives
  \begin{align*}
     & 1-2^{-2\alpha}\ge 2 \alpha f(1) + (1-2\alpha) f(0) = \alpha \, ,\\
     & 1-2^{2\alpha} \geq  2 \alpha f(-1) + (1-2\alpha) f(0)=-2\alpha.
  \end{align*}
  Consequently,
  \begin{align*}
    T_{\alpha}(r) & \leq  C_2 \frac{r^{-2\alpha}}{1-2^{-2\alpha}} \leq \frac{C_2}{\alpha} r^{-2\alpha} ;               \\
    T_{\alpha}(r) & \geq 2^{-d} C_1 \frac{r^{-2\alpha}}{2^{2\alpha} -1}\ge 2^{-d} C_1 \frac{r^{-2\alpha}}{2\alpha}
    \ge 2^{-d-1} \frac{C_1}{\alpha}r^{-2\alpha}.
  \end{align*}
  Renaming $2^{-d-1} C_1$ as $C_1$ gives us the desired result.
\end{proof}

The following lemma provides an estimate with kernel $\min\{1,\frac{|l|}{\gamma|k|}\}$ in line $\Gamma_{q_0,l}$.

\begin{lemma}[Unified affine-line congestion]
  \label{lem:line-congestion}
  Let $l\in\mathbb Z_0^d$ and $q_0\in\mathbb Z^d$, and define the line
  \[
    \Gamma_{q_0,l}:=q_0+\langle l\rangle
    =\{q_0+ml;m\in\mathbb Z\}.
  \]
  For $\gamma\ge1$ and $k \in \Z_0^d$, set
  \[
    \Phi_\gamma(k,l)
    :=\min\left\{1,\frac{|l|}{\gamma|k|}\right\}.
  \]
  Then, for every $0<A\le B$, the following estimate holds
  \begin{equation} \label{eq:line-congestion}
    \sum_{\substack{k \in\Gamma_{q_0,l}\\     A \le \max\{|k|,|l|\} \le B}}  \Phi_\gamma(k,l)  \leq   \sum_{\substack{k \in\Gamma_{q_0,l}\\
        A \le \max\{|k|,|l|\}  \\  |k|\le B}}  \Phi_\gamma(k,l)  \leq  9+\frac2\gamma\log \frac BA.
  \end{equation}
\end{lemma}

\begin{proof}
  Distinct points of $\Gamma_{q_0,l}$ are separated by at least $|l|$. More
  precisely, for every $r>0$,
  \begin{equation*}
    \#\{k \in\Gamma_{q_0,l}:|k| \leq r\}
    \le 1+\frac{2r}{|l|}.
  \end{equation*}
  Thus, using the fact $\Phi_\gamma(k,l)\le1$, we have
  \begin{equation}  \label{eq:inner-line-part}
    \sum_{\substack{k \in\Gamma_{q_0,l}, \; |k| \leq |l|\\   A \le \max\{|k|,|l|\}\le B}} \Phi_\gamma(k,l)  
    \leq \#\{k \in\Gamma_{q_0,l}:|k| \leq |l|\} \leq 3.
  \end{equation}

  It remains to deal with the region  $|k|>|l|$. In this case, since $\gamma\ge1$, we have
  \[
    \Phi_\gamma(k,l) =\min\left\{1,\frac{|l|}{\gamma|k|}\right\} =\frac{|l|}{\gamma|k|},
  \]
  and $\max\{|k|,|l|\}=|k|$. Thus, we have
  \begin{align*}
    \sum_{\substack{k \in\Gamma_{q_0,l}, \; |k| > |l| \\   A \le \max\{|k|,|l|\}\le B}}  \Phi_\gamma(k,l) 
    = \sum_{\substack{k \in\Gamma_{q_0,l}, \; |k| > |l| \\   A \le |k|\le B}} \frac{|l|}{\gamma|k|}.
  \end{align*}

Let $\Pi_{l^{\perp}} q_0:=q_0-\frac{q_0 \cdot l}{|l|^2} l$ denote the orthogonal projection of $q_0$ on the hyperplane $l^{\perp}$. Define the following two rays:
\begin{equation*}
    \Gamma_{q_0,l}^{+}:=\{k=\Pi_{l^{\perp}} q_0+r l; \, r\in \R_+  \} \cap \Gamma_{q_0,l}; \quad \Gamma_{q_0,l}^{-}:=\{k=\Pi_{l^{\perp}} q_0-r l; \, r \in \R_+\} \cap \Gamma_{q_0,l}.
\end{equation*}
Let $k_{+}$ (resp. $k_{-}$) be the point on the ray $\Gamma_{q_0,l}^{+}$ (resp. $\Gamma_{q_0,l}^{-}$) with the smallest norm but strictly greater than $|l|$. Consequently, $k_{+} \cdot l>0$ and $k_{-} \cdot l <0$; it holds that
  \begin{equation*}
    \{ k \in\Gamma_{q_0,l} ; |k| >|l| \} =\{k=k_{+} + n l ; \, n \in \N  \} \cup \{k=k_{-} - n l ; \, n \in \N  \}.
  \end{equation*}
  Thus, for any $k=k_{+}+n l$, we have
  \begin{equation*}
    \frac{|k|}{|l|} \geq  \frac{|k \cdot l|}{|l|^2} =  \frac{|n l \cdot l|}{|l|^2}+\frac{|k_{+} \cdot l|}{|l|^2} \geq n, \quad \frac{|k|}{|l|} \geq \frac{|k_{+}|}{|l|} \geq 1;
  \end{equation*}
  for any $k=k_{-}-n l$, similar results also hold. Then, we obtain
  \begin{align*}
    \sum_{\substack{k \in\Gamma_{q_0,l}, \; |k| > |l| \\   A \le |k|\le B}} \frac{|l|}{\gamma|k|} 
    & = \sum_{ \{n \in \N; A \leq |k_{+}+n l| \leq B \} } \frac{|l|}{\gamma|k_{+}+n l|} + \sum_{ \{n \in \N; A \leq |k_{-}-n l| \leq B \} } \frac{|l|}{\gamma|k_{-}-n l|} \\
    & \leq  \sum_{ \{n \in \N; A \leq |k_{+}+n l| \leq B \} } \frac{1}{\gamma \, \max\{1,  n, \frac{A}{|l|} \}} + \sum_{ \{n \in \N ; A \leq |k_{-}-n l| \leq B \} } \frac{1}{\gamma \, \max\{1, n, \frac{A}{|l|} \}}.
  \end{align*}
Due to $\{n \in \N;  A \leq |k_{+}+n l| \leq B \}  \subset \{n \in \N;  n \leq B/|l|\}$, we have
  \begin{align*}
    \sum_{ \{n \in \N; A \leq |k_{+}+n l| \leq B \} } \frac{1}{\gamma \, \max\{1, n, \frac{A}{|l|} \}} 
    & \leq \frac{2}{\gamma} + \sum_{\{n \in \N; 2 \leq n \leq B/|l|\}}  \frac{1}{\gamma \, \max\{n, \frac{A}{|l|} \}} \\
    & \leq \frac{2}{\gamma} + \int_{1}^{\frac{B}{|l|}}  \frac{1}{\gamma \, \max\{r, \frac{A}{|l|} \}} \, \d r \\
    & \leq \frac{2}{\gamma} + \Big| \int_{1}^{\frac{A}{|l|}}  \frac{1}{\gamma \, \frac{A}{|l|}} \, \d r \Big|+   \int_{\frac{A}{|l|}}^{\frac{B}{|l|}}  \frac{1}{\gamma \,r} \, \d r = \frac{3}{\gamma}+ \frac{1}{\gamma} \log \frac{B}{A}.
  \end{align*}
The estimate of the other part is similar. Thus, we have
  \begin{equation*}
    \sum_{\substack{k \in\Gamma_{q_0,l}, \; |k| > |l| \\ A \le \max\{|k|,|l|\}\le B}}  \Phi_\gamma(k,l) 
    = \sum_{\substack{k \in\Gamma_{q_0,l}, \; |k|>|l| \\ A \le \max\{|k|,|l|\}  \\  |k|\le B}} \Phi_\gamma(k,l) \leq \frac{6}{\gamma}+ \frac{2}{\gamma} \log \frac{B}{A}.
  \end{equation*}
Since $\gamma \geq 1$, combining the above estimate and \eqref{eq:inner-line-part}, we get \eqref{eq:line-congestion}.
\end{proof}

\section{Weighted discrete Poincar\'e inequality} \label{sec: poincare inequality}

In this section, we present an improved version of the weighted discrete Poincaré inequality, building upon those in \cite{Luo_2024_Elementary} and \cite{Rowan_2025_Obukhovcorrsin}. Our improved inequality provides a necessary tool for proving the anomalous regularity and the scaling limits in positive Sobolev spaces.

Firstly, we define a new weighted function. For any fixed $R\ge1$ and $k,l\in\mathbb{Z}_0^d$, we write
\[
  \mu_{R}(k,l):=\min\left\{\frac{|l|}{R},|k|\right\}.
\]
Thus $\mu_{R}(k,l)$ is the scale-truncated size of the increment $l$ relative to the mode $k$.  Define
\[
  W_{R}(k):= \sum_{l\in\mathbb{Z}_0^d}  \mu_{R}(k,l)^2 \theta_l^2  = \sum_{l\in\mathbb{Z}_0^d}  \min\left\{\frac{|l|^2}{R^2},|k|^2\right\} \theta_l^2, \quad \forall \, k \in \Z_0^d.
\]
The weight $W_{R}(k)$ measures the effective coercivity seen by the Fourier mode $k$:
\begin{equation*}
  W_{R}(k)= \frac{1}{R^2} \sum_{|l| \leq R |k|} |l|^2 \theta_{l}^2 + \bigg( \sum_{|l|> R |k|} \theta_{l}^2 \bigg) |k|^2.
\end{equation*}
Compared with the weight function in \cite{Rowan_2025_Obukhovcorrsin}, which consists only of the first summand above, $W_R(k)$ also captures information from high frequencies. This has two advantages: (i) $W_R(k)$ never vanishes and is in fact positive for every $k \in \Z_0^d$; (ii) for a suitable choice of $R$, the high-frequency part ensures the uniform lower bound on the weight functions required by the scaling limit in Section \ref{sec: scaling limit}; see the following Example for details.

\begin{example} \label{exm: power law}
Fix $\alpha\in(0,1)$. For $n\in\N$, let
\[
\theta_{n,l}=Z_{n,\alpha}|l|^{-\frac d2-\alpha}\ind_{\{|l|\ge n\}},\quad l\in\Z_0^d,
\]
where $Z_{n,\alpha}$ is a normalizing constant. Then Lemma \ref{lem: power tail} shows that, for every $k\in\Z_0^d$, the weight function $W_{n,\alpha}(k):=\sum_l\theta_{n,l}^2  \min \big\{ \frac{|l|^2}{n^2},|k|^2 \big\}$ satisfies
\begin{align*}
  W_{n,\alpha}(k)
  &= \frac{1}{n^2}\sum_{|l|\le n|k|} |l|^2 \theta_{n,l}^2
  + \bigg(\sum_{|l|\ge n|k|}\theta_{n,l}^2 \bigg) |k|^2 \\
  &\ge |k|^2\sum_{|l|\ge n|k|}\theta_{n,l}^2
  = |k|^2\frac{T_\alpha(n|k|)}{T_\alpha(n)}
  \ge \frac{C_1}{C_2}|k|^{2-2\alpha} , \quad k\in\Z_0^d,
\end{align*}
where $C_1,C_2$ are the constants from Lemma \ref{lem: power tail}, uniform in $n\in\N$. In fact, $W_{n,\alpha}(k)\sim |k|^{2-2\alpha}$.

\end{example}

In order to establish the weighted discrete Poincar\'e  inequality, we make the following assumption on the coefficients $\{ \theta_{k} \}_{k} \in \ell^2(\Z_0^d)$. 

\begin{assumption} \label{ass: non-degenerate}
  The noise coefficient $\{ \theta_{k} \}_{k} \in \ell^2(\Z_0^d)$ satisfies $\theta_{k}=\theta_{-k}$ for all $k \in \Z_0^d$; moreover, it satisfies the non-degeneracy condition: there is $\delta>0$ such that for any $R\ge 1$,
  \begin{equation*}
    \sum_{l\in\mathbb{Z}_0^d} \frac{|\Pi_{l^\perp}k|}{|k|} \mu_{R}(k,l)^2  \theta_l^2 \ge \delta \sum_{l\in\mathbb{Z}_0^d}  \mu_{R}(k,l)^2 \theta_l^2 =  \delta W_{R}(k), \quad \forall\, k \in \Z_0^d.
  \end{equation*}
\end{assumption}

Assumption \ref{ass: non-degenerate} is a consequence of the following Assumption \ref{ass: non-degenerate-L}, which is similar to the conditions given by Rowan, cf. \cite[Assumption 1.15]{Rowan_2025_Obukhovcorrsin}.

\begin{assumption} \label{ass: non-degenerate-L}
  The noise coefficient $\{ \theta_{k} \}_{k} \in \ell^2(\Z_0^d)$ satisfies $\theta_{k}=\theta_{-k}$ for all $k \in \Z_0^d$. In addition, it satisfies the following non-degeneracy condition: there exists  $\delta>0$ such that
  \begin{equation*}
    \inf_{x \in \R^d, |x|=1 } \sum_{|k|\leq  r} |\Pi_{k^{\perp}} x | \,  |k|^2 \theta_{k}^{2} \geq \delta \sum_{|k|\leq r} \, |k|^2\theta_k^2, \quad \forall \, r \geq 1.
  \end{equation*}
\end{assumption}

\begin{lemma} \label{lem: sufficient ass}
  If $\{ \theta_{k} \}_{k} \in \ell^2(\Z_0^d)$ satisfies Assumption \ref{ass: non-degenerate-L}, then for any $R \geq 1$, it satisfies Assumption \ref{ass: non-degenerate} with the same constant $\delta>0$.
\end{lemma}

\begin{proof}
  For any fixed $k,l\in \Z_0^d$, it holds
  \begin{equation} \label{eq-mu-R}
  \mu_{R}(k,l)^2 =\min\left\{\frac{|l|^2}{R^2},|k|^2\right\} = 2|k|^2 |l|^2 \int_{R |k|}^\infty  \frac{\ind_{\{ r \geq |l|\}}}{r^3}\, \d r.    
  \end{equation}
  Thus, Fubini's theorem and Assumption \ref{ass: non-degenerate-L} yield, for every unit vector $x\in \R^d$,
  \begin{align*}
    \sum_l |\Pi_{l^\perp}x| \mu_{R}(k,l)^2 \theta_l^2
     & = 2|k|^2 \sum_l|\Pi_{l^\perp}x|\, |l|^2 \theta_l^2 \int_{R |k|}^\infty  \frac{\ind_{\{ r \geq |l|\}}}{r^3}\, \d r \\
     & = 2|k|^2 \int_{R |k|}^\infty \frac1{r^3} \sum_{|l|\leq r} |\Pi_{l^\perp}x|\, |l|^2 \theta_l^2 \, \d r \\
     & \ge 2 \delta |k|^2 \int_{R |k|}^\infty \frac1{r^3}
    \sum_{|l| \leq r} |l|^2 \theta_l^2  \, \d r \\ 
    &= 2\delta |k|^2 \sum_l |l|^2 \theta_l^2 \int_{R |k|}^\infty \frac{\ind_{\{ r \geq |l|\}}}{r^3} \,\d r.
  \end{align*}
  Again by \eqref{eq-mu-R} and the definition of $W_R(k)$, we arrive at
  $$\sum_l |\Pi_{l^\perp}x| \mu_{R}(k,l)^2 \theta_l^2 \ge \delta W_R(k) . $$
  Finally, taking $x=k/|k|$, we obtain Assumption \ref{ass: non-degenerate}.
\end{proof}

When $\{\theta_{k}\}_{k \in \Z_0^d}$ are symmetric in $k \in \Z_0^d$, then both Assumptions \ref{ass: non-degenerate} and \ref{ass: non-degenerate-L} hold. 

\begin{lemma} \label{lem: symmetry}
  Assume that $\{\theta_{k}\}_{k\in \Z_0^d} $ are radially symmetric, i.e.
  \begin{equation*}
    \theta_{k} = \theta_{l}, \quad \forall\, k, l \in \Z_0^d, \quad \mbox{with} \quad |k| = |l|.
  \end{equation*}
  Then Assumptions \ref{ass: non-degenerate} and \ref{ass: non-degenerate-L} are satisfied with $\delta=(d-1)/d$.
\end{lemma}

\begin{proof}
Due to $|\Pi_{k^{\perp}} x | \leq 1$ for $|x|=1$, we have
  \begin{align*}
    \sum_{|k|\leq  r} |\Pi_{k^{\perp}} x | |k|^{2} \theta_{k}^{2}
     & \geq \sum_{|k|\leq  r} |\Pi_{k^{\perp}} x |^2 |k|^{2} \theta_{k}^{2} = \sum_{|k| \leq r} |k|^2 \theta_{k}^2-\sum_{|k|\leq  r} \Big|x \cdot \frac{k}{|k|} \Big|^2 |k|^2 \theta_{k}^{2}.
  \end{align*}
Using the symmetry property of $\{\theta_{k}\}_{k}$ and \cite[Lemma 2.1]{Luo_2024_Elementary}, we obtain, for any $|x|=1$, 
  \begin{equation*}
    \sum_{|k|\leq  r} |\Pi_{k^{\perp}} x |\, |k|^{2} \theta_{k}^{2} \geq \frac{d-1}{d} \sum_{|k| \leq r} |k|^2 \theta_{k}^2 |x|^2 = \frac{d-1}{d} \sum_{|k| \leq r} |k|^2 \theta_{k}^2 .
  \end{equation*}
Thus, the non-degeneracy condition in Assumption \ref{ass: non-degenerate-L} holds with $\delta = (d-1)/d$. Using Lemma \ref{lem: sufficient ass}, Assumption \ref{ass: non-degenerate} also holds with the same $\delta = (d-1)/d$.
\end{proof}

The following proposition is the core ingredient of our energy-spectrum analysis and is used in the proofs of both the anomalous regularity and the scaling limit. It extends \cite[Lemma 6.3]{Rowan_2025_Obukhovcorrsin}, with slight differences in formulation and in the assumptions on the noise coefficients; see Corollary \ref{cor: rowan-recover}. The new weight $W_R(k)$, together with the adaptive frequency decomposition $\mathcal{T}_{M_{k,l},2M_{k,l}-1}^{k,l}$ (rather than the shell decomposition), yields a uniform discrete Poincar\'e inequality suitable for the scaling limit. By contrast, in Rowan's result the weight function tends to zero under the scaling and hence does not provide a uniform lower bound.

\begin{proposition} \label{prop: weighted Poincare}
  Let $p>1$. If noise coefficients $\{\theta_{k}\}_{k } \in \ell^2 (\Z_0^d)$ satisfy Assumption \ref{ass: non-degenerate}, then for any $R \geq 1$ and nonnegative sequence $\{ a_{k} \}_{k} \in \ell^p( \Z_0^d)$ verifying $\sum_{k} W_{R}(k) a_{k}^{p}<+\infty$, it holds
  \begin{equation} \label{eq: weighted Poincare}
    \sum_{k\in\mathbb{Z}_0^d}W_{R}(k)a_k^p  \leq
    \frac{p^2H_\delta^2}{(p-1)\delta^2}  \sum_{k,l \in \Z_0^d}  \theta_{l}^2 |\Pi_{l^\perp} k|^2  \big(a_{k+l}^{p-1}- a_{k}^{p-1} \big)(a_{k+l}-a_{k}),
  \end{equation}
  where $\delta$ is given in Assumption \ref{ass: non-degenerate} and the constant $H_{\delta} \lesssim \log (1+\frac{1}{\delta})$ is independent of $R$.
\end{proposition}

\begin{proof}
  Estimate \eqref{eq: weighted Poincare} naturally holds in the following case:
  \begin{equation*}
    \sum_{k,l \in \Z_0^d}  \theta_{l}^2 |\Pi_{l^\perp} k|^2 \big(a_{k+l}^{p-1}- a_{k}^{p-1} \big)(a_{k+l}-a_{k}) = +\infty.
  \end{equation*}
  Therefore, in the subsequent proof we will assume that the series is convergent. For simplicity, we denote $b_{k}:=a_{k}^{p}$ for all $k \in \Z_0^d$ and
  \begin{equation}\label{def-B_R}
    B_{R} :=  \sum_{k\in\Z_0^d} W_{R}(k) b_k<+\infty .
  \end{equation}

  \textbf{Step 1: estimate of $b_{k}$.} For any $k, l \in \Z_0^d$ and $c \in \N$, one has
  \[b_k = b_{k+c l} - \sum_{m=0}^{c-1} (b_{k+m l + l} - b_{k+m l} ).\]
  Define $\epsilon_{k,l}=1$ if $k \cdot l\ge 0$, and $\epsilon_{k,l}=-1$ otherwise. Then for every $m\ge 0$, $|k + m\epsilon_{k,l} l| \ge |k|$. For any $M \in \Z_{+}$, applying the triangle inequality and sum over $ M \leq c \leq 2M-1 $, we obtain
  \begin{align}
    M \, b_k & \leq \sum_{c=M}^{2M-1} \Big(b_{k+c \epsilon_{k,l} l} +  \sum_{m=0}^{c-1}| b_{k+(m+1)\epsilon_{k,l} l} - b_{k+m \epsilon_{k,l} l}|\Big) \nonumber      \\
             & \leq \sum_{c=M}^{2M-1} b_{k+c \epsilon_{k,l} l} +  M  \sum_{m=0}^{2M-2} | b_{k+(m+1)\epsilon_{k,l} l} - b_{k+m \epsilon_{k,l} l}|.  \label{eq: sum_b}
  \end{align}
  For any $a,b \in \N$ with $a\leq b$, define the segment set
  \[\mathcal{T}_{a,b}^{k,l}:= \big\{q \in \Z_0^d; q=k+m \epsilon_{k,l} l, \, m \in [a,b] \big\}.\]
  With $M=M_{k,l}$ to be determined later, dividing both sides of \eqref{eq: sum_b} by $M_{k,l}$ gives
  \[b_k \leq \frac{1}{M_{k,l}} \sum_{q \in \mathcal{T}_{M_{k,l},2M_{k,l}-1}^{k,l} } b_{q} + \sum_{q \in \mathcal{T}_{0,2M_{k,l}-2}^{k,l} } | b_{q+ \epsilon_{k,l} l} - b_{q}|.\]
  Multiplying both sides by $ \frac{|\Pi_{l^\perp} k|}{|k|} \mu_{R}(k,l)^2 \theta_{l}^2 $ and summing over $k,l \in \Z_0^d$, we obtain
  \begin{align*}
    \I & := \sum_{k\in \Z_0^d}  \sum_{l\in\mathbb{Z}_0^d} \frac{|\Pi_{l^\perp}k|}{|k|} \mu_{R}(k,l)^2 \theta_l^2 b_{k} \\
       & \, \leq \sum_{k,l} \frac{|\Pi_{l^\perp}k|}{|k|}  \frac{\mu_{R}(k,l)^2  \theta_l^2 }{M_{k,l}} \sum_{q \in \mathcal{T}_{M_{k,l},2M_{k,l}-1}^{k,l}} b_{q} \\
       & \quad + \sum_{k,l} \frac{|\Pi_{l^\perp}k|}{|k|} \mu_{R}(k,l)^2  \theta_l^2  \sum_{q \in \mathcal{T}_{0, 2M_{k,l} -2}^{k,l}} | b_{q+ \epsilon_{k,l} l} - b_{q}| \\
       & =: \II +\III.
  \end{align*}
  It is clear that $\I \le \sum_{k,l} \mu_{R}(k,l)^2 \theta_l^2 b_{k} = \sum_k W_R(k) b_k = B_R <\infty$. Assumption \ref{ass: non-degenerate} gives us $\delta B_{R} \leq \I$; hence, it suffices to estimate $\II$ and $\III$ below.

  \textbf{Step 2: estimate of $\II$.} For every $m\ge 0$, we have $|k + m\epsilon_{k,l} l| \ge |k|$, thus
  \begin{equation*}
    \mu_{R}(k,l)^2   =  \min\Big\{\frac{|l|^2}{R^2},|k|^2\Big\} \leq  \min\Big\{\frac{|l|^2}{R^2},|q|^2\Big\} =\mu_{R}(q,l)^2, \quad \forall \, q \in \mathcal{T}_{M_{k,l},2M_{k,l}-1}^{k,l}.
  \end{equation*}
  Using $\frac{|\Pi_{l^\perp}k|}{|k|} \leq 1$ and Fubini theorem, we have
  \begin{equation}\label{estim-II}
    \II \leq \sum_{k,l} \frac{\mu_{R}(k,l)^2  \theta_l^2 }{M_{k,l}} \sum_{q \in \mathcal{T}_{M_{k,l},2M_{k,l}-1}^{k,l}} b_{q} \leq \sum_{q,l} \mu_{R}(q,l)^2 \theta_l^2 b_q \sum_{k \in \Z_0^d} \ind_{ \big\{ q\in \mathcal{T}_{M_{k,l}, 2M_{k,l}-1}^{k,l} \big\} } \frac{1}{M_{k,l}}.
  \end{equation}

  Let $K \geq 1$ be a constant to be determined later; we choose $M_{k,l}$ as follows
  \begin{equation*}
    M_{k,l}:= \left\lceil K \max \Big\{1 , \frac{|k|}{|l|}  \Big\} \right\rceil.
  \end{equation*}
  For any $ q \in \mathcal{T}_{M_{k,l}, 2M_{k,l}-1}^{k,l} $, there is $ M_{k,l} \leq m \leq 2 M_{k,l}- 1$ such that $q=k+m \epsilon_{k,l} l$. Then
  \begin{align}
    |q| & \leq |k| + m |l| \leq |k| + (2 M_{k,l}-1) |l| \leq  (2 K+2) \max \{|k|,|l| \},  \label{estim-q-1} \\
    |q| & \geq \frac{|q \cdot l|}{|l|} = \frac{|(k+m \epsilon_{k,l} l ) \cdot l|}{|l|} \geq m |l| \geq K \max \{|k|,|l| \}. \nonumber
  \end{align}
  Consequently, $ \frac{|q|}{2K+2} \leq \max \{|k|,|l|  \} \leq \frac{|q|}{K}$, and we have
  \begin{equation*}
    \big\{k\in \Z_0^d; \,  q \in \mathcal{T}_{M_{k,l}, 2M_{k,l}-1}^{k,l} \big\} \subset \big\{k=q+ml; m\in\mathbb Z \big\}  \cap \Big\{k\in \Z_0^d; \frac{|q|}{2K+2} \leq \max \{|k|,|l|  \} \leq \frac{|q|}{K} \Big\}.
  \end{equation*}
  Applying Lemma \ref{lem:line-congestion} with $\gamma=1$, we obtain
  \begin{align*}
   \sum_{k \in \Z_0^d} \ind_{ \big\{ q\in \mathcal{T}_{M_{k,l}, 2M_{k,l}-1}^{k,l} \big\} } \frac{1}{M_{k,l}} \leq  \sum_{\substack{ k \in \{q+ml; m\in\mathbb Z \} \\  \frac{|q|}{2K+2} \le \max\{|k|,|l|\}\le \frac{|q|}{K}}} \frac{1}{K} \min \Big\{1,\frac{|l|}{|k|} \Big\} \leq \frac{1}{K} \Big[ 9+ 2\log \Big(2+ \frac{2}{K} \Big) \Big].
  \end{align*}
  Taking $K \geq \max\{ \frac{2}{e-2}, \frac{22}{\delta} \}$, one has 
  $$ \frac{1}{K} \Big[ 9+ 2\log \Big(2+ \frac{2}{K} \Big) \Big] \le \frac\delta{22} \big[ 9+ 2\log (2+ e-2 ) \big] = \frac\delta 2. $$
  Substituting these estimates into \eqref{estim-II} and using Assumption \ref{ass: non-degenerate}, we have
  \begin{equation*}
    \II \leq \frac{\delta}{2}  \sum_{q,l} \mu_{R}(q,l)^2 \theta_l^2 b_q \leq \frac{1}{2} \I.
  \end{equation*}

  \textbf{Step 3: estimate of $\III$.} Using Lemma \ref{lem: bound p}, we have
  \begin{equation*}
    | b_{q+ \epsilon_{k,l} l} - b_{q}|= \big|a_{q+ \epsilon_{k,l} l}^p - a_{q}^p \big| \leq \frac{p}{\sqrt{p-1}} \sqrt{a_{q+ \epsilon_{k,l} l}^p + a_{q}^p} \sqrt{\big(a_{q+ \epsilon_{k,l} l}^{p-1}- a_{q}^{p-1} \big) (a_{q+ \epsilon_{k,l} l}-a_{q})}.
  \end{equation*}
  Applying the Cauchy--Schwarz inequality to $\III$, we obtain
  \begin{align*}
    \III & = \sum_{k,l} \frac{|\Pi_{l^\perp}k|}{|k|} \mu_{R}(k,l)^2  \theta_l^2 \sum_{q\in \mathcal{T}_{0,2M_{k,l}-2}^{k,l}} | b_{q+ \epsilon_{k,l} l} - b_{q}| \\
    & \leq \bigg( \sum_{k,l} \mu_{R}(k,l)^2 \theta_l^2 \, \frac{\mu_{R}(k,l)}{|k|} \sum_{q \in \mathcal{T}_{0,2M_{k,l}-2}^{k,l}} (b_{q+ \epsilon_{k,l} l}+ b_{q}) \bigg)^{1/2} \\
    & \quad \times \bigg( \frac{p^2}{p-1} \sum_{k,l} |\Pi_{l^\perp}k|^2 \theta_l^2 \, \frac{\mu_{R}(k,l)}{|k|} \sum_{q \in \mathcal{T}_{0,2M_{k,l}-2}^{k,l}} \big(a_{q+ \epsilon_{k,l} l}^{p-1}- a_{q}^{p-1}\big) (a_{q+ \epsilon_{k,l} l}-a_{q}) \bigg)^{1/2} .
  \end{align*}
  For any $ q \in \mathcal{T}_{0,2M_{k,l}-1}^{k,l}$, we have $\mu_{R}(k,l) \leq \mu_{R}(q,l)$ and $|\Pi_{l^\perp}k|^2 =|\Pi_{l^\perp}q |^2$. Then, using the fact $\theta_{l}=\theta_{-l}$ and Fubini's theorem, we get $\III \leq S_1^{1/2} S_2^{1/2}$, where
  \begin{align*}
    S_1 & := 2 \sum_{q,l} \mu_{R}(q,l)^2  \theta_l^2  b_{q} \sum_{k\in \Z_0^d} \ind_{ \big\{q \in \mathcal{T}_{0, 2M_{k,l} -1}^{k,l} \big\} } \frac{\mu_{R}(k,l)}{|k|}, \\
    S_2 & : = 2 \frac{p^2}{p-1} \sum_{q,l} |\Pi_{l^\perp}q|^2 \theta_l^2 \big(a_{q+ l}^{p-1}- a_{q}^{p-1}\big) (a_{q+ l}-a_{q}) \sum_{k\in \Z_0^d} \ind_{ \big\{q \in \mathcal{T}_{0, 2M_{k,l} -1}^{k,l} \big\} } \frac{\mu_{R}(k,l)}{|k|}.
  \end{align*}
  For any $q \in \mathcal{T}_{0, 2M_{k,l} -1}^{k,l} $, there is $ 0 \leq m \leq 2 M_{k,l}- 1$ such that $q=k+m \epsilon_{k,l} l$ and the estimate \eqref{estim-q-1} holds; moreover, by the definition of $\epsilon_{k,l}$, we have $|q| \geq |k|$. Thus,
  \begin{equation*}
    \big\{k\in \Z_0^d; \,  q \in \mathcal{T}_{0, 2M_{k,l} -1}^{k,l} \big\} \subset \big\{k=q+ml; m\in\mathbb Z \big\}  \cap \Big\{k; \;  \frac{|q|}{2K+2} \leq \max \{|k|,|l|  \}  , \,  |k| \leq |q|\Big\}.
  \end{equation*}
  Due to $\frac{\mu_{R}(k,l)}{|k|} = \min\big\{1,\frac{|l|}{R |k|} \big\} $, applying Lemma \ref{lem:line-congestion} with $\gamma=R$, we get
  \begin{align*}
    \sum_{k\in \Z_0^d} \ind_{ \big\{q \in \mathcal{T}_{0, 2M_{k,l} -1}^{k,l} \big\} } \frac{\mu_{R}(k,l)}{|k|} \leq  \sum_{\substack{ k \in \{q+ml; m\in\mathbb Z \} \\  \frac{|q|}{2K+2} \le \max\{|k|,|l| \} , \; |k|\le |q|}}  \min \Big\{1,\frac{|l|}{R |k|} \Big\} \leq 9+ \frac{2}{R}\log (2 K+2) .
  \end{align*}
  Substituting this estimate into $S_1$ and $S_2$, we get
  \begin{equation*}
    \III \leq \Big(18+ \frac{4}{R}\log (K+2) \Big) B_{R}^{1/2} \Big( \frac{p^2}{p-1} \sum_{q,l} |\Pi_{l^\perp}q|^2 \theta_l^2 \big(a_{q+ l}^{p-1}- a_{q}^{p-1} \big) (a_{q+ l}-a_{q})  \Big)^{1/2},
  \end{equation*}
  where $B_R$ is defined in \eqref{def-B_R}.

  Combining the above estimates on $\I$, $\II$, $\III$ and Assumption \ref{ass: non-degenerate}, we get
  \begin{equation*}
    \frac{\delta}{2} B_{R} \leq \frac{1}{2} \I \leq \Big(18+ \frac{4}{R}\log (K+2) \Big) B_{R}^{1/2} \Big( \frac{p^2}{p-1} \sum_{q,l} |\Pi_{l^\perp}q|^2 \theta_l^2 \big(a_{q+ l}^{p-1}- a_{q}^{p-1} \big) (a_{q+ l}-a_{q}) \Big)^{1/2}.
  \end{equation*}
  Due to $B_{R} <+\infty$, we take $K = \max\{ \frac{2}{e-2}, \frac{22}{\delta} \}$ and get
  \begin{equation*}
    B_{R} \leq \frac{H_{\delta}^2}{\delta^2} \, \frac{p^2}{p-1} \sum_{q,l} |\Pi_{l^\perp}q|^2 \theta_l^2 \big(a_{q+ l}^{p-1}- a_{q}^{p-1} \big) (a_{q+ l}-a_{q}),
  \end{equation*}
  where $H_{\delta}= \sup_{R \geq 1}\big[ 36+ \frac{8}{R}\log \big( \max\{ \frac{2}{e-2}, \frac{22}{\delta} \}+2 \big) \big]^2 \lesssim \log (1+\frac{1}{\delta})$.
\end{proof}

The following corollary recovers the result of Rowan \cite[Lemma 6.3]{Rowan_2025_Obukhovcorrsin}. Moreover, thanks to our finer estimates, we remove conditions (1.24) and (1.26) from \cite[Assumption 1.15]{Rowan_2025_Obukhovcorrsin}.

\begin{corollary} \label{cor: rowan-recover}
  Let $p>1$ and $\alpha \in (0,1)$. Assume that the noise coefficients $\{\theta_{k}\}_k \in \ell^2 (\Z_0^d)$ satisfy $\sum_{k} |k|^{2\alpha} \theta_{k}^2<+\infty$, and there exists some $\delta_0>0$ such that the weight function $S(r):=\sum_{|k| \leq r} |k|^{1+\alpha} \theta_{k}^2$ verifies
  \begin{equation} \label{ass: non-degenerate-S}
    \inf_{x \in \R^d, |x|=1 } \sum_{|k|\leq  r} |\Pi_{k^{\perp}} x | \,  |k|^{1+\alpha} \theta_{k}^{2} \geq \delta_0 S(r) \quad (\forall \, r \geq 1),
  \end{equation}
  then there exists a constant $C>0$ such that for any $R \geq 1$ and any nonnegative sequence $\{ a_{k} \}_{k} $ satisfying $\sum_{k} |k|^{2(1-\alpha)} a_{k}^{p}<+\infty$, it holds
  \begin{equation} \label{eq: weighted Poincare-S}
    \sum_{k\in\mathbb{Z}_0^d} S(R |k|)^2 a_k^p  \leq
    \frac{p^2 C R^2}{(p-1)\delta_0^4}  \sum_{k,l \in \Z_0^d}  \theta_{l}^2 |\Pi_{l^\perp} k|^2  \big(a_{k+l}^{p-1}- a_{k}^{p-1} \big)(a_{k+l}-a_{k}) .
  \end{equation}
\end{corollary}

\begin{proof}
  We define new noise coefficients $\{ \widetilde{\theta}_{k} \}_{k}$ as
  \begin{equation*}
    \widetilde{\theta}_{k}^2 := S(|k|) |k|^{\alpha-1} \theta_{k}^2, \quad \forall \, k \in \Z_0^d;
  \end{equation*}
  the corresponding function $L$ is given as
  \begin{equation*}
    L(r):=\sum_{|k| \leq r} |k|^2 \widetilde{\theta}_{k}^2=\sum_{|k| \leq r} S(|k|) |k|^{1+\alpha} \theta_{k}^2, \quad \forall\, r \geq 1.
  \end{equation*}
  Using the Fubini theorem, we obtain
  \begin{align*}
    L(r) 
    & = \sum_{|k| \leq r} |k|^{1+\alpha} \theta_{k}^2 \sum_{|l| \leq |k|}  |l|^{1+\alpha} \theta_{l}^2 \\
    & =\frac{1}{2} \Big( \sum_{|k| \leq r} |k|^{1+\alpha} \theta_{k}^2 \sum_{|l| \leq |k|}  |l|^{1+\alpha} \theta_{l}^2 + \sum_{|l| \leq r} |l|^{1+\alpha} \theta_{l}^2  \sum_{|l| \leq |k| \leq r}  |k|^{1+\alpha} \theta_{k}^2 \Big)   \\
    & \geq \frac{1}{2} \sum_{|k| \leq r} |k|^{1+\alpha} \theta_{k}^2 \sum_{|l| \leq r}  |l|^{1+\alpha} \theta_{l}^2 =\frac{1}{2} S(r)^2.
  \end{align*}
  By a similar argument as above, we get, for any unit vector $x \in \R^d$,
  \begin{align*}
    \sum_{|k|\leq  r} |\Pi_{k^{\perp}} x | \, |k|^2 \widetilde{\theta}_{k}^2 
    & = \sum_{|k|\leq  r} |\Pi_{k^{\perp}} x | \, S(|k|) |k|^{1+\alpha} \theta_{k}^2  \\
    & \geq \sum_{|k|\leq  r} |\Pi_{k^{\perp}} x | \, |k|^{1+\alpha} \theta_{k}^2 \sum_{|l| \leq |k|} |\Pi_{l^{\perp}} x | \, |l|^{1+\alpha} \theta_{l}^2 \\
    & \geq \frac{1}{2} \bigg( \sum_{|k|\leq r} |\Pi_{k^{\perp}} x | \, |k|^{1+\alpha} \theta_{k}^2 \bigg)^2 .
  \end{align*}
  Using the condition \eqref{ass: non-degenerate-S} and the simple fact $L(r) \leq S(r)^2$, we have
  \begin{equation*}
    \inf_{x \in \R^d, |x|=1 } \sum_{|k|\leq  r} |\Pi_{k^{\perp}} x | \,  |k|^2 \widetilde{\theta}_{k}^2 \geq \frac{1}{2} \Big( \inf_{x \in \R^d, |x|=1 } \sum_{|k|\leq  r} |\Pi_{k^{\perp}} x | \,   |k|^{1+\alpha} \theta_{k}^2  \Big)^2  \geq \frac{\delta_0^2}{2} S(r)^2 \geq \frac{\delta_0^2}{2} L(r).
  \end{equation*}
  Thus, Lemma \ref{lem: sufficient ass} implies that $\{ \widetilde{\theta}_{k} \}_{k}$ satisfies Assumption \ref{ass: non-degenerate} with $\delta=\frac{\delta_0^2}{2}$.

  Next we will show that $\sum_{k} \widetilde{W}_{R}(k) a_{k}^{p}<+\infty$, where
  \begin{equation*}
    \widetilde{W}_{R}(k) = \frac{1}{R^2} L(R |k|) + \bigg( \sum_{|l|> R |k|} \widetilde{\theta}_{l}^2 \bigg) |k|^2, \quad \forall \, k \in \Z_0^d.
  \end{equation*}
  Setting $M:=\sum_{k} |k|^{2\alpha} \theta_{k}^2<+\infty $, then we have
  \begin{equation*}
    S(r):=\sum_{|k| \leq r} |k|^{1+\alpha} \theta_{k}^2 = \sum_{|k| \leq r} |k|^{1-\alpha} |k|^{2\alpha} \theta_{k}^2 \leq M r^{1-\alpha}, \quad  \forall \,  r \geq 1.
  \end{equation*}
  Recalling the definition of $\widetilde{\theta}_{k}^2$, this implies $\widetilde{\theta}_{k}^2 \leq M \theta_{k}^2$ for any $k\in \Z_0^d$, and
  \begin{align*}
    \widetilde{W}_{R}(k) 
    & \leq  \frac{1}{R^2} S(R |k|)^2+ M \bigg( \sum_{|l|> R |k|} |l|^{-2\alpha} |l|^{2\alpha} \theta_{l}^2 \bigg) |k|^2 \\
    & \leq \frac{1}{R^2} M^2 (R |k|)^{2-2\alpha} + \frac{M^2}{ (R|k|)^{2\alpha} }  |k|^2 = \frac{2M^2}{R^{2\alpha}} |k|^{2-2\alpha}.
  \end{align*}
  Since $\sum_{k} |k|^{2(1-\alpha)} a_{k}^{p}<+\infty$, we get $\sum_{k} \widetilde{W}_{R}(k) a_{k}^{p}<+\infty$.

  Thus, we can apply Proposition \ref{prop: weighted Poincare} with the new noise coefficients $\{ \widetilde{\theta}_{k} \}_{k}$:
  \begin{equation*}
    \sum_{k\in\mathbb{Z}_0^d} \widetilde{W}_{R}(k)a_k^p  \leq
    \frac{p^2H_{\delta}^2}{(p-1)\delta^2}  \sum_{k,l \in \Z_0^d}  \widetilde{\theta}_{l}^2 |\Pi_{l^\perp} k|^2  \big(a_{k+l}^{p-1}- a_{k}^{p-1} \big)(a_{k+l}-a_{k})
  \end{equation*}
  with $\delta=\frac{\delta_0^2}{2}$. Using $\widetilde{W}_{R}(k) \geq \frac{1}{R^2} L(R |k|) \geq \frac{1}{2 R^2} S(R |k|)^2$ and $\widetilde{\theta}_{k}^2 \leq M \theta_{k}^2$, we get \eqref{eq: weighted Poincare-S}.
\end{proof}

\section{Summataion reordering and cutoff function} \label{sec: summatation reordering}

Let $\chi\in C_c^\infty(\mathbb R_+;[0,1])$ be a decreasing cutoff function, with $\chi(r)=1$ if $r\in [0,1]$ and $\chi(r)=0$ if $r \geq 2$; for any $N\ge 1$, we define $\chi_N(|k|):= \chi(|k|/N)$ for all $k \in \Z_0^d$.
In this section, we justify the interchange of the summation order involving the cutoff functions $\chi_N$. To this end, we first establish a uniform estimate for the cutoff functions.

\begin{lemma} \label{lem-cutoff}
  There exists $C_{\chi} > 0$ such that the following estimate holds uniformly in $N\geq1$:
  \begin{equation} \label{claim: cutoff}
    \big|\chi_N(|k+l|)+\chi_N(|k-l|)-2\chi_N(|k|)\big|
    \leq
    C_\chi\min\Big\{1,\frac{|l|^2}{|k|^2}\Big\}, \quad \forall\, k,l\in \Z_0^d
  \end{equation}
\end{lemma}
\begin{proof}

  Due to $\chi$ is smooth, there exists $L_{\chi} > 0$ such that $\|D^2 \chi_{N}\|_{\infty} \leq L_{\chi} N^{-2}$. Using the fact $0\le \chi \leq 1$, the second order Taylor expansion of $\chi_N$ gives us
  \begin{equation*}
    \big|\chi_N(|k+l|)+\chi_N(|k-l|)-2\chi_N(|k|)\big| \leq \min\big\{2, L_{\chi} N^{-2} |l|^2 \big\} \leq \min\Big\{2, 16 L_{\chi} \frac{|l|^2}{|k|^2} \Big\} 
  \end{equation*}
  for all $|k|\leq 4N$. When $|k|>4N$ and $|k|>2|l|$, then $\chi_N(|k|)=0$ and
  \begin{equation*}
    \min\{|k+l|,|k-l| \} \geq |k|-|l| \geq \frac{|k|}{2} > 2N.
  \end{equation*}
  Thus, $\chi_N(|k+l|)=\chi_N(|k-l|)=0$ and the left-hand side of \eqref{claim: cutoff} is zero. If $|k| \leq 2 |l|$, then $\frac{|l|}{|k|} \geq \frac{1}{2}$; we have
  \begin{equation*}
    \big|\chi_N(|k+l|)+\chi_N(|k-l|)-2\chi_N(|k|)\big| \leq 2 \leq \min\Big\{2, 8 \frac{|l|^2}{|k|^2} \Big\}.
  \end{equation*} 
  Combining the above estimates, we get the dsired result \eqref{claim: cutoff} with $C_\chi = 16 \max\{ 1, L_{\chi}\}$.
\end{proof}

The following lemma shows that, under suitable conditions on the coefficients $\{\theta_k \}_k$, we can exchange the order of summation involving the cutoff function $\chi_N$ and take the limit $N \rightarrow +\infty$. This result plays an important role in the proofs of Propositions \ref{prop-Y_k^p-diff} and \ref{prop-optimal}.

\begin{lemma} \label{lem-sum-cutoff}
  For $R \geq 1$, let $W_{R}(k)$ be the weight function defined in Section \ref{sec: poincare inequality}. Assume that
  \begin{itemize}
    \item the coefficients $\{\theta_k \}_k \in \ell^2(\Z_0^d)$ satisfy $\theta_{k}=\theta_{-k}$ for all $k \in \Z_0^d$;
    \item $\{X_{k} \}_{k \in \Z_0^d}$ is a sequence of nonnegative functions on $[0,T]$ satisfying, for some $p \geq 1$,
          \begin{equation*}
            \int_{0}^{T} \sum_{k} W_{R}(k) X_{k}^p (t) \, \d t< +\infty.
          \end{equation*}
  \end{itemize}
  Denote the summation with cutoff function $\chi_N$ as
  \begin{equation*}
    S_N := 2 \int_{0}^{T} \sum_k \sum_l \theta_l^2 |\Pi_{l^\perp} k|^2 \chi_{N}(|k|) X_k^{p-1} (X_{k-l} -X_k) \, \d  s,
  \end{equation*}
  then it holds
  \begin{equation*}
    \lim_{N \rightarrow \infty} S_N = -\int_{0}^{T} \sum_{k,l} \theta_l^2 |\Pi_{l^\perp} k|^2 \big(X_{k+l}^{p-1} -X_k^{p-1}\big) (X_{k+l} -X_k) \, \d s.
  \end{equation*}
  In particular, if $p=1$, then we have $\lim_{N\to \infty} S_N=0$.
\end{lemma}

\begin{proof}
  Due to the cutoff function $\chi_{N}$ and $\{\theta_{k}\}_k \in \ell^2(\Z_0^d)$, for any fixed $N$ the series $S_N$ converges absolutely. By Fubini's theorem, we have
  \begin{equation}\label{eq: Fourier sums-cuttoff}
    S_N = 2 \int_{0}^{T} \sum_l \sum_k  \theta_l^2 |\Pi_{l^\perp} k|^2 \chi_{N}(|k|) X_k^{p-1} (X_{k-l} -X_k) \, \d  s.
  \end{equation}
  Analogous to the proof of Proposition \ref{prop-Y_k^p}, we will use the facts that $\theta_{k}=\theta_{-k}$ and $\Pi_{l^\perp} l=0$. For any fixed $l \in \Z_0^d$, substituting $k$ for $\tilde{k}+l$ on the right-hand side of the above equality yields
  \begin{equation*}
    S_N = 2 \int_{0}^{T} \sum_l \sum_{\tilde k}  \theta_l^2 |\Pi_{l^\perp} \tilde k|^2 \chi_{N}(|\tilde{k}+l|)  X_{\tilde{k}+l}^{p-1} \big( X_{\tilde k} -X_{\tilde{k}+l} \big) \, \d  s.
  \end{equation*}
  Substituting $l$ for $-\tilde{l}$ in the identity \eqref{eq: Fourier sums-cuttoff} and using $\theta_{l}=\theta_{-l}$, we get
  \begin{equation*}
    S_N = 2 \int_{0}^{T} \sum_{\tilde{l}} \sum_k  \theta_{\tilde{l}}^2 |\Pi_{{\tilde{l}}^\perp} k|^2 \chi_{N}(|k|) X_k^{p-1} (X_{k+\tilde{l}} -X_k) \, \d  s.
  \end{equation*}
  Changing $\tilde k, \tilde l$ back to $k,l$ in the two equalities above and summing them up, we obtain
  \begin{equation*}
    S_N = - \int_{0}^{T} \sum_l \sum_{ k}  \theta_l^2 |\Pi_{l^\perp}  k|^2 \big( \chi_{N}(|k+l|) X_{k+l}^{p-1}  -\chi_{N}(|k|)  X_{k}^{p-1} \big) \big( X_{k+l}-X_{ k}  \big) \, \d  s.
  \end{equation*}

  Define a nonnegative function $\Phi_{p}: \R_+ \times \R_+ \rightarrow  \R_+$ as
  \begin{equation*}
    \Phi_{p} (x,y):= \begin{cases} 
    \frac{x^p-y^p}{p}- y^{p-1}(x-y), & p > 1; \\ 
    0, &  p=1.
    \end{cases}
  \end{equation*}
  Then, we have $S_N=D_N+E_{N,1} + E_{N,2}$, where
  $$\aligned
    D_N & := -  \int_{0}^{T} \sum_l \sum_{ k}  \theta_l^2 |\Pi_{l^\perp}  k|^2 \Big( \chi_{N}(|k+l|) \Phi_{p}(X_{k},X_{k+l})  + \chi_{N}(|k|) \Phi_{p}(X_{k+l},X_{k}) \Big)  \, \d  s, \\
    E_{N,1} &:=- \frac1p \int_{0}^{T} \sum_l \sum_{ k}  \theta_l^2 |\Pi_{l^\perp}  k|^2 \big( \chi_{N}(|k+l|) - \chi_{N}(|k|) \big)  X_{k+l}^{p}  \, \d  s, \\
    E_{N,2} &:= \frac1p\int_{0}^{T} \sum_l \sum_{ k}  \theta_l^2 |\Pi_{l^\perp}  k|^2 \big( \chi_{N}(|k+l|) - \chi_{N}(|k|) \big) X_{k}^{p}  \, \d  s.
    \endaligned $$
  Note that $D_N=0$ in the case $p=1$. For $p>1$, every summand in $D_N$ is nonnegative and
  \begin{equation*}
    \Phi_{p}(X_{k},X_{k+l})  + \Phi_{p}(X_{k+l},X_{k}) = \big(X_{k+l}^{p-1} -X_k^{p-1}\big) (X_{k+l} -X_k).
  \end{equation*}
  Thus, the monotone convergence theorem gives
  \begin{equation*}
    \lim_{N \rightarrow +\infty} D_{N} = -\int_{0}^{T} \sum_{l} \sum_{k} \theta_l^2 |\Pi_{l^\perp} k|^2 \big(X_{k+l}^{p-1} -X_k^{p-1}\big) (X_{k+l} -X_k) \, \d s,
  \end{equation*}
 which vanishes for $p=1$.

  For any fixed $l \in \Z_0^d$, substituting $k$ for $\tilde{k}-l$ in $E_{N,1}$ yields
  \begin{equation*}
    E_{N,1} =- \frac1p \int_{0}^{T} \sum_l \sum_{\tilde k}  \theta_l^2 |\Pi_{l^\perp}  \tilde{k}|^2 \big( \chi_{N}(|\tilde{k}|) - \chi_{N}(|\tilde{k}-l|) \big)  X_{\tilde{k}}^{p}  \, \d s.
  \end{equation*}
  Then $E_N:=E_{N,1}+ E_{N,2}$ is
  \begin{equation*}
    E_{N} = \frac1p \int_{0}^{T} \sum_l \sum_{ k}  \theta_l^2 |\Pi_{l^\perp}  k|^2 \big( \chi_{N}(|k+l|) + \chi_{N}(|k-l|) -2 \chi_{N}(|k|) \big) X_{k}^{p}  \, \d s.
  \end{equation*}
  Next, we will prove $E_N \rightarrow 0$. Note that $|\Pi_{l^\perp}k|\leq|k|$, by estimate \eqref{claim: cutoff}, we obtain
  \begin{align*}
     & \sum_l\theta_l^2|\Pi_{l^\perp}k|^2 \big|\chi_{N}(|k+l|) + \chi_{N}(|k-l|) -2 \chi_{N}(|k|) \big| \\
     & \quad \leq C_\chi \sum_l\theta_l^2|k|^2  \min\Big\{1,\frac{|l|^2}{|k|^2}\Big\}                      \\
     & \quad \leq C_\chi R^2 \sum_l \mu_{R}(k,l)^2 \theta_l^2 = C_\chi R^2 W_{R}(k),
  \end{align*}
  for any $R \geq 1$, where $ \mu_{R}(k,l)$ is given in Section \ref{sec: poincare inequality}. Suppose we have shown that for any fixed $k\in \Z_0^d$, 
  \begin{equation}\label{appen-DCT}
  \lim_{N\to \infty} \sum_l\theta_l^2|\Pi_{l^\perp}k|^2 \big|\chi_{N}(|k+l|) + \chi_{N}(|k-l|) -2 \chi_{N}(|k|) \big| =0, 
  \end{equation}
  then, due to $\int_{0}^{T} \sum_{ k} W_{R}(k)  X_{k}^{p}   \, \d  t <+\infty $, the dominated convergence theorem implies that $E_N \rightarrow 0$ as $N \rightarrow \infty$. Since $S_N=D_N+E_N$, the above limit results give the desired assertion.

  It remains to show the limit \eqref{appen-DCT} for any fixed $k\in \Z_0^d$. As a function of $l\in \Z_0^d$, the summand is bounded by $2 \theta_l^2 |k|^2$, which is summable due to $\{\theta_l \}_l \in \ell^2(\Z_0^d)$; moreover, for each $l\in \Z_0^d$, it is clear that $\chi_{N}(|k+l|) + \chi_{N}(|k-l|) -2 \chi_{N}(|k|) \to 0$ as $N\to \infty$. Hence, we can use again the dominated convergence theorem to conclude \eqref{appen-DCT}. This finishes the proof.
\end{proof}

\bigskip

\noindent \textbf{Acknowledgements.} The first author is grateful to the financial supports of the National Natural Science Foundation of China (No. 12595281) and the Open Foundation of the State Key Laboratory of Mathematical Sciences (Grant No. 09). The second author is grateful to the financial supports of the National Natural Science Foundation of China (No. 12601266) and the Stability Support Program for Basic Research (Grant No. WZJC202601005).


\end{document}